\documentclass[10pt]{article}
\usepackage{amssymb,amsmath,amsthm}
\usepackage{bm}
\usepackage{}
\usepackage{amsbsy}
\usepackage{psfrag}
\usepackage{graphics}
\usepackage{graphicx}
\usepackage{pgfplots}
\usepackage{relsize}
\usepackage{color}
\usepackage{xcolor}
\usepackage{mathtools}
\usepackage{pdfpages}
\usepackage{url}
\usepackage{makecell}
\usepackage[nohead,margin=1.0in]{geometry}
 \usepackage{hyperref}
\usepackage{graphicx}
\usepackage{caption}
\usepackage{subcaption}
\usepackage{verbatim}
\usepackage{booktabs}
\usepackage{threeparttable}
\hypersetup{
	colorlinks=true,
	linkcolor=blue,
	filecolor=magenta,
	urlcolor=cyan,
}

\theoremstyle{plain}
\theoremstyle{theorem}
\newtheorem{theorem}{Theorem}[section]
\newtheorem{proposition}{Proposition}[section]
\newtheorem{definition}{Definition}[section]

\newtheorem{lemma}{Lemma}[section]

\newtheorem{example}[theorem]{Example}
\newcommand{\norm}[1]{\left\Vert #1 \right\Vert}
\newcommand{\brb}[1]{\Bigl( #1 \Bigr)}
\newcommand{\abs}[1]{\left| #1 \right|}
\newcommand{\bs}{\boldsymbol}
\newcommand{\nd}{n_d}
\newcommand{\gd}{g_d}
\newtheorem{remark}[theorem]{Remark}

\begin{document}
\allowdisplaybreaks
\numberwithin{equation}{section}

\title{Solving Elliptic Optimal Control Problems via Neural Networks and Optimality System\thanks{The work of B. Jin is supported by Hong Kong RGC General Research Fund (Project 14306423), and a start-up fund from The Chinese University of Hong Kong. The work of  Z. Zhou is supported by Hong Kong Research Grants Council (15303021) and an internal grant of Hong Kong Polytechnic University (Project ID: P0038888, Work Programme: ZVX3).}}

\author{Yongcheng Dai\thanks{Department of Applied Mathematics,
The Hong Kong Polytechnic University, Kowloon, Hong Kong, P.R. China (\texttt{yongcheng.dai@connect.polyu.hk,zhizhou@polyu.edu.hk})}
\and Bangti Jin\thanks{Department of Mathematics, The Chinese University of Hong Kong, Shatin, New Territories, Hong Kong, P.R. China (\texttt{bangti.jin@gmail.com, b.jin@cuhk.edu.hk}).}
\and Ramesh Chandra Sau\thanks{Department of Mathematics, Indian Institute of Technology, Bombay, India (\texttt{rcsau1994@gmail.com})} 
\and Zhi Zhou\footnotemark[2]}

\maketitle

\begin{abstract}
In this work, we investigate a neural network based solver for optimal control problems (without / with box constraint) for linear and semilinear second-order elliptic problems. It utilizes a coupled system derived from the first-order optimality system of the optimal control problem, and employs deep neural networks to represent the solutions to the reduced system. We present an error analysis of the scheme, and provide $L^2(\Omega)$ error bounds on the state, control and adjoint in terms of  neural network parameters (e.g., depth, width, and parameter bounds) and the numbers of sampling points. The main tools in the analysis include offset Rademacher complexity and boundedness and Lipschitz continuity of neural network functions. We present several numerical examples to illustrate the method and compare it with two existing ones.
\end{abstract}


\section{Introduction}

In the past decade, neural networks (NNs) have demonstrated remarkable performance across a wide range of applications, e.g., computer vision, imaging, and natural language processing. These successes have permeated into scientific computing, and NNs have become popular tools to approximate solutions to PDEs, which represent one important class of mathematical models for describing physical phenomena. There have been intensive activities in developing neural solvers  for solving PDEs, e.g., deep Ritz method \cite{EYu:2018}, deep Galerkin method \cite{Sirignano:2018dgm}, weak adversarial network \cite{ZangBaoZhou:2020}, deep least-squares method \cite{CaiChenLiu:2021} and physics informed neural networks \cite{RaissiPerdolarisKarniadakis:2019}; see the reviews   \cite{Karniadakis:2021nature,BeckJentzen:2023,TanyuNing:2023} for the details of the methods. All these solvers employ NNs as the ansatz space to approximate the PDE solutions directly. These solvers have shown promising empirical results \cite{Karniadakis:2021nature,TanyuNing:2023}, and hold potentials for high-dimensional problems.

Thus, it is a very natural step forward to explore the use of NNs for solving optimal control problems subject to PDE constraint, which encompass an important class of problems in practice and has important applications, e.g., fluid flow, heat conduction, structural optimization, microelectronics, crystal growth and vascular surgery. The mathematical theory of distributed / boundary optimal control subject to elliptic / parabolic PDEs has been well developed (see, e.g., \cite{Lions:1971,Trolzstch:2005,ManzoniQuarteroniSalsa:2021}). Traditionally, one popular numerical approach utilizes the first-order optimality system, using either adjoint state or Lagrangian multiplier, and then iteratively updates the control using established optimization algorithms, e.g., gradient descent and Newton method. In practice, the resulting continuous formulation is often discretized by finite element method (FEM) and spectral methods \cite{ManzoniQuarteroniSalsa:2021}. 

We investigate the use of deep neural networks for solving high-dimensional elliptic optimal control problems. 
The method is based on the well-known Karush-Kuhn-Tucker (KKT) system \cite{Lions:1971,Trolzstch:2005}, which consists of a state equation, an adjoint equation, and an equality / variational inequality for the unconstrained and constrained cases, respectively. The (in)equality allows recovering the optimal control from the optimal adjoint, and the KKT system can be reduced a coupled system for the state and adjoint. Then one applies deep neural networks to solve the reduced optimality system, and hence the name OSNN. We describe the discrete formulations in detail for the unconstrained, constrained, and semilinear cases. Moreover, we give an error analysis of the approach for linear elliptic optimal control problems, and  derive an $L^2(\Omega)$ error estimate for the approximate control, state, and adjoint. The convergence analysis represents the main technical contribution of the work, and it is carried out in two steps. First, using the PDE regularity theory and the idea of reconstruction, we derive weak coercivity estimates which bound the state etc. using the population loss, cf. Lemmas \ref{lem:fund-est1} and \ref{lem:fund-est2}. Second, using statistical learning theory, we bound the difference between the population and empirical loss, which involves two parts: approximation and statistical errors. The former arises from approximating the solution of the problem using NNs and is bounded using approximation theory of NNs \cite{GuhrinRaslan:2021}, and the latter arises from using Monte Carlo methods to approximate the (potentially) high-dimensional integrals and is bounded using offset Rademacher complexity \cite{LiangRakhlin:2015,DuanJiao:2023fast} and suitable continuity estimates on NN functions. This analysis yields sharper bounds on the statistical error, improving several existing works (see, e.g., \cite{JiaoLai:2021error,HuJinZhou:2024}) that are based on Rademacher complexity, and thus it is of independent interest. Finally, we present several numerical experiments to illustrate features of the method, and give a comparative study with two existing neural solvers. To the best of our knowledge, OSNN represents the first NN based solver with convergence guarantees for elliptic optimal control problems.

In the literature, several neural solvers have been proposed for elliptic optimal control problems. One class of methods is based on penalty, including standard penalty method (soft constraint) \cite{MowlaviNabi:2023}, and augmented Lagrangian
method \cite{LuJohnson:2021}. The idea is to transform the PDE-constrained optimization into an unconstrained
one by penalizing the PDE residual. However, the resulting optimization problem becomes increasingly ill-conditioned
(and hence more challenging to optimize) as the penalty weight increases \cite{Bertsekas:1982,LuJohnson:2021}, but imposing the constraint accurately does require a large weight. Moreover,  the handling of multiple constraints is nontrivial. The augmented Lagrangian method can partly alleviate ill-conditioning. The second class of methods employs the adjoint technique \cite{Cea:1986}, which solves the adjoint equation and computes the derivative of the reduced objective with respect to the control. One method in the class is the direct adjoint looping (DAL) \cite{Lions:1971,Jameson:1988}, which iteratively updates the control using a gradient descent scheme with the gradient computed via the adjoint technique. Based on DAL, Yin et al \cite{YinYang:2023} proposed the adjoint-oriented neural network (AONN). It is an NN
realization of the adjoint-based gradient descent scheme, and involves solving the state and adjoint and then updates the control sequentially. The authors reported excellent empirical performance on a number of challenging settings, e.g., (parametric) Navier-Stokes problem. See Section \ref{ssec:existing} for further descriptions. Demo et al \cite{DemoStrazzullo:2023} proposed an approach that augments NN inputs with parameters so that the KKT system and neural networks are combined for parametric optimal control problems; see also \cite{BarrySarshar:2022}. Recently, Song et al \cite{SongYuan:2023} develop an NN and ADMM based approach for efficiently solving nonsmooth PDE constrained optimization. In comparison, we develop an easy-to-implement numerical scheme, thoroughly analyze its convergence, and evaluate the method against existing deep-learning based techniques. Note that there is also a rich literature on using DNNs for other types of optimal control problems; see, e.g., the references \cite{NakamuraGong:2021,MengZhangKarniadakis:2022,Borovykh:2022,ZhaoHan:2024} for an incomplete list and references therein.

The rest of the paper is organized as follows. 
In Section \ref{sec:contProbPINN}, we describe the details of OSNN for three cases. In Section \ref{sec:erranalysis}, we present an error analysis for both constrained and unconstrained cases, which relies on technical estimates in the appendix. Finally, in Section \ref{sec:numexp}, we present numerical examples to illustrate the method. Throughout $c$ denotes a generic constant which may change at each occurrence, but it is always independent of NN parameters.

\section{Solving elliptic optimal control using DNNs}\label{sec:contProbPINN}

Now we develop the numerical scheme for solving elliptic optimal control problems via the coupled optimality system, and term the resulting method as OSNN. We also survey three existing neural solvers.

\subsection{Optimal control without control constraint}\label{subsec:unconst}
Let $\Omega\subset (-1,1)^d\subset\mathbb{R}^d$ ($d\geq 1$) be a bounded domain with a smooth boundary $\partial\Omega$. Consider the following distributed optimal control problem
\begin{align}\label{eq:min-J:Diffu}
\min\; J(y,u):=\tfrac{1}{2}\|y-y_d\|_{L^2(\Omega)}^2+\tfrac{\lambda}{2}\| u\|_{L^2(\Omega)}^2,
\end{align}
subject to the following elliptic PDE constraint
\begin{align}\label{eqn:ellipt}
    - \Delta y = f +u \text{ in }\Omega\quad \mbox{with}\quad
	y=0 \mbox{ on } \partial\Omega,
\end{align}
where $u\in L^2(\Omega)$ is the control, $f\in L^2(\Omega)$ and
$y_{d} \in L^2(\Omega)$ are the given source and target function, and the scalar $\lambda>0$ balances the two terms in $J(y,u)$. Problem \eqref{eq:min-J:Diffu}--\eqref{eqn:ellipt} has a unique solution
$(\bar{y},\bar{u})$ \cite[Theorem 2.17]{Trolzstch:2005}. The goal is to approximate the solution $(\bar{y},
\bar{u})$ by an NN pair $(\bar{y}_{\theta}, \bar{u}_{\kappa})$.

\subsubsection{Reduced optimality system}
Problem \eqref{eq:min-J:Diffu}-\eqref{eqn:ellipt} is equivalent to the following  KTT system for $(\bar y,\bar p,\bar u)$ \cite[Chap. 2, Theorem 1.4]{Lions:1971}:
\begin{align}\label{eqn:opt-unconstrained}
\left\{\begin{aligned}
	- \Delta {y} &= f +{u}\text{ in }
	\Omega && \mbox{with}\quad
	{y}=0 \text{ on } \partial\Omega,\\
	- \Delta {p} &= {y}-y_d \text{ in }
	\Omega&& \mbox{with}\quad	{p}=0 \text{ on }\partial\Omega,\\
   {u}&= -\lambda^{-1} p \mbox{ in }\Omega.
\end{aligned}\right.
\end{align}
The last line is the first-order necessary optimality condition of problem \eqref{eq:min-J:Diffu}--\eqref{eqn:ellipt}. Using the condition $u = -\lambda^{-1} p$, problem \eqref{eq:min-J:Diffu}-\eqref{eqn:ellipt} reduces to the following system in $y$ and $p$:
\begin{align}\label{eqn:coupled}
\left\{\begin{aligned}
	- \Delta {y} &= f -\lambda^{-1} p \text{ in }	\Omega && \mbox{with}\quad
	{y}=0 \text{on } \partial\Omega,\\
	- \Delta {p} &= {y}-y_d  \text{ in }	\Omega && \mbox{with}
	\quad {p}=0 \text{ on } \partial\Omega.
\end{aligned}\right.
\end{align}
This problem has a unique solution tuple $(\bar y,\bar p)$, and the optimal $\bar u$ can be recovered ast $\bar u =-\lambda^{-1}\bar p$.

In view of \eqref{eqn:coupled} and the principle of PDE residual minimization, we define the following continuous loss
\begin{align}\label{eqn:loss-pop-uncons}
    \mathcal{L}(y,p)= \norm{\Delta y+f- \lambda^{-1} p }^2_{L^2(\Omega)}+ \alpha_{i} \norm{\Delta p+y-y_d }^2_{L^2(\Omega)}+\alpha^y_b\norm{y}^2_{L^2(\partial\Omega)}+\alpha^p_b \norm{p}^2_{L^2(\partial\Omega)},
\end{align}
where $\alpha_i>0$, $\alpha_b^y>0$ and $\alpha_b^p>0$ are penalty parameters to approximately enforce the PDE residual (of $p$) and zero Dirichlet boundary conditions of $y$ and $p$. Below we let $\bs\alpha = (\alpha_i,\alpha_b^y,\alpha_b^p)$. The idea of utilizing NNs to parameterize the solution and then learning the NN parameters using a loss given by the least-squares of the PDE residual and boundary fitting was first proposed in the 1990s \cite{LeeKang:1990,LagarisLikas:1998,LagarisLikas:2000}, and later greatly popularized by 
\cite{RaissiPerdolarisKarniadakis:2019}, under the name of physics informed neural networks (PINN). It has been successfully applied to a range of problems \cite{Karniadakis:2021nature,CaiMao:2021}. We  explore the use of NN for solving elliptic optimal control.

\subsubsection{OSNN}
First we recall fully connected feedforward NNs. Fix the depth $L\in\mathbb{N}$, and integers $\{n_\ell\}_{\ell=0}^L\subset \mathbb{N}$ ($W:=\max_{\ell=1,\ldots,L}(n_\ell)$ is the width), with $n_0=d$ and $n_L=1$. An NN $v:\mathbb{R}^d\to \mathbb{R}$ is defined recursively by
$v^{(0)}(x) = x$, 
$v^{(\ell)}(x) =  \rho (A^{(\ell)}v^{(\ell-1)}+ b^{(\ell)})$ for $\ell = 1,\ldots, L-1$, and
$v:= v^{(L)}(x) = A^{(L)} v^{(L-1)} + b^{(L)}$,
where $A^{(\ell)}=[a_{ij}^{(\ell)}]\in \mathbb{R}^{n_\ell\times n_{\ell-1}}$, and $b^{(\ell)}=(b_i^{(\ell)})\in
\mathbb{R}^{n_\ell}$ are the weight matrix and bias vector at the $\ell$th layer, collectively denoted by $\theta$. 
$\mathcal{N}_\rho(L,\boldsymbol{n}_L,R)$ denotes the collection of NN functions generated by $\rho$ of depth $L$, total number  $\boldsymbol n_L$ of nonzero weights and each weight bounded by $R$. We focus on  hyperbolic tangent $\rho(t)=\frac{e^{t}-e^{-t}}{e^{t}+e^{-t}}$
and sigmoid $\rho(t)=\frac{1}{1+e^{-t}}$.

We employ two NNs $y_\theta \in \mathcal{Y}=\mathcal{N}_\rho(L,\boldsymbol{n}_L,R)$ and $p_\sigma \in \mathcal{P}=\mathcal{N}_\rho(L,\boldsymbol{n}_L,R)$ to approximate $y$ and $p$, respectively, and then approximate relevant integrals via the standard Monte Carlo method. Let $U(D)$ be the uniform distribution over a set $D$, $|D|$ the Lebesgue measure of $D$, and $\mathbb{E}_\nu$ taking expectation with respect to $\nu$. The loss $\mathcal{L}(y_\theta,p_\sigma)$ can be rewritten as
 \begin{align}\label{Loss:residual}
    \mathcal{L}(y_\theta,p_\sigma) =&|\Omega|\mathbb{E}_{X\sim U(\Omega)} \left[\big(\Delta y_\theta(X)+f(X)- \lambda^{-1} p_\sigma(X)\big)^2\right]+ \alpha_{i}|\Omega|\mathbb{E}_{X\sim U(\Omega)} \left[\big(\Delta p_\sigma(X)+y_\theta(X)-y_d(X) )^2\right]\nonumber\\
    &+ \alpha^y_b|\partial \Omega|\mathbb{E}_{Y\sim U(\partial\Omega)}\left[\big(y_\theta(Y)\big)^2\right]+ \alpha^p_b|\partial \Omega|\mathbb{E}_{Y\sim U(\partial\Omega)}\left[\big(p_\sigma(Y)\big)^2\right].
\end{align}
Let $\{X_i\}_{i=1}^{n_d}$ and $\{Y_j\}_{j=1}^{n_b}$ be i.i.d. samples drawn from $U(\Omega)$ and $U(\partial\Omega)$, respectively, i.e., $\mathbb{X}=\{X_i\}_{i=1}^{n_d}\sim U(\Omega)$ and $\mathbb{Y}=\{Y_j\}_{j=1}^{n_b}\sim U(\partial\Omega)$. Then the empirical loss $\widehat{\mathcal{L}}(y_\theta,p_\sigma)$ is given by
\begin{align}\label{Loss:emperical}
    \widehat{\mathcal{L}}(y_{\theta},p_{\sigma}) =& \frac{|\Omega|}{\nd} \sum_{i=1}^{\nd} \big(\Delta y_{\theta}(X_i)+f(X_i)- \lambda^{-1} p_{\sigma }(X_i)\big)^2+ \alpha_{i} \frac{|\Omega|}{\nd} \sum_{i=1}^{\nd} \big(\Delta p_{\sigma}(X_i)+y_{\theta}(X_i)-y_d(X_i) )^2\nonumber
    \\&+  \alpha^y_b\frac{|\partial \Omega|}{n_b} \sum_{j=1}^{n_b}\big(y_{\theta}(Y_j)\big)^2+ \alpha^p_b\frac{|\partial \Omega|}{n_b} \sum_{j=1}^{n_b}\big(p_{\sigma}(Y_j)\big)^2.
\end{align}
 Let $({y}_{\theta^*},{p}_{\sigma^*})
\in \mathcal{Y}\times \mathcal{P}$ be a minimizer of $ \widehat{\mathcal{L}}(y_{\theta},p_{\sigma})$ over $\mathcal{Y}\times \mathcal{P}$, i.e., with respect to the NN parameter vectors $\theta$ and $\sigma$. In practice, this is commonly achieved by off-shelf optimizers, e.g., L-BFGS \cite{ByrdLu:1995} and Adam \cite{KingmaBa:2015}.  Evaluating the loss $\widehat{\mathcal{L}}(y_{\theta},p_{\sigma})$ requires the derivative of $y_\theta$ and $p_\sigma$ with respect to $x$, and the optimizer requires computing the derivative of $\widehat{\mathcal{L}}(y_\theta,p_\sigma)$ with respect to $\theta$ and $\sigma$. Both can be realized using automatic differentiation \cite{Baydin:2018}. Given a (local) minimizer $(\theta^*,\sigma^*)$, the approximations of $\bar y$ and $\bar p$ are given by $y_{\theta^*}$ and $p_{\sigma^*}$, respectively, and  the approximation $u^*$ of $\bar u$ is given by $-\lambda^{-1} p_{\sigma^*}$.

\begin{remark}
Demo et al \cite{DemoStrazzullo:2023} proposed an approach based on the first-order KKT system for parametric optimal control problems. In the standard elliptic case,  the loss is given by
\begin{align*}
    \mathcal{L}(y,p,u)= &\norm{\Delta y+f+u }^2_{L^2(\Omega)}+ \alpha_{i} \norm{\Delta p+y-y_d }^2_{L^2(\Omega)}\\
     & +\alpha^y_b\norm{y}^2_{L^2(\partial\Omega)} +\alpha^p_b \norm{p}^2_{L^2(\partial\Omega)} + \alpha_i^u \|u+ \lambda^{-1} p\|_{L^2(\Omega)}^2.
\end{align*}
Compared with the proposed scheme \eqref{eqn:loss-pop-uncons}, the loss $\mathcal{L}(y,p,u)$ requires also one NN $u_{\kappa}\in\mathcal{U}$ for approximating the control $u$, and tuning the weight $\alpha_i^u$, which might be nontrivial. Also, the work \cite{DemoStrazzullo:2023} does not discuss the constrained case, where the control $u$ is inherently nonsmooth, and thus numerically challenging to approximate using DNNs \cite{YinYang:2023}. A similar approach to \cite{DemoStrazzullo:2023} was suggested in \cite{BarrySarshar:2022}.
\end{remark}

\subsection{Optimal control with box constraint}\label{subsec:boxconst.}

Now we turn to the constrained case, i.e., the control $u$ is from a box constrain set
$U=\{u\in L^2(\Omega): u_a \leq u\leq u_b \text{ a.e. in } \Omega\},$ with $u_a<u_b$.
Then the optimal tuple $(\bar y,\bar p,\bar u)$ satisfies the following first-order necessary optimality system \cite[p. 67]{Trolzstch:2005}
\begin{align}
	\left\{\begin{aligned}- \Delta {y} &= f + u  \text{ in }
	\Omega && \mbox{with}\quad
	{y}=0\text{ on } \partial\Omega,\\
	- \Delta {p} &= {y}-y_d \text{ in }
	\Omega && \mbox{with}\quad
	{p}=0 \text{ on } \partial\Omega,\\
 (u+\lambda^{-1}&p,v-u) \geq 0,&& \forall v\in U.
 \end{aligned}\right.\label{eqn:coupled-cons-first}
\end{align}
Using the pointwise projection operator $P_{U}$ into the set $U$, defined by
$P_{U}(v)(x) = \min(\max(v(x),u_a),u_b)$, we can rewrite the variational inequality in \eqref{eqn:coupled-cons-first} as
$ u = P_U (-\lambda^{-1}p)$,
and the system is reduced to 
\begin{align}\label{eqn:coupled-cons}
	\left\{\begin{aligned}- \Delta y &= f +P_U(-\lambda^{-1} p) \text{ in }
	\Omega&&
\mbox{with}\quad	{y}=0  \text{ on } \partial\Omega,\\
	- \Delta {p} &= {y}-y_d \text{ in }
	\Omega
&&\mbox{with}\quad	{p}=0 \text{ on } \partial\Omega.
 \end{aligned}\right.
\end{align}

The system \eqref{eqn:coupled-cons} naturally motivates the following continuous loss:
\begin{align*}
    \mathcal{L}(y_\theta,p_\sigma) =& | \Omega| \mathbb{E}_{X\sim U(\Omega)} \left[\big(\Delta y_\theta(X)+f(X)+ P_{U}(-\lambda^{-1} p_\sigma(X))\big)^2\right] + \alpha^y_b  |\partial \Omega| \mathbb{E}_{Y\sim U(\partial\Omega)}\left[\big(y_\theta(Y)\big)^2\right] \\
   & + \alpha_{i} |\Omega| \mathbb{E}_{X\sim U(\Omega)} \left[\big(\Delta p_\sigma(X)+y_\theta(X)-y_d(X) )^2\right]+ \alpha^p_b |\partial \Omega |\mathbb{E}_{Y\sim U(\partial\Omega)}\left[\big(p_\sigma(Y)\big)^2\right].
\end{align*}
Upon approximating the expectations with Monte Carlo using i.i.d. samples $\mathbb{X}=\{X_i\}_{i=1}^{\nd}\sim U(\Omega)$ and $\mathbb{Y}=\{Y_j\}_{j=1}^{n_b}\sim U(\partial\Omega) $, we obtain the following empirical loss
\begin{align*}
    \widehat{\mathcal{L}}(y_{\theta},p_{\sigma}) =& \frac{|\Omega|}{\nd} \sum_{i=1}^{\nd} \Bigl(\Delta y_{\theta}(X_i)+f(X_i)+P_{U}(-\lambda^{-1}) p_{\sigma}(X_i)\Bigr)^2 +  \alpha^y_b\frac{| \partial \Omega|}{n_b}   \sum_{j=1}^{n_b}\big(y_{\theta}(Y_j)\big)^2 \\
    &+ \alpha_{i} \frac{|\Omega |}{n_d} \sum_{i=1}^{\nd} \big(\Delta p_{\sigma}(X_i)+y_{\theta}(X_i)-y_d(X_i) )^2
    + \alpha^p_b\frac{| \partial \Omega |}{n_b} \sum_{j=1}^{n_b}\big(p_{\sigma}(Y_j)\big)^2.
\end{align*}
The empirical loss $\widehat{\mathcal{L}}(y_\theta,p_\sigma)$ is then minimized with respect to the NN parameters $(\theta,\sigma)$, with a minimizer $(\theta^*,\sigma^*)$. The approximation to the optimal control $\bar u$ is given by ${P}_U(-\lambda^{-1}p_{\sigma^*})$. Note that the change from the unconstrained case is rather minimal in terms of implementation.

\subsection{Optimal control for semilinear problems}\label{subsec;semilinprob}
The proposed OSNN applies also to more complex problems and is fairly flexible with the governing PDEs. We illustrate this with semilinear elliptic problems, one popular model class in PDE optimal control. Consider the following optimal control problem
\begin{align}\label{eq:min-J:semi}
\min\; J(y,u):=\tfrac{1}{2}\|y-y_d\|_{L^2(\Omega)}^2+\tfrac{\lambda}{2}\| u\|_{L^2(\Omega)}^2,
\end{align}
subject to the following semilinear elliptic PDE constraint
\begin{align}\label{eqn:semi-ellipt}
    - \Delta y = f(y) +u  \text{ in } \Omega\quad \mbox{with}\quad
	y=0 \text{ on } \partial\Omega.
\end{align}
The function $f:\mathbb{R}\to \mathbb{R}$ should satisfy suitable conditions (e.g., monotone) so that  \eqref{eqn:semi-ellipt} is well-posed \cite[Theorem 4.4]{Trolzstch:2005}. Then the optimal control problem has a locally unique solution $(\bar{y},\bar{u})$ \cite[p. 208]{Trolzstch:2005}. The KKT system for problem \eqref{eq:min-J:semi}--\eqref{eqn:semi-ellipt} is given by \cite[p. 216]{Trolzstch:2005}:
\begin{align*}
\left\{\begin{aligned}
	- \Delta {y} &= f( y) +{u} \text{ in }	\Omega&& \mbox{with}\quad
	{y}=0 \text{ on } \partial\Omega,\\
	- \Delta {p} &=f'(y) p + {y}-y_d \text{ in }
	\Omega&& \mbox{with}\quad
	{p}=0 \text{ on } \partial\Omega,\\
   {u}&= -\lambda^{-1}  p \mbox{ in }\Omega.
\end{aligned}\right.
\end{align*}
Like before, the empirical loss $\widehat{L}(y_\theta,p_\sigma)$ for problem \eqref{eq:min-J:semi}--\eqref{eqn:semi-ellipt} is given by
\begin{align*}
    \widehat{\mathcal{L}}(y_{\theta},p_{\sigma}) =& \frac{|\Omega|}{\nd} \sum_{i=1}^{\nd} \Bigl(\Delta y_{\theta}(X_i)+f(y_\theta(X_i))-\lambda^{-1}p_{\sigma}(X_i)\Bigr)^2 +  \alpha^y_b\frac{| \partial \Omega|}{n_b}   \sum_{j=1}^{n_b}\big(y_{\theta}(Y_j)\big)^2 \\
    &+ \alpha_{i} \frac{|\Omega |}{\nd} \sum_{i=1}^{\nd} \big(\Delta p_{\sigma}(X_i)+f'(y_\theta(X_i))p_\sigma(X_i)+y_{\theta}(X_i)-y_d(X_i) )^2
    + \alpha^p_b\frac{| \partial \Omega |}{n_b} \sum_{j=1}^{n_b}\big(p_{\sigma}(Y_j)\big)^2.
\end{align*}
The adaptation of the approach to semilinear problems with box constraint on the control $U$ is direct.

\subsection{Existing NN based approaches to elliptic optimal control} \label{ssec:existing}

In the literature, three methods have been proposed, i.e., penalty method (PM) \cite{MowlaviNabi:2023}, 
augmented Lagrangian method (ALM) \cite{LuJohnson:2021}, and adjoint oriented neural network (AONN) \cite{YinYang:2023}. We describe the them briefly below. 

PM transforms \eqref{eq:min-J:Diffu}--\eqref{eqn:ellipt} into an unconstrained problem by including the loss $\mathcal{L}_{\rm nn}$:
\begin{align*}
\mathcal{L}_{\rm pm}(y_\theta,u_\kappa)
:=&J(y_\theta,u_\kappa)+\mu \mathcal{L}_{\rm nn}(y_\theta,u_\kappa)\\
=&\big(\tfrac{1}{2}\|y_\theta-y_d\|_{L^2(\Omega)}^2+\tfrac{\lambda}{2}\| u_\kappa\|_{L^2(\Omega)}^2\big) + \mu \big(\tfrac{1}{2}\|F(y_\theta,u_\kappa)\|_{L^2(\Omega)}^2 + \tfrac{\alpha}{2}\|y_\theta\|_{L^2(\partial\Omega)}^2\big),
\end{align*}
where $\mu>0$ is the penalty weight for the loss $\mathcal{L}_{\rm nn}(y_\theta,u_\kappa)=\tfrac{1}{2}\|F(y_\theta,u_\kappa)\|_{L^2(\Omega)}^2 + \tfrac{\alpha}{2}\|y_\theta\|_{L^2(\partial\Omega)}^2$, and $F(y_\theta,u_\kappa)=\Delta y_\theta + f + u_\kappa $. The weight $\mu$ can be fixed \textit{a priori} or determined by numerical continuation \cite{LuJohnson:2021,HuJinZhou:2022}. The loss $\mathcal{L}_{\rm pm}$ becomes more ill-conditioned as $\mu\to\infty$ \cite{Bertsekas:1982,LuJohnson:2021}, and the case of multiple constraints requires lengthy and time-consuming tuning of weights \cite{Krishnapriyan:2021,YinYang:2023}.

ALM mitigates ill-conditioning using Lagrangian multipliers \cite{Bertsekas:1982}. Let $(\cdot,\cdot)$ and $(\cdot,\cdot)_{L^2(\partial\Omega)}$ denote the $L^2(\Omega)$ and $L^2(\partial\Omega)$ inner products. Then the loss $\mathcal{L}_{\rm alm}$ reads
$$
\mathcal{L}_{\rm alm}(y_\theta,u_\kappa,\eta^d,\eta^b):= J(y_\theta,u_\kappa)+\mu \mathcal{L}_{\rm nn}(y_\theta,u_\kappa)+(\eta^d,F(y_\theta,u_\kappa))+ (\eta^b,y_\theta)_{L^2(\partial\Omega)},
$$
where $\eta^d\in L^2(\Omega)$ and $\eta^b\in L^2(\partial\Omega)$ are Lagrangian multipliers for the constraints $F(y,u)=0$ in $\Omega$ and $y=0$ on $\partial\Omega$, respectively, which may also be realized by NNs.
Given $(\eta^d_k,
\eta^b_k)$, we minimize $\mathcal{L}_{\rm alm}(y_\theta,u_\kappa,\eta_k^d,\eta_k^b)$ in $(y_\theta,u_\kappa)$:
$(y_{\theta^k},u_{\kappa^k}) = \arg\min_{(y_\theta,u_\kappa)\in\mathcal{Y}\times\mathcal{P}}  \mathcal{L}_{\rm alm}(y_\theta,u_\kappa,\eta^d_k,\eta^b_k)$,
and then update the Lagrangian multipliers $(\eta^d_{k+1},\eta^b_{k+1})$ by
$\eta^d_{k+1}=\eta^d_k + \mu F(y_{\theta^k},u_{\kappa^k})$ and $ \eta^b_{k+1}=\eta^b_k + \mu \alpha y_{\theta^k}.$
In practice, the integrals are approximated by evaluations on fixed collocation points, and it suffices to discretize the multipliers point-wise \cite{LuJohnson:2021}.

AONN is also based on the KKT system \eqref{eqn:opt-unconstrained}, but utilizes the adjoint $p$ to compute the gradient of the reduced cost $J(u)=J(y(u),u)$ in $u$ only. The total derivative ${\rm d}_uJ(y,u)=\alpha u + p$ is used to update $u$: \begin{equation}\label{eqn:aonn-gd}
    u_{\kappa^{k+1}}= \arg\min_{u_\kappa\in \mathcal{U}}\|u_\kappa - (u_{\kappa^k} - s^k{\rm  d}_uJ(y_{\theta^k},u_{\kappa^k}))\|^2_{L^2(\Omega)},
\end{equation}
where the step size $s^k>0$ should be suitably chosen. The minimization step projects the gradient update into the set $\mathcal{U}$, i.e., $u$ is approximated by an NN  $u_\kappa\in\mathcal{U}$. 

The presence of box constraint  $u\in U$ requires a slightly different treatment. In PM, one can add an additional penalty term:
\begin{equation}\label{eqn:pm-constraint}
\mathcal{L}_{\rm pm\text{-}c}(y_\theta,u_\kappa)= J(y_\theta,u_\kappa)+\mu \mathcal{L}_{\text{nn}}(y_\theta,u_\kappa)+\tfrac{\mu'}{2}\|u_\kappa - {P}_{U}(u_\kappa)\|_{L^2(\Omega)}^2,
\end{equation}
where $\mu'>0$, and $u_\kappa\in \mathcal{U}$ is an NN approximation of $u$. However, the approximation $u_\kappa$ may be infeasible.
One can also enforce the constraint by applying  $P_U$ on $u_\kappa$ directly:
$
\mathcal{L}_{\rm pm}(y_\theta,u_\kappa)= J(y_\theta,P_U(u_\kappa))+\mu \mathcal{L}_{\text{nn}}(y_{\theta},P_U(u_{\kappa})).
$
This may lead to vanishing gradient when the constraint is active, causing numerical issues \cite{HaghighatRaissi:2021,LuJohnson:2021}.
AONN treats the box constraint $u\in U$ by
$$
u_{\kappa^{k+1}}=\arg\min_{u_\kappa\in\mathcal{U}}\|u_\kappa - P_U(u_{\kappa^k}-s^k d_uJ(y_{\theta^k},u_{\kappa^k}))\|_{L^2(\Omega)}^2,
$$
where the initial guess for the NN parameters $\kappa^{k+1}$ is $\kappa^{k}$. For OSNN, the constraint is enforced by substituting  $P_U(p_\sigma)$ into the PDE system. 

\section{Error analysis}\label{sec:erranalysis}
Now we analyze the proposed OSNN for linear elliptic optimal control problems, by deriving error bounds on the NN approximation $(y_{\hat\theta^*},p_{\hat\sigma^*},u_{\hat\sigma^*})$, under the assumption that $(\hat\theta^*,\hat\sigma^*)$ is a global minimizer of the empirical loss $\widehat{\mathcal{L}}$. The error analysis of neural PDE solvers (including physics informed neural networks and deep Ritz method) for second-order elliptic PDEs has received much attention \cite{JiaoLaiLi:2022pinn,ShinDarbon:2020,MishraMolinaro:2023,LuChenLu:2021,LuLu:2021,HuJinZhou:2022,HuJinZhou:2024}. The analysis below builds on these works, but with two distinct contributions. First, we employ offset Rademacher complexity \cite{LiangRakhlin:2015,DuanJiao:2023fast} to derive faster convergence rates than existing works for the statistical error. Second, we derive the weak coercivity estimates, for the coupled system (which may include the projection operator $P_U$) Lemmas \ref{lem:fund-est1} and \ref{lem:fund-est2} below. These are the main technical contributions of the work. We leave the technically more involved semilinear case to future works, especially the weak coercivity estimate as Lemma \ref{lem:fund-est1}. Moreover, in practice, the optimizer only finds a local minimizer of the empirical loss $\widehat{\mathcal{L}}(y_\theta, p_\sigma)$, due to its highly nonconvex landscape. However, the analysis of the optimization error is largely open, and we do not further pursue the issue below. 

\subsection{Fundamental estimates}
First, we provide two fundamental stability results bounding the error of the NN approximations (relative to the optimal tuple $(\bar y,\bar p,\bar u)$) in terms of the population loss. Such weak type coercivity results enable reducing the error analysis to properly estimating the loss function. We first discuss the unconstrained case. The main challenge lies in the weak imposition of the Dirichlet boundary condition.
\begin{lemma}\label{lem:fund-est1}
Let $(\bar y,\bar p)$ be the solution tuple to the system \eqref{eqn:coupled} with $ \bar u=-\lambda^{-1}\bar p$. Then for any $(y_{\theta},p_\sigma)\in \mathcal{Y}\times \mathcal{P}$, with $u_{\sigma}=-\lambda^{-1}p_{\sigma}$, the following estimate holds
\begin{align*}
 \norm{\bar{y}-y_{\theta}}^2_{L^2(\Omega)}+\norm{\bar{p}-p_{\sigma}}^2_{L^2(\Omega)}+\norm{\bar{u}-u_{\sigma}}^2_{L^2(\Omega)}\leq c(\bs{\alpha},\lambda) \mathcal{L}(y_{\theta},p_{\sigma}).
\end{align*}
\end{lemma}
\begin{proof}
To deal with the nonzero Dirichlet boundary conditions with $y_\theta$ and $p_\sigma$, we define
\begin{align}\label{HarmonicExtn}
\left\{\begin{aligned}
    -\Delta \zeta^{y} &= 0, && \text{in} \; \Omega,\\
    \zeta^{y} & = y_{\theta}, && \text{on} \; \partial\Omega,
    \end{aligned}\right.
\quad\mbox{and}\quad\left\{\begin{aligned}
    -\Delta \zeta^{p} &= 0, && \text{in} \; \Omega,\\
    \zeta^{p} & = p_{\sigma}, && \text{on} \; \partial\Omega.
\end{aligned}\right.
\end{align}
That is, $\zeta^y$ and $\zeta^p$ are harmonic extensions of the $L^2(\partial\Omega)$ boundary data $y_\theta$ and $p_\sigma$, respectively. By the standard elliptic regularity theory \cite[Theorem 4.2, p. 870]{Berggren:2004},
the following stability estimates hold
\begin{equation}\label{eqn:stab-aux}
\|\zeta^y\|_{L^2(\Omega)}\leq c\|y_{\theta}\|_{L^2(\partial\Omega)}\quad\mbox{and}\quad \|\zeta^p\|_{L^2(\Omega)}\leq c \|p_{\sigma}\|_{L^2(\partial\Omega)}.
\end{equation}
Let $e_y =\bar{y}-y_{\theta}$ and $e_p =\bar{p}-p_{\sigma}$, also $\tilde{e}_y=e_y+\zeta^{y}$
and $\tilde{e}_p=e_p+\zeta^{p}$. Then $\tilde{e}_y$ and $\tilde{e}_p$ satisfy
\begin{equation*}
\begin{aligned}
    -\Delta \tilde{e}_y +\lambda^{-1} \tilde{e}_p &=f-\lambda^{-1}p_{\sigma}+\Delta y_{\theta}+\lambda^{-1}\zeta^p\text{ in } \Omega&& \mbox{with}
\quad    \tilde{e}_y = 0\text{ on } \partial\Omega,\\
    -\Delta \tilde{e}_p -\tilde{e}_y &= y_{\theta}-y_d+\Delta p_{\sigma}-\zeta^y\text{ in } \Omega&& \mbox{with}\quad
    \tilde{e}_p=0\text{ on } \partial\Omega. 
\end{aligned}
\end{equation*}
Multiplying the first line by $\tilde{e}_y$ and the second by $\lambda^{-1}\tilde{e}_p$ and integrating over $\Omega$ yield
\begin{align*}
    \|\nabla \tilde{e}_y\|^2_{L^2(\Omega)}+\lambda^{-1}\|\nabla \tilde{e}_p\|^2_{L^2(\Omega)}
    \leq &\|f-\lambda^{-1}p_{\sigma}+\Delta y_{\theta}+\lambda^{-1}\zeta^p\|_{L^2(\Omega)}\|\tilde{e}_y\|_{L^2(\Omega)}\\
    &+\lambda^{-1}\|y_{\theta}-y_d+\Delta p_{\sigma }-\zeta^y\|_{L^2(\Omega)}\|\tilde{e}_p\|_{L^2(\Omega)}.
\end{align*}
Then by the triangle inequality and the stability estimates in \eqref{eqn:stab-aux}, we deduce
\begin{align*}
    \|\nabla \tilde{e}_y\|^2_{L^2(\Omega)}+\lambda^{-1}\|\nabla \tilde{e}_p\|^2_{L^2(\Omega)}
  \leq&\|f-\lambda^{-1}p_{\sigma}+\Delta y_{\theta}\|_{L^2(\Omega)}\|\tilde{e}_y\|_{L^2(\Omega)}+c\lambda^{-1}\|p_{\sigma}\|_{L^2(\partial\Omega)} \|\tilde{e}_y\|_{L^2(\Omega)}\\
    &+\lambda^{-1}\|y_{\theta}-y_d+\Delta p_{\sigma}\|_{L^2(\Omega)}\|\tilde{e}_p\|_{L^2(\Omega)}+c\lambda^{-1} \|y_{\theta}\|_{L^2(\partial\Omega)}\|\tilde{e}_p\|_{L^2(\Omega)}.
\end{align*}
Now by Poincar\'e inequality, trace inequality, and Young's inequality, we obtain
\begin{align}\label{ErrEstm1}
\|\tilde{e}_y\|^2_{L^2(\Omega)}+\lambda^{-1}\|\tilde{e}_p\|^2_{L^2(\Omega)}\leq& c\big(\|f-\lambda^{-1}p_{\sigma }+\Delta y_{\theta}\|^2_{L^2(\Omega)}+\lambda^{-1}\|y_{\theta}-y_d+\Delta p_{\sigma}\|^2_{L^2(\Omega)}\nonumber\\
     &+\lambda^{-1} \|y_{\theta}\|^2_{L^2(\partial\Omega)}+\lambda^{-1}\|p_{\sigma}\|^2_{L^2(\partial\Omega)}\big).
\end{align}
Meanwhile, by the triangle inequality, we have
\begin{align*}
     \|e_y\|^2_{L^2(\Omega)} \leq c(\|\tilde{e}_y\|^2_{L^2(\Omega)}+\|\zeta^y\|^2_{L^2(\Omega)})
      \leq& c(\|\tilde{e}_y\|^2_{L^2(\Omega)}+\|y_{\theta}\|^2_{L^2(\partial\Omega)})
      \leq c(\boldsymbol{\alpha},\lambda)\mathcal{L}(y_{\theta},p_{\sigma}),\\
    \|e_p\|^2_{L^2(\Omega)} \leq c(\|\tilde{e}_p\|^2_{L^2(\Omega)}+\|\zeta^p\|^2_{L^2(\Omega)})\leq& c(\|\tilde{e}_p\|^2_{L^2(\Omega)}+\|p_\sigma\|^2_{L^2(\partial\Omega)})\leq c(\boldsymbol\alpha,\lambda) \mathcal{L}(y_{\theta},p_{\sigma}).
\end{align*}
Moreover, for the error $e_u=\bar u-u_\sigma$, from the identity $\|e_u\|_{L^2(\Omega)}=\lambda^{-1}\|e_p\|_{L^2(\Omega)}$, we deduce
$\|e_u\|^2_{L^2(\Omega)} \leq c(\boldsymbol\alpha,\lambda) \mathcal{L}(y_{\theta},p_{\sigma})$.
Combining the preceding three estimates completes the proof of the lemma.
\end{proof}

\begin{remark}
Unlike the standard FEM discretization, where the zero Dirichlet boundary condition can be imposed directly, NNs often do not satisfy the condition exactly, due to their nonlocal nature. We have adopted an $L^2(\partial\Omega)$ penalty to enforce approximately the boundary condition. This induces a type of consistency error, which also precludes obtaining error bounds in the $H^1(\Omega)$ norm. The coercivity estimate bounds only the difference between the approximate state $y_\theta$ and optimal one $\bar y$, but not directly the state $y(u_\sigma)$ (corresponding to the approximate control $u_\sigma$) to $\bar y$.
Moreover, it follows from the estimate \eqref{ErrEstm1} that one should choose $\alpha_i$ as $\lambda^{-1}$. The boundary penalty parameters $\alpha_b^y$ and $\alpha_b^p$ should be scaled accordingly.
\end{remark}

Next, we discuss the constrained case, which is more involved, due to the presence of the projection operator $P_U$. To overcome the challenge, we employ the idea of reconstruction.
\begin{lemma}\label{lem:fund-est2}
Let $(\bar y, \bar p)$ be the solution tuple of the system \eqref{eqn:coupled-cons}, with $\bar u = P_U(-\lambda^{-1}\bar p)$. Then for any $(y_{\theta},p_\sigma)\in \mathcal{Y}\times \mathcal{P}$, and $u_{\sigma}={P}_{U}(-\lambda^{-1}p_{\sigma})$, the following estimate holds
\begin{align*}
\norm{\bar{y}-y_{\theta}}^2_{L^2(\Omega)}+\norm{\bar{p}-p_{\sigma}}^2_{L^2(\Omega)}+\norm{\bar{u}-u_{\sigma}}^2_{L^2(\Omega)}\leq c(\boldsymbol\alpha,\lambda) \mathcal{L}(y_{\theta},p_{\sigma}).
\end{align*}
\end{lemma}
\begin{proof}
The proof employs crucially the following auxiliary functions $R\bar y,R\bar p \in H_0^1(\Omega)$ satisfying
\begin{align}
\left\{\begin{aligned}\label{AuxSysCons}
    -\Delta R\bar{y} &= f+u_{\sigma}, && \text{in} \; \Omega,\\
    R\bar{y} & = 0, && \text{on} \; \partial\Omega,
    \end{aligned}\right.
\quad\mbox{and}\quad\left\{\begin{aligned}
    -\Delta R\bar{p} &= y_{\theta}-y_d, && \text{in} \; \Omega,\\
    R\bar{p} & = 0, && \text{on} \; \partial\Omega.
\end{aligned}\right.
\end{align}
Then subtracting  \eqref{AuxSysCons} from \eqref{eqn:coupled} yields
\begin{equation*}
\begin{aligned}
   -\Delta (\bar{y}-R\bar{y}) &= \bar{u}-u_{\sigma} \text{ in }  \Omega& \mbox{with}\quad
    \bar{y}-R\bar{y}  = 0  \text{ on } \partial\Omega,\\
   -\Delta (\bar{p}-R\bar{p}) &= \bar{y}-y_{\theta}  \text{ in } \Omega&\mbox{with}\quad
    \bar{p}-R\bar{p}  = 0  \text{ on }  \partial\Omega.
\end{aligned}
\end{equation*}
Now multiplying the first line by $(\bar{p}-R\bar{p})$ and the second line by $(\bar{y}-R\bar{y})$, and integrating over $\Omega$ yield
\begin{align}\label{AuxEstm11}
 (\bar{p}-R\bar{p},\bar{u}-u_{\sigma})=(\bar{y}-R\bar{y},\bar{y}-y_{\theta}).
\end{align}
It follows from the first-order necessary optimality condition in \eqref{eqn:coupled-cons-first} that
\begin{align}\label{ContVI}
    (\bar{p}+\lambda \bar{u},v-\bar{u})\geq 0 ,\quad\forall v\in U.
\end{align}
 Also, since $u_{\sigma}={P}_{U}(\lambda^{-1}p_{\sigma})\in U$,  we have
 \begin{align}\label{SemiDiscVI}
   (p_{\sigma}+\lambda u_{\sigma},v-u_{\sigma})\geq 0,\quad \forall v \in U.
 \end{align}
Putting $v=u_{\sigma}\in U$ in \eqref{ContVI}, $v=\bar{u}\in U$ in \eqref{SemiDiscVI} and adding the resulting inequalities give
\begin{align}\label{eqn:coerc-cc}
    (\bar{p}-p_{\sigma},u_{\sigma}-\bar{u})-\lambda \norm{\bar{u}-u_{\sigma}}^2_{L^2(\Omega)}\geq 0.
\end{align}
Meanwhile, it follows from the identity \eqref{AuxEstm11} that
\begin{align*}
    (\bar{p}-p_{\sigma},u_{\sigma}-\bar{u}) & = (\bar{p}-R\bar{p},u_{\sigma}-\bar{u})+(R\bar{p}-p_{\sigma},u_{\sigma}-\bar{u}) \\
    & = -(\bar{y}-R\bar{y},\bar{y}-y_{\theta})+(R\bar{p}-p_{\sigma},u_{\sigma}-\bar{u})\\
    & = -\norm{\bar{y}-R\bar{y}}^2_{L^2(\Omega)}-(\bar{y}-R\bar{y},R\bar{y}-y_{\theta})+(R\bar{p}-p_{\sigma},u_{\sigma}-\bar{u}) .
\end{align*}
Together with the inequality \eqref{eqn:coerc-cc} and Cauchy-Schwarz inequality, we obtain
\begin{align*}
        \norm{\bar{y}-R\bar{y}}^2_{L^2(\Omega)}+\lambda \norm{\bar{u}-u_{\sigma}}^2_{L^2(\Omega)}&\leq (\bar{y}-R\bar{y},y_{\theta}-R\bar{y})+(R\bar{p}-p_{\sigma},u_{\sigma}-\bar{u})\\
        &\leq\|\bar{y}-R\bar{y}\|_{L^2(\Omega)}\|y_{\theta}-R\bar{y}\|_{L^2(\Omega)}+\|R\bar{p}-p_{\sigma}\|_{L^2(\Omega)}\|u_{\sigma}-\bar{u}\|_{L^2(\Omega)}.
\end{align*}
Then an application of Young's inequality leads to
\begin{align*}
    \norm{\bar{y}-R\bar{y}}^2_{L^2(\Omega)}+\lambda \norm{\bar{u}-u_{\sigma}}^2_{L^2(\Omega)}\leq\|y_{\theta}-R\bar{y}\|_{L^2(\Omega)}^2+\lambda^{-1}\|R\bar{p}-p_{\sigma}\|_{L^2(\Omega)}^2.
\end{align*}
By the triangle inequality, we have
\begin{align}\label{Aus:EstmState}
    \norm{\bar{y}-y_{\theta}}_{L^2(\Omega)}&\leq c(\lambda) \big({ \norm{R\bar{y}-y_{\theta}}_{L^2(\Omega)}+\norm{R\bar{p}-p_{\sigma}}_{L^2(\Omega)}}\big).
\end{align}
By the standard elliptic regularity theory, we have $\norm{\bar{p}-R\bar{p}}_{L^2(\Omega)}\leq c \norm{\bar{y}-y_{\theta}}_{L^2(\Omega)}$.
Then using \eqref{Aus:EstmState}, we obtain the following estimate
\begin{align*}
    &\norm{\bar{p}-p_{\sigma}}_{L^2(\Omega)}\leq \norm{\bar{p}-R\bar{p}}_{L^2(\Omega)}+\norm{R\bar{p}-p_{\sigma}}_{L^2(\Omega)}\\
    \leq &c \norm{\bar{y}-y_{\theta}}_{L^2(\Omega)}+\norm{R\bar{p}-p_{\sigma}}_{L^2(\Omega)}
    \leq c(\lambda) \big({ \norm{R\bar{y}-y_{\theta}}_{L^2(\Omega)}+\norm{R\bar{p}-p_{\sigma}}_{L^2(\Omega)}}\big).
\end{align*}
Next, let $e^R_y=R\bar{y}-y_{\theta}$ and $\tilde{e}^R_y=e^R_y+\zeta^y$, and similarly, let $e^R_p=R\bar{p}-p_{\sigma}$ and $\tilde{e}^R_p=e^R_p+\zeta^p$, where $\zeta^y$ and $\zeta^p$ are harmonic extensions defined in \eqref{HarmonicExtn}. Clearly, $\tilde{e}^R_y$ and $\tilde{e}^R_p$ satisfy
\begin{equation*}
\begin{aligned}
    -\Delta \tilde{e}^R_y &= f+u_{\sigma}+\Delta y_{\theta}\text{ in } \Omega && \mbox{with}\quad
    \tilde{e}^R_y  = 0 \text{ on } \partial\Omega,\\
    -\Delta \tilde{e}^R_p &= y_{\theta}-y_d+\Delta p_{\sigma}  \text{ in } \Omega &&\mbox{with}\quad
    \tilde{e}^R_p  = 0 \mbox{ on } \partial\Omega.
\end{aligned}
\end{equation*}
By the standard elliptic regularity estimate, we have
\begin{align*}
\norm{\tilde{e}^R_y}_{L^2(\Omega)} \leq c\norm{f+u_{\sigma}+\Delta y_{\theta}}_{L^2(\Omega)}\quad\mbox{and}\quad \norm{\tilde{e}^R_p}_{L^2(\Omega)}\leq c \norm{y_{\theta}-y_d+\Delta p_{\sigma}}_{L^2(\Omega)}.
\end{align*}
Hence, we have
\begin{equation*}
\begin{aligned}
    \norm{R\bar{y}-y_{\theta}}_{L^2(\Omega)}&=\norm{e^R_y}_{L^2(\Omega)}=\norm{\tilde{e}^R_y-\zeta^y}_{L^2(\Omega)}\leq \norm{\tilde{e}^R_y}_{L^2(\Omega)}+\norm{\zeta^y}_{L^2(\Omega)}\\
    &\leq c\norm{f+u_{\sigma}+\Delta y_{\theta}}_{L^2(\Omega)}+c\norm{y_{\theta}}_{L^2(\partial\Omega)},\\
    \norm{R\bar{p}-p_{\sigma}}_{L^2(\Omega)}&=\norm{e^R_p}_{L^2(\Omega)}=\norm{\tilde{e}^R_p-\zeta^p}_{L^2(\Omega)}\leq \norm{\tilde{e}^R_p}_{L^2(\Omega)}+\norm{\zeta^p}_{L^2(\Omega)}\\
    &\leq c\norm{y_{\theta}-y_d+\Delta p_{\sigma}}_{L^2(\Omega)}+c\norm{p_{\sigma}}_{L^2(\partial\Omega)}.
\end{aligned}
\end{equation*}
Therefore, we obtain the desired estimates and complete the proof of the lemma.
\end{proof}

\subsection{Generalization error}

Note that the coercivity estimates in Lemmas \ref{lem:fund-est1} and \ref{lem:fund-est2} involve the population loss $\mathcal{L}$, whereas in practice we only minimize the empirical loss $\widehat{\mathcal{L}}$, which has a different global minimizer $(\widehat y,\widehat p)$ than that of $\mathcal{L}$. Moreover, due to the randomness of the sampling points $\mathbb{X}$ and $\mathbb{Y}$, the minimizer $(\widehat y,\widehat p)$ is random. The analysis of the quantity $\mathcal{L}(\widehat y,\widehat p)$ is commonly known as generalization error. One commonly used tool in statistical learning theory for analyzing generalization error is Rademacher complexity \cite{AnthonyBartlett:1999,BartlettMendelson:2002}, which measures the complexity of a collection of functions by the correlation between function values and the Rademacher random variable. In this work, we bound the generalization error via offset Rademacher complexity. This tool was first
introduced by Liang et al \cite{LiangRakhlin:2015} to prove sharp bounds for a two-step star estimator, and recently further
developed in \cite{DuanJiao:2023fast}. It is a penalized version of
Rademacher complexity which localizes $\mathcal{F}$ adaptively according to the magnitude of $f^2$. This technique will allow us to derive sharper estimates on the statistical error. 
\begin{definition}
Let $\mathcal{F}$
be a class of measurable functions from $\mathcal{X}$ to $\mathbb{R}$, $\mu_X$ a probability
distribution of $X$ and $\mathbb{X}:=\{X_i\}_{i=1}^n$ i.i.d. samples draw from $\mu_X$. Let
$\{\tau_i\}_{i=1}^n$ be i.i.d. Rademacher random variables, i.e., $\mathbb{P}(\tau_i=1)=\mathbb{P}(\tau_i=-1)=\frac12$. Then the empirical
offset Rademacher complexity of $\mathcal{F}$ is defined as
\begin{equation*}
  \mathcal{R}_n^{\rm off}(\mathcal{F},\beta|\mathbb{X}):=\mathbb{E}_\tau \Big[\sup_{f\in\mathcal{F}}\frac{1}{n}\sum_{i=1}^n\tau_if(X_i)-\beta f(X_i)^2|\mathbb{X}\Big],
\end{equation*}
for some $\beta>0$, and the offset Rademacher complexity of $\mathcal{F}$ is defined by
\begin{align*}
  \mathcal{R}_n^{\rm off}(\mathcal{F},\beta): = \mathbb{E}_{\mathbb{X}} \mathcal{R}_n^{\rm off}(\mathcal{F},\beta|\mathbb{X})
  = \mathbb{E}_{\mathbb{X},\tau} \Big[\sup_{f\in\mathcal{F}}\frac{1}{n}\sum_{i=1}^n\tau_if(X_i)-\beta f(X_i)^2\Big].
\end{align*}
\end{definition}

First, we give a fundamental estimate of the generalization error in terms of the approximation error and statistical error (using offset Rademacher complexity). We define two sets $\mathcal{G}_d$ and $\mathcal{G}_b$ of mappings from $X\in\Omega \mapsto \mathbb{R}$ and $Y\in \partial\Omega \mapsto \mathbb{R}$, respectively, by
\begin{align*}
   \mathcal{G}_d&=\{|\Omega|((\Delta y+f-\lambda^{-1} p)(X))^{2}+\alpha_{i}|\Omega|((\Delta p+y-y_d)(X))^{2}: y\in\mathcal{Y},p\in\mathcal{P}\},\\
   \mathcal{G}_b&=\{\alpha_{b}^{y}|\partial\Omega|y(Y)^{2}+\alpha_{b}^{p}|\partial\Omega|p(Y)^{2}:y\in\mathcal{Y},p\in\mathcal{P}\}.
\end{align*}
Below we write $g_d(y, p, X)=|\Omega|((\left(\Delta y+f-\lambda^{-1} p\right)(X))^{2}+\alpha_{i}((\Delta p+y-y_d)(X))^{2})\in \mathcal{G}_d$, and $g_b(y,p,Y)=|\partial\Omega|(\alpha_{b}^{y}y(Y)^{2}+\alpha_{b}^{p}p(Y)^{2})\in\mathcal{G}_b$, and let $b_d=\sup_{g\in\mathcal{G}_d}\|g\|_{L^\infty(\Omega)}$ and $b_b=\sup_{g\in\mathcal{G}_b}\|g\|_{L^\infty(\partial\Omega)}$.
\begin{theorem}\label{thm:err-decomp}
For any minimizer $(\hat y,\hat p)\in\mathcal{Y}\times \mathcal{P}$ of  $\widehat{\mathcal{L}}(y,p)$, the following estimate holds
\begin{align*}
  \mathbb{E}_{\mathbb{X},\mathbb{Y}}[\mathcal{L}(\hat y,\hat p)] \leq 4\mathcal{R}_{\nd}^{\rm off}\big(\mathcal{G}_d,(2b_d)^{-1}\big) + 4\mathcal{R}_{n_b}^{\rm off}\big(\mathcal{G}_b,(2b_b)^{-1}\big) + 3\inf_{(y,p)\in\mathcal{Y}\times\mathcal{P}}\mathcal{L}(y,p).
\end{align*}
\end{theorem}
\begin{proof}
Clearly, we have the following decomposition for any $(y,p)\in \mathcal{Y}\times\mathcal{P}$,
\begin{align*}
  &\mathbb{E}_{\mathbb{X},\mathbb{Y}}[\mathcal{L}(\hat y,\hat p)] = \mathbb{E}_{\mathbb{X},\mathbb{Y}}\left[\mathbb{E}_{X,Y} \left[{\gd}(\hat y,\hat p,X)+ g_b(\hat y,\hat p,Y)\right]\right]  \\
  = & \mathbb{E}_{\mathbb{X},\mathbb{Y}}\Big[\mathbb{E}_{X,Y} \left[{\gd}(\hat y,\hat p,X)+ g_b(\hat y,\hat p,Y)\right]-\frac{3}{{\nd}}\sum_{i=1}^{{\nd}} {\gd}(\hat y,\hat p,X_i)-\frac{3}{n_b}\sum_{i=1}^{n_b} g_b(\hat y,\hat p,Y_i)\Big]\\
   &+ \mathbb{E}_{\mathbb{X},\mathbb{Y}}\Big[\frac{3}{{\nd}}\sum_{i=1}^{{\nd}}
 {\gd}(\hat y,\hat p,X_i)+\frac{3}{n_b}\sum_{i=1}^{{\nd}} g_b(\hat y,\hat p,Y_i)\Big] := {\rm I} + {\rm II}.
\end{align*}
Since $(\hat y,\hat p)\in\mathcal{Y}\times\mathcal{P}$ is a minimizer, i.e.,
$\widehat{\mathcal{L}}(\hat y,\hat p)\leq \widehat{\mathcal{L}}(y,p)$ for any $ (y,p)\in\mathcal{Y}\times\mathcal{P}$, we deduce
\begin{align*}
  {\rm II} \leq  \inf_{(y,p)\in\mathcal{Y}\times\mathcal{P}} \mathbb{E}_{\mathbb{X},\mathbb{Y}}\Big[\frac{3}{{\nd}}\sum_{i=1}^{{\nd}} {\gd}(y,p,X_i)+\frac{3}{n_b}\sum_{i=1}^{n_b} g_b(y,p,Y_i)\Big] = 3\inf_{(y,p)\in\mathcal{Y}\times\mathcal{P}}\mathcal{L}(y,p).
\end{align*}
Meanwhile, the convexity of supremum and Jensen's inequality imply
\begin{align*}
 {\rm I} \leq& \sup_{(y,p)\in\mathcal{Y}\times\mathcal{P}}\mathbb{E}_{\mathbb{X},\mathbb{Y}}\Big[\mathbb{E}_{X,Y} [{\gd}(y,p,X)+ g_b(y,p,Y)]-\frac{3}{{\nd}}\sum_{i=1}^{{\nd}} {\gd}(y,p,X_i)-\frac{3}{n_b}\sum_{i=1}^{n_b} g_b(y, p,Y_i)\Big]\\
  \leq & \mathbb{E}_{\mathbb{X},\mathbb{Y}} \sup_{(y,p)\in \mathcal{Y}\times\mathcal{P}}\Big[\mathbb{E}_{X} [{\gd}( y, p,X)]-\frac{3}{{\nd}}\sum_{i=1}^{{\nd}} {\gd}(y,p,X_i)+\mathbb{E}_Y[g_b( y, p,Y)]-\frac{3}{n_b}\sum_{i=1}^{n_b} g_b(y, p,Y_i)\Big]\\
  \leq & \mathbb{E}_{\mathbb{X}} \sup_{(y,p)\in \mathcal{Y}\times\mathcal{P}}\Big[\mathbb{E}_{X} [{\gd}( y, p,X)]-\frac{3}{{\nd}}\sum_{i=1}^{{\nd}} {\gd}(y,p,X_i)\Big]
   +\mathbb{E}_{\mathbb{Y}} \sup_{(y,p)\in \mathcal{Y}\times\mathcal{P}}\Big[\mathbb{E}_{Y} [g_b( y, p,Y)]-\frac{3}{n_b}\sum_{i=1}^{n_b} g_b(y, p,Y_i)\Big].
\end{align*}
Next we bound the two terms, denoted by ${\rm I}_d$ and ${\rm I}_b$. Note that $0\leq {\gd}(\cdot,\cdot,X)\leq b_d$ for all $X\in\Omega$, which implies ${\gd}^2(\cdot,\cdot,X)\leq b_d {\gd}(\cdot,\cdot,X)$. Consequently,
\begin{align*}
 {\rm I}_d \leq & \mathbb{E}_{\mathbb{X}} \sup_{(y,p)\in \mathcal{Y}\times\mathcal{P}}\Big[2\mathbb{E}_{X} {\gd}( y, p,X)-\frac{1}{b_d}\mathbb{E}_{X} {\gd}( y, p,X)^2-\frac{2}{{\nd}}\sum_{i=1}^{{\nd}} {\gd}( y, p,X_i)-\frac{1}{b_d{\nd}}\sum_{i=1}^{{\nd}} g_d( y, p,X_i)^2\Big].
\end{align*}
Now we introduce independent copies of $\mathbb{X}$ and $\mathbb{Y}$, i.e., $\mathbb{X}':=\{X_i'\}_{i=1}^{{\nd}}$ and $\mathbb{Y}'=\{Y_i'\}_{i=1}^{n_b}$, and let $\tau=\{\tau_i\}_{i=1}^{n_d}$ be a sequence of i.i.d. Rademacher random variables independent of $\mathbb{X}$ and $\mathbb{X}'$. By the technique of symmetrization, the convexity of supremum, and Jensen's inequality, we obtain
\begin{align*}
{\rm I}_d  \leq & \mathbb{E}_{\mathbb{X}} \sup_{(y,p)\in \mathcal{Y}\times\mathcal{P}}\Big[\mathbb{E}_{\mathbb{X}'}\Big[\frac{2}{{\nd}}\sum_{i=1}^{{\nd}} {\gd}(y,p,X_i')-\frac{1}{b_d{\nd}}\sum_{i=1}^{\nd} {\gd}( y, p,X_i')^2\Big]-\frac{2}{\nd}\sum_{i=1}^{\nd} {\gd}( y, p,X_i)-\frac{1}{b_d {\nd}}\sum_{i=1}^{{\nd}}{\gd}( y, p,X_i)^2\Big]\\
  \leq & \mathbb{E}_{\mathbb{X}} \mathbb{E}_{\mathbb{X}'}\sup_{(y,p)\in \mathcal{Y}\times\mathcal{P}}\Big[\frac{2}{{\nd}}\sum_{i=1}^{{\nd}}( {\gd}(y,p,X_i')- {\gd}(y,p,X_i))-\frac{1}{b_d\nd}\sum_{i=1}^{{\nd}}({\gd}( y, p,X_i')^2+ {\gd}( y, p,X_i)^2)\Big]\\
  = & \mathbb{E}_{\mathbb{X}} \mathbb{E}_{\mathbb{X}'}\mathbb{E}_\tau\sup_{(y,p)\in \mathcal{Y}\times\mathcal{P}}\Big[\frac{2}{{\nd}}\sum_{i=1}^{{\nd}}\tau_i( {\gd}(y,p,X_i')- {\gd}(y,p,X_i))-\frac{1}{b_d{\nd}}\sum_{i=1}^{{\nd}}({\gd}( y, p,X_i')^2+ {\gd}( y, p,X_i)^2)\Big]\\
  =& 2 \mathbb{E}_{\mathbb{X}'}\mathbb{E}_\tau\sup_{(y,p)\in \mathcal{Y}\times\mathcal{P}}\frac{1}{{\nd}}\sum_{i=1}^{{\nd}}\Big(\tau_i {\gd}(y,p,X_i')- \frac{1}{2b_d} {\gd}(y,p,X_i')^2\Big)\\
   &+ 2\mathbb{E}_\mathbb{X}\mathbb{E}_\tau\sup_{(y,p)\in\mathcal{Y}\times\mathcal{P}}\frac{1}{{\nd}}\sum_{i=1}^{{\nd}}\Big(-\tau_i{\gd}( y, p,X_i)-\frac{1}{2b_d} {\gd}( y, p,X_i)^2\Big)\\
  =  & 2 \mathbb{E}_{\mathbb{X}'} \mathcal{R}_n^{\rm off}\big(\mathcal{G}_d, (2b_d)^{-1}|\mathbb{X}'\big) + 2\mathbb{E}_\mathbb{X}\mathcal{R}_{{\nd}}^{\rm off}\big(\mathcal{G}_d,(2b_d)^{-1}|\mathbb{X}\big) = 4\mathcal{R}_{{\nd}}^{\rm off}\big(\mathcal{G}_d, (2b_d)^{-1}\big).
\end{align*}
Similarly, we can deduce
${\rm I}_b \leq 4\mathcal{R}_{n_b}^{\rm off}\big(\mathcal{G}_b, (2b_b)^{-1}\big)$.
This completes the proof of the theorem.
\end{proof}

By Theorem \ref{thm:err-decomp}, the generalization error can be decomposed into the approximation error $\mathcal{E}_{app}=3\inf_{(y,p) \in\mathcal{Y} \times\mathcal{P}} \mathcal{L}(y,p)$ and statistical error $\mathcal{E}_{stat}:=4\mathcal{R}_{{\nd}}^{\rm
off} (\mathcal{G}_d, (2b_d)^{-1})+4\mathcal{R}_{n_b}^{\rm
off} (\mathcal{G}_b, (2b_b)^{-1})$. The former arises from restricting the trial spaces of $(y,p)$ from the space $(H_0^1(\Omega)\cap H^2(\Omega))^2$ to $\mathcal{Y}\times\mathcal{P}$, and the latter is due to approximating the integrals with Monte Carlo methods. Throughout, the optimization error is ignored.

\subsubsection{Approximation Error}\label{sssec:approxerr}
We recall a result on the approximation error \cite[Prop. 4.8]{GuhrinRaslan:2021}.
\begin{proposition}\label{prop:approx}
Let $p\geq 1,s,k,d\in \mathbb{N}\cup\{0\}, s\geq k+1$, $\rho$ be the ${\rm logistic}$ or $\tanh$ function, and fix $\mu>0$ small. Then for any $\epsilon>0$, $f\in W^{s,p}([0,1]^d)$ with $\|f\|_{W^{s,p}([0,1]^d)}\leq 1$, there exists an NN $f_{\theta}\in \mathcal{N}_\rho(c\log(d+s),\, c(d,s,p,k)\epsilon^{-\frac{d}{s-k-\mu}},\, c(d,s,p,k)\epsilon^{-2-\frac{2(d/p+d+k+\mu)+d/p+d}{s-k-\mu }})$ such that
$\|f-f_{\theta}\|_{W^{k,p}([0,1]^d)}\leq \epsilon$.
\end{proposition}

Then we have the following bound on the approximation error $\mathcal{E}_{app}$.
\begin{lemma}\label{lem:approx}
Fix a tolerance $\epsilon\in (0,1)$, and $\mu>0$ arbitrarily small. If the optimal state $\bar y\in H^s(\Omega)\cap H_0^1(\Omega)$ and adjoint $\bar p\in H^s(\Omega)\cap H_0^1(\Omega)$, $s\geq 3$, then there exist NNs $( y_{\theta},p_{\sigma})\in \mathcal{Y} \times \mathcal{P}$ with
$\mathcal{Y} = \mathcal{P} =  \mathcal{N}_{\rho}(c\log(d+s),c(d,s)\epsilon^{-\frac{d}{s-2-\mu}},c(d,s)\epsilon^{-\frac{9d+4s}{2(s-2-\mu)}})$
and $c>0$ dependent on $\bf\alpha$, $\lambda$, $s$, $\|\bar{y}\|_{H^s(\Omega)}$ and $\|\bar{p}\|_{H^s(\Omega)}$ such that the approximation error
$\mathcal{E}_{app}\leq c \epsilon^2.$
\end{lemma}
\begin{proof}
For any $(y_\theta,p_\sigma)\in\mathcal{Y}\times\mathcal{P}$,
 since $\mathcal{L}(\bar y,\bar p)=0$, we deduce
\begin{align*}
 &\mathcal{L}(y_{\theta},p_{\sigma})-\mathcal{L}(\bar{y},\bar{p})=  \mathcal{L}(y_{\theta},p_{\sigma})\\
=& \|{\Delta y_{\theta}+f- \lambda^{-1} p_{\sigma} }\|^2_{L^2(\Omega)}+ \alpha_{i} \|{\Delta p_{\sigma}+y_{\theta}-y_d }\|^2_{L^2(\Omega)}+\alpha^y_b\|{y_{\theta}}\|^2_{L^2(\partial\Omega)}+\alpha^p_b \|{p_{\sigma}}\|^2_{L^2(\partial\Omega)}\\
    =& \|{\Delta (y_{\theta}-\bar{y})+ \lambda^{-1}(\bar{p}-p_{\sigma})}\|^2_{L^2(\Omega)}+\alpha_{i} \|{\Delta (p_{\sigma}-\bar{p})+y_{\theta}-\bar{y}}\|^2_{L^2(\Omega)}
    +\alpha^y_b\|{y_{\theta}-\bar{y}}\|^2_{L^2(\partial\Omega)}+\alpha^p_b \|{p_{\sigma}-\bar{p}}\|^2_{L^2(\partial\Omega)}\\
    \leq& c(\boldsymbol\alpha,\lambda) \big(\|{\bar{y}-y_{\theta}}\|^2_{H^2(\Omega)}+\|{\bar{p}-p_{\sigma}}\|^2_{H^2(\Omega)}\big),
\end{align*}
by the trace inequality
$\|{y_{\theta}-\bar{y}}\|_{L^2(\partial\Omega)}\leq c \|y_{\theta}-\bar{y}\|_{H^2(\Omega)}$.
With $c_y=\|\bar y\|_{H^s(\Omega)}$, we have
\begin{align*}
 \inf_{y_{\theta}\in \mathcal{Y}}\|\bar{y}-y_{\theta}\|^2_{H^2(\Omega)}&=c_y^2\inf_{y_{\theta}\in \mathcal{Y}}\Big\|\frac{\bar y}{c_y}-\frac{y_\theta}{c_y}\Big\|^2_{H^2(\Omega)}
 =c_y^2\inf_{y_{\theta}\in \mathcal{Y}}\Big\|\frac{\bar{y}}{c_y}-y_{\theta}\Big\|^2_{H^2(\Omega)}.
\end{align*}
By Proposition \ref{prop:approx}, there exists an NN $y_{\theta}\in \mathcal{Y}= \mathcal{N}_{\rho}(c\log(d+s),c(d,s)\epsilon^{-\frac{d}{s-2-\mu}},c(d,s)\epsilon^{-\frac{9d+4s}{2(s-2-\mu)}})$
such that
$\|\frac{\bar{y}}{c_y}-y_{\theta}\|_{H^2(\Omega)}\leq \epsilon$. Similarly, there exists an NN $p_\sigma\in \mathcal{P}$ with the requisite property. Hence,
\begin{align*}
    \mathcal{E}_{app}&= 3\inf_{(y_{\theta},p_{\sigma})\in \mathcal{Y}\times \mathcal{P}}\mathcal{L}(y_\theta,p_\sigma)\leq c(\boldsymbol{\alpha},\lambda)\inf_{(y_{\theta},p_{\sigma})\in \mathcal{Y}\times \mathcal{P}}(\norm{\bar{y}-y_{\theta}}^2_{H^2(\Omega)}+\norm{\bar{p}-p_{\sigma}}^2_{H^2(\Omega)})
    \leq c \epsilon^2,
\end{align*}
where $c$ depends on $\boldsymbol\alpha$, $\lambda$, $s$ and the $H^s(\Omega)$ norms of $\bar y$ and $\bar p$.
Then the desired assertion follows.
\end{proof}

\begin{remark}
The approximation error $\mathcal{E}_{app}$depends on the regularity of $\bar y$ and $\bar u$. 
If $f,y_d\in H^{s-2}(\Omega)\cap L^{\infty}(\Omega)$, by standard elliptic regularity theory, we have $\bar{y},\bar p\in H^s(\Omega)$, which satisfies the requirements of Lemma \ref{lem:approx}. In the unconstrained case, the optimal control $\bar u$ has the same regularity as the adjoint state $\bar p$, but in the constrained case, we have $\bar u\in H^1(\Omega)$ only, due to the presence of the projection operator $P_U$.
\end{remark}

\subsubsection{Statistical error}

To bound the statistical error via the offset Rademacher complexities, we employ the covering number \cite{AnthonyBartlett:1999}.
\begin{definition}
Let $(\mathcal{M},d)$ be a metric space, and the set $\mathcal{F}\subset \mathcal{M}$. A set $\mathcal{F}_\delta\subset \mathcal{M}$ is called a $\delta$-cover of $\mathcal{F}$ if for each $f\in\mathcal{F}$, there exists an $f_\delta\in \mathcal{F}_\delta$ such that $d(f,f_\delta)\leq \delta$. Moreover, $N(\delta,\mathcal{F},d):=\inf\{|\mathcal{F}_\delta|: \mathcal{F}_\delta \mbox{ is a $\delta$-cover of }  \mathcal{F}\}$ is called the $\delta$-covering number of $\mathcal{F}$.
\end{definition}

We also need Hoeffding's inequality \cite{Hoeffding:1963}.
\begin{lemma}\label{lem:Hoeffding}
Let $X_1$, ..., $X_n$ be independent random variables such that $a_{i}\leq X_{i}\leq b_{i}$ almost surely. Let
$S_n = X_1 + \cdots + X_n$. Then for all $t > 0$,
$\mathbb{P}\left(S_{n}-\mathbb{E}[S_{n}]\geq t\right)\leq \exp(-{\frac {2t^{2}}{\sum _{i=1}^{n}(b_{i}-a_{i})^{2}}}).$
\end{lemma}

Next, we bound the offset Rademacher complexity in terms of the covering number.
\begin{theorem}\label{thm:stat-err0}
Let $\kappa_d$ and $\kappa_b$ be the Lipschitz constants of ${\gd}\in \mathcal{G}_d$ and $g_b\in \mathcal{G}_b$ in $y$ and $p$ in the $W^{2,\infty}(\Omega)$ and $L^\infty(\partial\Omega)$ norms, respectively, i.e.,
\begin{align*}
  |{\gd}(y,p,X)-{\gd}(\tilde y,\tilde p,X)|&\leq \kappa_d(\|y-\tilde y\|_{W^{2,\infty}(\Omega)}+\|p-\tilde p\|_{W^{2,\infty}(\Omega)}),\quad \forall (y,p),(\tilde y,\tilde p)\in \mathcal{Y}\times\mathcal{P}, X\in\Omega,\\
  |g_b(y,p,Y)-g_b(\tilde y,\tilde p,Y)|&\leq \kappa_b(\|y-\tilde y\|_{L^{\infty}(\partial\Omega)}+\|p-\tilde p\|_{L^{\infty}(\partial\Omega)}),\quad \forall (y,p),(\tilde y,\tilde p)\in \mathcal{Y}\times\mathcal{P}, Y\in\partial\Omega.
\end{align*}
Then for any $\delta>0$, the offset Rademacher complexities $\mathcal{R}_n^{\rm off}(\mathcal{G}_d,\beta)$ and $\mathcal{R}_n^{\rm off}(\mathcal{G}_b,\beta)$ are bounded by
\begin{align*}
\mathcal{R}_n^{\rm off}(\mathcal{G}_d,\beta)
\leq & \frac{1+\log N(\delta,\mathcal{F}\times\mathcal{P},W^{2,\infty}(\Omega))}{2n\beta}
+2(1+2b_d\beta)\kappa_d\delta,\\
\mathcal{R}_n^{\rm off}(\mathcal{G}_b,\beta)
\leq & \frac{1+\log N(\delta,\mathcal{F}\times\mathcal{P},L^{\infty}(\partial\Omega))}{2n\beta}
+2(1+2b_b\beta)\kappa_b\delta.
\end{align*}
\end{theorem}
\begin{proof}
We give the proof only for $\mathcal{R}_n^{\rm off}(\mathcal{G}_d,\beta)$, since that for $\mathcal{R}_n^{\rm off}(\mathcal{G}_b,\beta)$ is identical.
Since $\tau=\{\tau_i\}_{i=1}^n$ is a sequence of i.i.d. random variables independent of $\mathbb{X}$, then conditionally on $\mathbb{X}$, we have
\begin{align*}
  & \mathbb{E}_{\tau|\mathbb{X}}\Big[\sup_{(y,p)\in\mathcal{Y}\times\mathcal{P}}\frac{1}{n}\sum_{i=1}^n\tau_i {\gd}(y,p,X_i)-\frac{\beta}{n}\sum_{i=1}^n{\gd}(y,p,X_i)^2|\mathbb{X}\Big]\\
  =& \mathbb{E}_{\tau}\sup_{(y,p)\in\mathcal{Y}\times\mathcal{P}}\Big[\frac{1}{n}\sum_{i=1}^n\tau_i {\gd}(y,p,X_i)-\frac{\beta}{n}\sum_{i=1}^n{\gd}(y,p,X_i)^2\Big].
\end{align*}
Let $\delta>0$ and let $\mathcal{Y}_\delta$ and $\mathcal{P}_\delta$  be a minimal $W^{2,\infty}(\Omega)$ $\delta$-cover of $\mathcal{Y}$ and $\mathcal{P}$, respectively.
For any fixed $(y,p)\in \mathcal{Y}\times \mathcal{P}$, there exists a tuple $(y_\delta,p_\delta)\in\mathcal{Y}_\delta\times\mathcal{P}_\delta$ such that $\|y-y_\delta\|_{W^{2,\infty}(\Omega)}\leq \delta$ and $\|p-p_\delta\|_{W^{2,\infty}(\Omega)}\leq \delta$.
Therefore, by the Lipschitz continuity of ${\gd}$ with respect to $(y,p)$ in the $W^{2,\infty}(\Omega)$ norm,
\begin{align*}
  \frac{1}{n}\sum_{i=1}^n\tau_i {\gd}(y,p,X_i) & \leq  \frac{1}{n}\sum_{i=1}^n\tau_i {\gd}(y_\delta,p_\delta,X_i) + \frac{1}{n}\sum_{i=1}^n|\tau_i||({\gd}(y,p,X_i)-{\gd}(y_\delta,q_\delta,X_i)|\\
  & \leq \frac{1}{n}\sum_{i=1}^n\tau_i {\gd}(y_\delta,p_\delta,X_i) + 2\kappa_d\delta.
\end{align*}
Since $|g_d(y,p,X)|\leq b_d$ and $|g_d(y_\delta,p_\delta,X)|\leq b_d$, and by Lipschitz continuity of ${\gd}(y,p,X)$, we obtain
\begin{align*}
  -\frac{1}{n}\sum_{i=1}^n{\gd}(y,p,X_i)^2 & = - \frac{1}{n}\sum_{i=1}^n{\gd}(y_\delta,p_\delta,X_i)^2 + \frac{1}{n}\sum_{i=1}^n({\gd}(y_\delta,p_\delta,X_i)^2-{\gd}(y,p,X_i)^2)\\
  & \leq - \frac{1}{n}\sum_{i=1}^n{\gd}(y_\delta,p_\delta,X_i)^2 +\frac{ 2b_d}{n}\sum_{i=1}^n|{\gd}(y_\delta,p_\delta,X_i)-{\gd}(y,p,X_i)|\\
  & \leq - \frac{1}{n}\sum_{i=1}^n{\gd}(y_\delta,p_\delta,X_i)^2 + 4b_d\kappa_d\delta.
\end{align*}
Hence, it follows that
\begin{align}
  &\mathbb{E}_{\tau}\sup_{(y,p)\in\mathcal{Y}\times\mathcal{P}}\Big[\frac{1}{n}\sum_{i=1}^n\tau_i {\gd}(y,p,X_i)-\frac{\beta}{n}\sum_{i=1}^n{\gd}(y,p,X_i)^2\Big]\nonumber\\
\leq & \mathbb{E}_\tau\max_{(y_\delta,p_\delta)\in\mathcal{F}_\delta\times\mathcal{P}_\delta} \Big[\frac{1}{n}\sum_{i=1}^n\tau_i {\gd}(y_\delta,p_\delta,X_i)-\frac{\beta}{n}\sum_{i=1}^n{\gd}(y_\delta,p_\delta,X_i)^2\Big]\label{eqn:est-Rad-off-basic}
+2(1+2b_d\beta)\kappa_d\delta.
\end{align}
Since $\{\tau_ig_d(y_\delta,p_\delta;X_i)\}_{i=1}^n$ are independent random variables conditioning on $\mathbb{X}$,
\begin{align*}
  \mathbb{E}_\tau [\tau_i{\gd}(y_\delta,p_\delta,X_i)]&=0, \quad i=1,\ldots,n,\\
   -{\gd}(y_\delta,p_\delta,X_i) \leq \tau_i{\gd}(y_\delta,p_\delta,X_i)&\leq {\gd}(y_\delta,p_\delta,X_i), \quad i=1,\ldots,n.
\end{align*}
By Hoeffding's inequality from Lemma \ref{lem:Hoeffding}, we deduce that for any $(y_\delta,p_\delta)\in\mathcal{Y}_\delta\times \mathcal{P}_\delta$ and $\xi>0$,
\begin{align*}
  &\mathbb{P}_\tau \Big\{\frac{1}{n}\sum_{i=1}^n\tau_i{\gd}(y_\delta,p_\delta,X_i)>\xi + \frac{\beta}{n}\sum_{i=1}^n{\gd}(y_\delta,p_\delta,X_i)^2\Big\} \\
   \leq& \exp\Big(-\frac{(n\xi+\beta\sum_{i=1}^n{\gd}(y_\delta,p_\delta,X_i)^2)^2}{2\sum_{i=1}^n{\gd}(y_\delta,p_\delta,X_i)^2}\Big)
  \leq \exp(-2\beta n\xi),
\end{align*}
where the last step follows from the elementary inequality $\frac{(a+y)^2}{y}\geq \frac{(a+a)^2}{a}=4a$, for any $y\in\mathbb{R}_+$.
Therefore, we may bound the tail probability by
\begin{align*}
& \mathbb{P}_\tau\left\{\max_{(y_\delta,p_\delta)\in\mathcal{Y}_\delta\times\mathcal{P}_\delta}\Big( \frac{1}{n}\sum_{i=1}^n\tau_i{\gd}(y_\delta,p_\delta,X_i)-\frac{\beta}{n}\sum_{i=1}^n{\gd}(y_\delta,p_\delta,X_i)^2\Big)>\xi \right\}\\
\leq & N(\delta,\mathcal{F}\times\mathcal{P},W^{2,\infty}(\Omega))
\max_{(y_\delta,p_\delta)\in\mathcal{Y}_\delta\times\mathcal{P}_\delta}\mathbb{P}_\tau \Bigl\{\frac{1}{n}\sum_{i=1}^n\tau_i{\gd}(y_\delta,p_\delta,X_i)>\xi + \frac{\beta}{n}\sum_{i=1}^n{\gd}(y_\delta,p_\delta,X_i)^2\Bigr\}\\
\leq &N(\delta,\mathcal{F}\times \mathcal{P},W^{2,\infty}(\Omega))\exp(-2\beta n\xi).
\end{align*}
Hence, for any $a>0$, we have
\begin{align*}
 & \mathbb{E}_\tau \Bigl[\max_{(y_\delta,p_\delta)\in\mathcal{Y}_\delta\times\mathcal{P}_\delta}\Big( \frac{1}{n}\sum_{i=1}^n\tau_i{\gd}(y_\delta,p_\delta,X_i)-\frac{\beta}{n}\sum_{i=1}^n{\gd}(y_\delta,p_\delta,X_i)^2\Big)\Bigr]\\
 \leq & \int_0^\infty \mathbb{P}_\tau\left\{\max_{(y_\delta,p_\delta)\in\mathcal{Y}_\delta\times\mathcal{P}_\delta}\Big( \frac{1}{n}\sum_{i=1}^n\tau_i{\gd}(y_\delta,p_\delta,X_i)-\frac{\beta}{n}\sum_{i=1}^n{\gd}(y_\delta,p_\delta,X_i)^2\Big)>\xi \right\}{\rm d}\xi\\
 \leq & a + \int_a^\infty N(\delta,\mathcal{F}\times\mathcal{P},W^{2,\infty}(\Omega))\exp(-2\beta n\xi){\rm d}\xi
 \leq  a + \frac{N(\delta,\mathcal{F}\times\mathcal{P},W^{2,\infty}(\Omega))}{2\beta n}\exp(-2\beta n a).
\end{align*}
Setting $a= \frac{\log N(\delta,\mathcal{F}\times\mathcal{P},W^{2,\infty}(\Omega))}{2\beta n}$ leads to
\begin{align*}
  &\mathbb{E}_\tau \Big[\max_{(y_\delta,p_\delta)\in\mathcal{Y}_\delta\times\mathcal{P}_\delta}\Big( \frac{1}{n}\sum_{i=1}^n\tau_i{\gd}(y_\delta,p_\delta,X_i)-\frac{\beta}{n}\sum_{i=1}^n{\gd}(y_\delta,p_\delta,X_i)^2\Big]
  \leq \frac{1+\log N(\delta,\mathcal{F}\times \mathcal{P},W^{2,\infty}(\Omega))}{2\beta n}.
\end{align*}
Combining this inequality with the estimate \eqref{eqn:est-Rad-off-basic} yields the desired assertion.
\end{proof}

The next result gives the statistical error, by exploiting the bound on the NNs in the appendix \ref{sec:techestm}.
\begin{theorem}
Let $\mathcal{Y}=\mathcal{P}=\mathcal{N}_\rho(L,\boldsymbol{n}_L,R)$, with depth $L$, $\boldsymbol{n}_L$ nonzero NN parameters and maximum bound $R$. Then the statistical error $\mathcal{E}_{stat}$ is bounded by
\begin{align*}
  \mathcal{E}_{stat} \leq \frac{cL\boldsymbol{n}_L^{4L+1}R^{4L}(L \log R + L\log \boldsymbol{n}_L + \log \nd)}{\nd}+\frac{c\boldsymbol{n}_L^3R^2(L\log R + L\log \boldsymbol{n}_L + \log n_b)}{n_b},
\end{align*}
where the constant $c$ depends on $\boldsymbol\alpha$, $\lambda$, $d$, $\|f\|_{L^\infty(\Omega)}$ and $\|y_d\|_{L^\infty(\Omega)}$ at most polynomially.
\end{theorem}
\begin{proof}
By Lemmas \ref{lem:est-f} and \ref{lem:est-f''-bdd}, for any $v\in \mathcal{Y}$,
$\|v\|_{C(\overline{\Omega})}\leq \boldsymbol{n}_LR $ and $\|\Delta v\|_{L^\infty(\Omega)} \leq dL\boldsymbol{n}_L^{2L} R^{2L}$, which directly give
$b_d=cL\boldsymbol{n}_L^{4L}R^{4L}$ and $b_b=c(\boldsymbol{n}_LR)^2$.
Next, for any $(y,p),(\tilde y,\tilde p)\in\mathcal{Y}\times\mathcal{P}$, we have
\begin{align*}
|{\gd}(y,p,X)-{\gd}(\tilde y,\tilde p,X)|
 \leq &   \kappa_d (\|y -\tilde y\|_{W^{2,\infty}(\Omega)} +\|p-\tilde p\|_{W^{2,\infty}(\Omega)}), \quad\mbox{with }\kappa_d=c\boldsymbol{n}_L^{2L}R^{2L},\\
|g_b(y,p,Y)-g_b(\tilde y,\tilde p,Y)|
  \leq & \kappa_b(\|y-\tilde y\|_{L^\infty(\partial\Omega)} + \|p-\tilde p\|_{L^\infty(\partial\Omega)}),\quad \mbox{with }\kappa_b=c\boldsymbol{n}_LR.
\end{align*}
Next we bound $N(\delta,\mathcal{Y},W^{2,\infty}(\Omega))$.  By Lemmas \ref{lem:est-f} and \ref{lem:est-f''-Lip}, for any
$v_\theta,v_{\tilde\theta}\in\mathcal{Y}$,
\begin{align*}
\|v_\theta-v_{\tilde\theta}\|_{C(\overline{\Omega})}&\leq \boldsymbol{n}_L^{L}R^{L-1}\|\theta-\tilde \theta\|_{\ell^2},\\
\|\partial_{x_ix_i}^2 v_\theta-\partial_{x_ix_i}^2 v_{\tilde\theta}\|_{C(\overline{\Omega})} &\leq4L^2\bs{n}_{L}^{3L-2}R^{3L-3}\|{\theta -\tilde{\theta}}\|_{\ell^2},\quad i=1,\ldots,d.
\end{align*}
Hence, with $\Lambda_d=cL^2\bs{n}_{L}^{3L-2}R^{3L-3}$ and $\Lambda_b=\boldsymbol{n}_L^{L}R^{L}$,
\begin{equation}\label{eqn:NN-Lip}
  \|v_\theta-v_{\tilde\theta}\|_{W^{2,\infty}(\Omega)}\leq \Lambda_d \|{\theta -\tilde{\theta}}\|_{\ell^2}\quad\mbox{and}\quad
   \|v_\theta-v_{\tilde\theta}\|_{L^\infty({\partial\Omega})}\leq \Lambda_b\|\theta-\tilde \theta\|_{\ell^2},\quad \forall v_\theta,v_{\tilde\theta}\in\mathcal{Y}.
\end{equation}
Note that for any $m \in \mathbb{N }$, $r \in [1, \infty)$, $\epsilon \in  (0,1)$,
and $ B_r := \{x\in\mathbb{R}^m:\ \|x\|_{\ell^2}\leq r\}$, then by a simple counting
argument (see, e.g., \cite[Proposition 5]{CuckerSmale:2002} or \cite[Lemma 5.5]{JiaoLai:2021error}), we have
$\log N(\epsilon,B_r,\|\cdot\|_{\ell^2})\leq m\log (2r\sqrt{m}\epsilon^{-1})$.
This, Lipschitz continuity of NN functions in \eqref{eqn:NN-Lip} and the estimate $\|\theta\|_{\ell^2}\leq \sqrt{\boldsymbol{n}_L}\|\theta\|_{\ell^\infty}\leq \sqrt{\boldsymbol{n}_L}R$ imply
\begin{align*}
\log N(\delta,\mathcal{Y},W^{2,\infty}(\Omega))&\leq \log N(\Lambda_d^{-1}\delta,\mathcal{N}_Y,\|\cdot\|_{\ell^2})\leq \boldsymbol{n}_L\log(2R\boldsymbol{n}_L\Lambda_d\delta^{-1}),\\
\log N(\delta,\mathcal{Y},L^{\infty}(\partial\Omega))&\leq \log N(\Lambda_b^{-1}\delta,\mathcal{N}_Y,\|\cdot\|_{\ell^2})\leq \boldsymbol{n}_L\log(2R\boldsymbol{n}_L\Lambda_b\delta^{-1}),
\end{align*}
where $\mathcal{N}_Y$ denotes the parameter space for $\mathcal{Y}$.
These estimates and Theorem \ref{thm:stat-err0}
with $\beta=(2b_d)^{-1}$ yield
\begin{align*}
\mathcal{R}_n^{\rm off}\big(\mathcal{G}_d,(2b_d)^{-1}\big) &\leq \frac{b_d(1+\log N(\delta,\mathcal{F}\times\mathcal{P},W^{2,\infty}(\Omega)))}{n}
+4\kappa_d\delta\\
  & \leq \frac{cL\boldsymbol{n}_L^{4L+1}R^{4L}\log(2R\boldsymbol{n}_L\Lambda_d\delta^{-1})}{n}+ c\boldsymbol{n}_L^{2L}R^{2L}\delta.
\end{align*}
An analogous bound on $\mathcal{R}_n^{\rm off}\big(\mathcal{G}_b,(2b_b)^{-1}\big)$ holds.
Setting $\delta=1/n$, substituting $\Lambda_d$ and $\Lambda_b$ and simplifying the resulting expressions yield the desired estimate.
\end{proof}

\subsubsection{Final error estimate}
Now we state the error of the approximation $(y_{\theta^*},p_{\sigma^*},u_{\sigma^*})$ given by OSNN. Thus, with the parameters in the loss $\widehat{\mathcal{L}}(y_\theta,p_\sigma)$ chosen suitably, OSNN can yield an accurate approximation.
\begin{theorem}\label{thm:conv:main}
Let the tuple $(\bar y,\bar p,\bar u)$ solve the optimality system \eqref{eqn:opt-unconstrained}
/ \eqref{eqn:coupled-cons-first} such that $\bar y\in H^s(\Omega)\cap H_0^1(\Omega)$ and $\bar p\in H^s(\Omega)\cap H_0^1(\Omega)$ with $s\geq 3$, and $(y_{\theta^*},p_{\sigma^*},u_{\sigma^*})$ be the OSNN approximation. Fix the
tolerance $\epsilon>0$, and take $\mathcal{Y}=\mathcal{P}=\mathcal{N}_\rho(c\log(d+s),c(d,s)\epsilon^{-\frac{d}{s-2-\mu}},
c(d,s)\epsilon^{-\frac{9d+4s}{2(s-2-\mu)}})$. Then by choosing ${\nd}=c(d,s)\epsilon^{-\frac{c(d+s)\log(d+s)}{s-2-\mu}}$ and $n_b=c(d,s)\epsilon^{-\frac{4(3d+s)}{s-2-\mu}-2}$ sampling points
in $\Omega$ and on $\partial\Omega$, there holds
\begin{align*}
  \mathbb{E}_{\mathbb{X},\mathbb{Y}}\big[\|\bar y-y_{\theta^*}\|_{L^2(\Omega)} ^2 + \|\bar p - p_{\sigma^*}\|_{L^2(\Omega)}^2 + \|\bar u - u_{\sigma^*}\|_{L^2(\Omega)}^2\big] \le c\epsilon^2,
\end{align*}
where the constant $c$ depends on $\boldsymbol\alpha$, $\lambda$, $d$, $s$, $\|f\|_{L^\infty(\Omega)}$, $\|y_d\|_{L^\infty(\Omega)}$, $\|\bar y\|_{H^s(\Omega)}$, and $\|\bar p\|_{H^s(\Omega)}$.
\end{theorem}
\begin{proof}
The fundamental estimates in Lemmas \ref{lem:fund-est1} and \ref{lem:fund-est2} imply
\begin{align*}
  \mathbb{E}_{\mathbb{X},\mathbb{Y}}[\|\bar y-y_{\theta^*}\|_{L^2(\Omega)} ^2 + \|\bar p - p_{\sigma^*}\|_{L^2(\Omega)}^2 + \|\bar u - u_{\sigma^*}\|_{L^2(\Omega)}^2] \leq c(\boldsymbol\alpha,\gamma)\mathbb{E}_{\mathbb{X},\mathbb{Y}}[\mathcal{L}(y_{\theta^*},p_{\sigma^*})].
\end{align*}
By the error decomposition in Theorem \ref{thm:err-decomp} and Lemma \ref{lem:approx} (with the given choice of $\mathcal{Y}$ and $\mathcal{P}$),
\begin{equation*}
  \mathbb{E}_{\mathbb{X},\mathbb{Y}}[\mathcal{L}(y_{\theta^*},p_{\sigma^*})] \leq c\epsilon^2 + \frac{cL\boldsymbol{n}_L^{4L+1}R^{4L}(L \log R + L\log \boldsymbol{n}_L + \log n_d)}{n_d} + \frac{c\boldsymbol{n}_L^3R^2(L\log R + L\log \boldsymbol{n}_L + \log n_b)}{n_b}.
\end{equation*}
Since $L=c\log(d+s)$, $\boldsymbol{n}_L=c(d,s)\epsilon^{-\frac{d}{s-2-\mu}}$ and $R=c(d,s)\epsilon^{-\frac{9d+4s}{2(s-2-\mu)}} $, the choices of the numbers ${\nd}$ and $n_b$ imply that the statistical error $\mathcal{E}_{stat}$ is also bounded by $c\epsilon^2$. This completes the proof of the theorem.
\end{proof}

\begin{remark}
The error bound on the OSNN approximation in Theorem \ref{thm:conv:main} is given in terms of the $L^2(\Omega)$ norm, due to the presence of the consistency error induced by the $L^2(\partial\Omega)$ boundary penalty, which precludes obtaining $H^1(\Omega)$ error estimates. In contrast, for the conventional Galerkin FEM discretization, one can obtain error estimates in $H^1(\Omega)$; see \cite{Falk:1973} for an early work, and \cite{ManzoniQuarteroniSalsa:2021} for many further results. See also \cite{HoffWendland:2021} for error estimates of discretizing optimal control problems with matrix-valued positive definite kernels.
\end{remark}

\section{Numerical experiments and discussions}\label{sec:numexp}
Now we present some numerical experiments to illustrate OSNN, and compare it with ALM and AONN, cf. Section \ref{ssec:existing}. All the examples are academic but not directly amenable with the standard Galerkin FEM, and we leave a more thorough experimental evaluation, including more advanced training algorithms and more complex problem settings to future work. Throughout, we fix the weight $\alpha_i=\lambda^{-1}$ in OSNN and set $\alpha^y_b=\alpha^p_b/\alpha_i$ (so there is only one hyper-parameter $\alpha^y_b$ in OSNN). This heuristic rule works reasonably well for OSNN. AONN and ALM involve more hyper-parameters. The boundary weight $\alpha_b$ of the PINN residual is determined by grid search \cite{MowlaviNabi:2023}: we apply PINN to solve a direct problem with analytic solution, with a range of  $\alpha_b$ values, and select the one $\alpha_b^*$ attaining the smallest validation error for use in the state and adjoint solvers. Other hyper-parameters, e.g., penalty factor $\mu$ in ALM, and step size $s$ and the number $K$ of gradient descent steps in AONN are all determined manually; see Table \ref{tab:hyper} for the values. The scalars $N_{\rm osnn}$ and $N_{\rm alm}$ denote the total number of iterations for OSNN or the number of iterations for each subproblem (ALM). For AONN, each PDE solve and the gradient descent step described in \eqref{eqn:aonn-gd} employs $N_{\rm aonn}$ iterations.
The function $\rho$ is taken to be $\tanh$. The NN parameters are initialized by the Xavier scheme (with  zero bias) via Flax (\url{https://github.com/google/flax}), and the multipliers in ALM are initialized to zero. To measure the accuracy of the approximate state $y^*$, we employ the relative $L^2(\Omega)$ error
$e(y^*)=\|y^*-\bar y\|_{L^2(\Omega)}/\|\bar y\|_{L^2(\Omega)}$,
where $\bar y$ denotes the exact state, and similarly for the control $u^*$. The relative errors are computed using a test set different from the training one. 
Throughout, the NN has $4$ hidden layers, each having 80 neurons. We employ the Adam algorithm \cite{KingmaBa:2015} provided by the Optax library \texttt{optax.adam} (version 2.0.0) from JAX (\url{https://github.com/google/jax}), using single floating-point precision. All the experiments are conducted on an NVIDIA RTX A100 GPU. The code for reproducing the experiments will be made publicly available at Github.

\begin{table}[hbt!]
\centering \setlength{\tabcolsep}{6pt}
\caption{Hyper-parameter settings for the Examples.}\label{tab:hyper} 
\begin{tabular}{c|ccc|c|c}
\toprule
Example & $\alpha_b$ &  $(s,K)_{\rm aonn}$  & $(\mu,K)_{\rm alm}$& Parameter&Value\\ 
\midrule
\ref{exam:hidim-unconstraint} & 200 & (20, 100) & (0.02, 20)&$N_{\rm osnn}$&6e4\\ 
\ref{exam:4d-constraint} & 100 & (10, 100) & (0.01, 20)& $N_{\rm aonn}$&6e2\\
\ref{exam:semilinear} & 100 & (20, 100) & (0.01, 10) &$N_{\rm alm}$&6e4\\
 \bottomrule
\end{tabular}
\end{table}

First, we consider the linear elliptic control problem:
\begin{equation*}
\begin{aligned}
\min_{y,u}\Big\{J(y,u)=\frac{1}{2}\|y-y_d\|^2_{L^2(\Omega)}+\frac{\lambda}{2}\|u\|^2_{L^2(\Omega)}\Big\},\quad
\text{subject to }\left\{\begin{aligned}-\Delta y&=f+u,&&\text{in }\Omega,\\
y&=g,&& \text{on }\partial\Omega,
\end{aligned}\right.
\end{aligned}
\end{equation*}
where $f$ and $g$ are the known problem data and $\lambda$ is the penalty weight.

\begin{example}\label{exam:hidim-unconstraint}
The domain $\Omega=(0,1)^4$, and $\lambda=0.01$. The data $f\equiv0$, $y_d(x)=(1+16\lambda\pi^4)\prod_{i=1}^4 \sin(\pi x_i)$, and $g=0$. The exact state $\bar y$ and control $\bar u$ are given by
$\bar y({x})=\prod_{i=1}^4 \sin(\pi x_i)$, and $\bar u({x})=4\pi^2 \prod_{i=1}^4 \sin(\pi x_i)$.
\end{example}

To form the empirical loss $\widehat{\mathcal{L}}(y_\theta,p_\sigma)$, we employ 60000 points in $\Omega$ and 5000 points on $\partial\Omega$, drawn i.i.d. from $U(\Omega)$ and $U(\partial\Omega)$, respectively. These numbers are determined in a trial-and-error manner. We employ 60000 Adam iterations for OSNN with a learning rate 1e-3, halved after every 6000 iterations; and 2000 Adam iterations for each PDE solve in AONN, with a learning rate $10^{-3}$ for the first subproblem (primal PDE, adjoint PDE and gradient descent), and then switch to $10^{-4}$ for the remaining sub-problems. For ALM, we employ 3000 iterations for each subproblem, with learning rate $10^{-3}$ for the first and also halved after every 6000 iterations. 
The relative errors of the approximate state $y^*$ and control $u^*$ are presented in Table \ref{tab:exam}(a).
The computing time (in second) is 4.28e3, 5.84e3 and 4.47e3 for OSNN, AONN, and ALM, respectively. These numbers should be interpreted only indicatively due to the presence of other influencing factors. Each subproblem in AONN does not require accurate resolution in the early stage: it takes only 600 iterations for each PDE solve, which is far less than that for OSNN and ALM. Table \ref{tab:exam}(a) shows that the choice can ensure accurate solutions with 200 gradient descent steps. OSNN and AONN take comparable computing time, but the approximate control of AONN converges much slower. Nonetheless, all three methods require many iterations to reach satisfactory results, and it is imperative to develop customized algorithms for these solvers and to further improve their accuracy. This phenomenon has been widely reported for standard neural PDE solvers \cite{Krishnapriyan:2021,WangYu:2022,Borovykh:2022}. 

\begin{table}[htp!]
  \centering
    \caption{Numerical results for the examples.\label{tab:exam}}
    \begin{tabular}{|c|ccc|c|ccc|}
    \multicolumn{4}{c}{(a) Example \ref{exam:hidim-unconstraint}}&\multicolumn{4}{c}{(b) Example \ref{exam:4d-constraint}}\\
    \cmidrule(lr){1-4} \cmidrule(lr){5-8}
        &  $e(y^*)$ & $e(u^*)$ & $J(u^*,y^*)$ & &  $e(y^*)$ & $e(u^*)$ & $J(u^*,y^*)$ \\
    \midrule
        \text{OSNN} &  8.24e-3 & 1.01e-2  & 8.088 & \text{OSNN}  & 2.12e-3 & 1.67e-1 & 2.27e-4 \\
        \text{AONN} & 8.28e-3 & 1.46e-2  & 8.038 & \text{AONN} & 3.07e-3 & 1.32e-1 & 3.26e-4 \\
        \text{ALM} & 1.03e-1 & 3.25e-2  & 8.074 & \text{ALM} & 2.1e-3 & 2.08e-1 & 2.47e-3 \\
    \midrule
    \multicolumn{4}{c}{(c) Example \ref{exam:semilinear}} & \multicolumn{4}{c}{(d) Example \ref{exam:highdimcpinn}} \\
    \cmidrule(lr){1-4} \cmidrule(lr){5-8}
           & $e(y^*)$ & $e(u^*)$ & $J(u^*,y^*)$ & PDE  & $e(y^*)$ & $e(u^*)$ & $J(u^*,y^*)$ \\
         \midrule
          \text{OSNN}  & 4.45e-4 & 3.36e-2 & 1.194 & LU  & 3.96e-2 & 6.05e-2 & 2.652e-1 \\
          \text{AONN} & 4.41e-4 & 4.26e-2 & 1.201 & LC & 5.78e-2 & 9.16e-2 & 2.619e-1 \\
          \text{ALM} & 9.11e-3 & 1.34e-1 & 1.231 & Semi & 3.96e-2 & 6.45e-2 & 2.321e-3  \\
    \bottomrule
    \end{tabular}
\end{table}

\begin{figure}[hbt!]
    \centering
    \setlength{\tabcolsep}{0pt}
    \begin{tabular}{cccc}
             \includegraphics[width=0.24\textwidth]{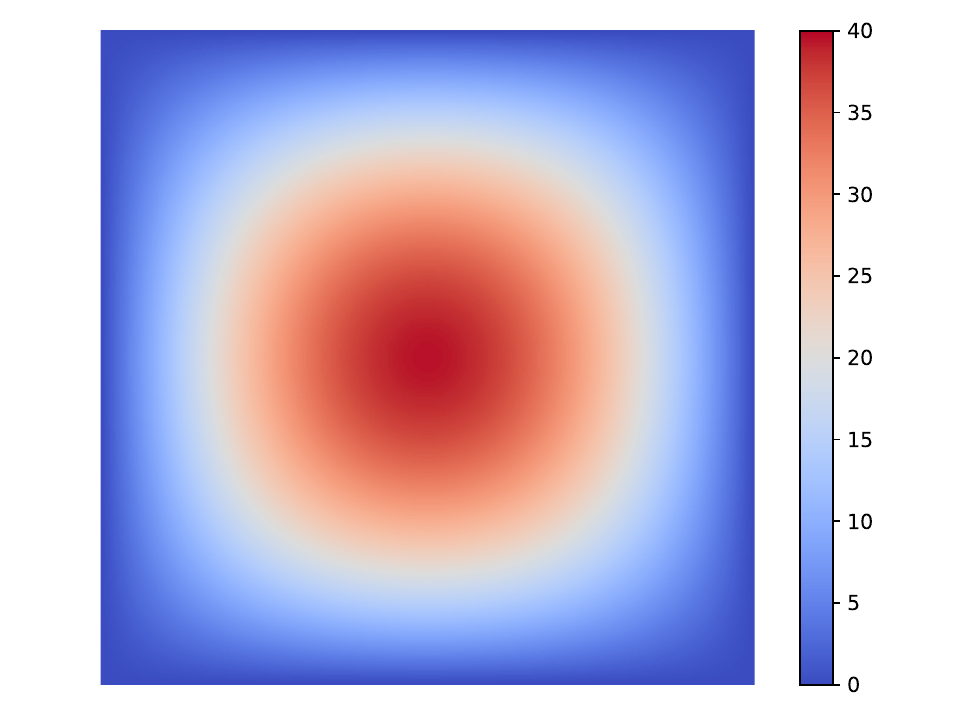}
            &\includegraphics[width=0.24\textwidth]{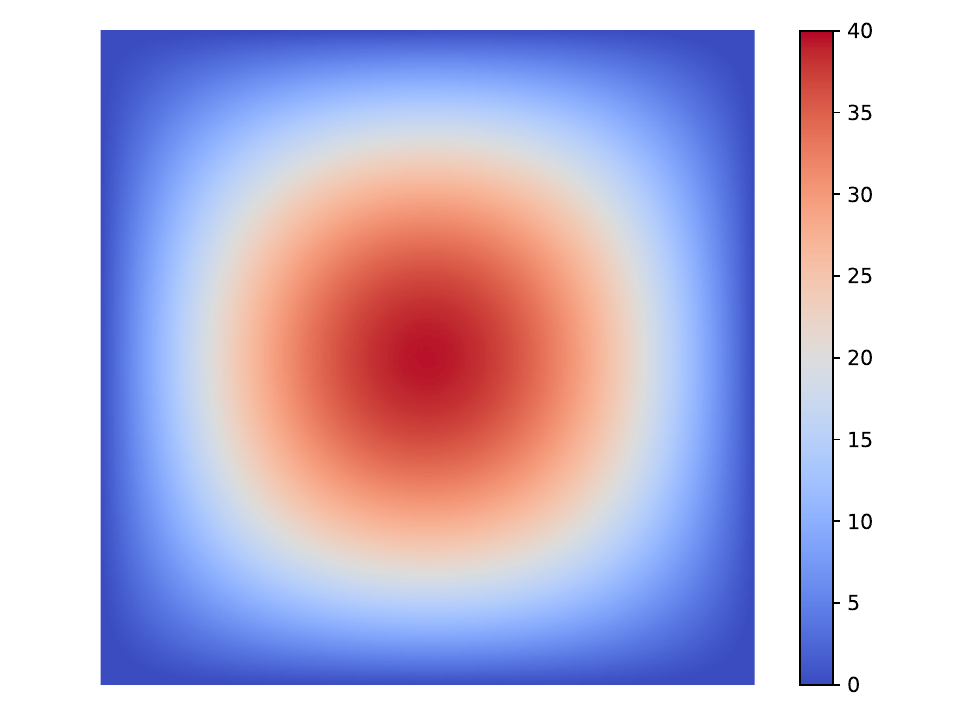}
            &\includegraphics[width=0.24\textwidth]{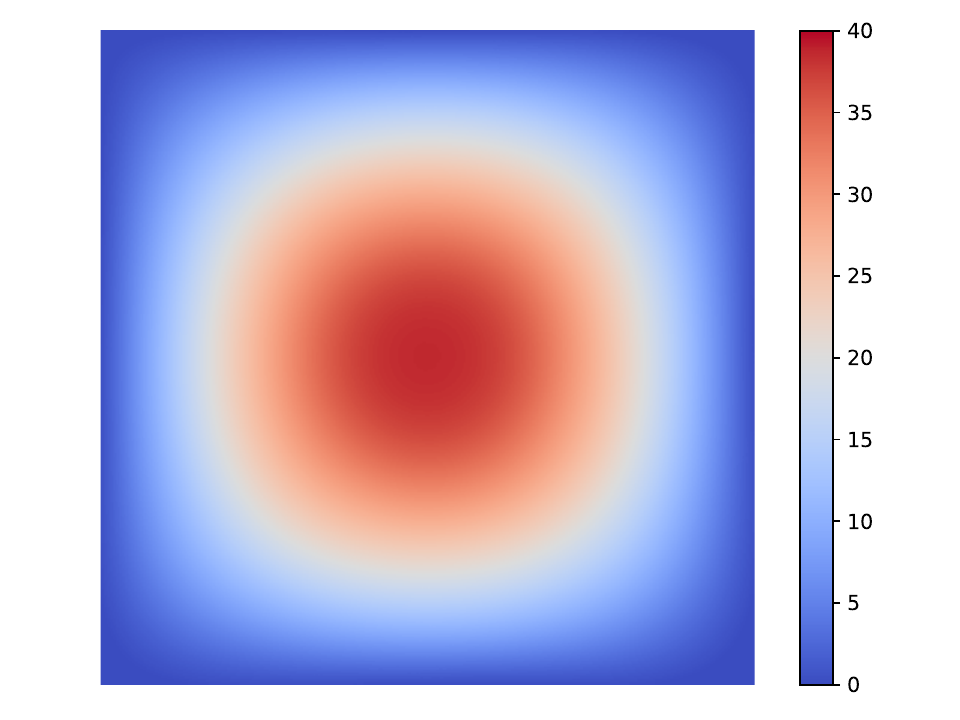}
            &\includegraphics[width=0.24\textwidth]{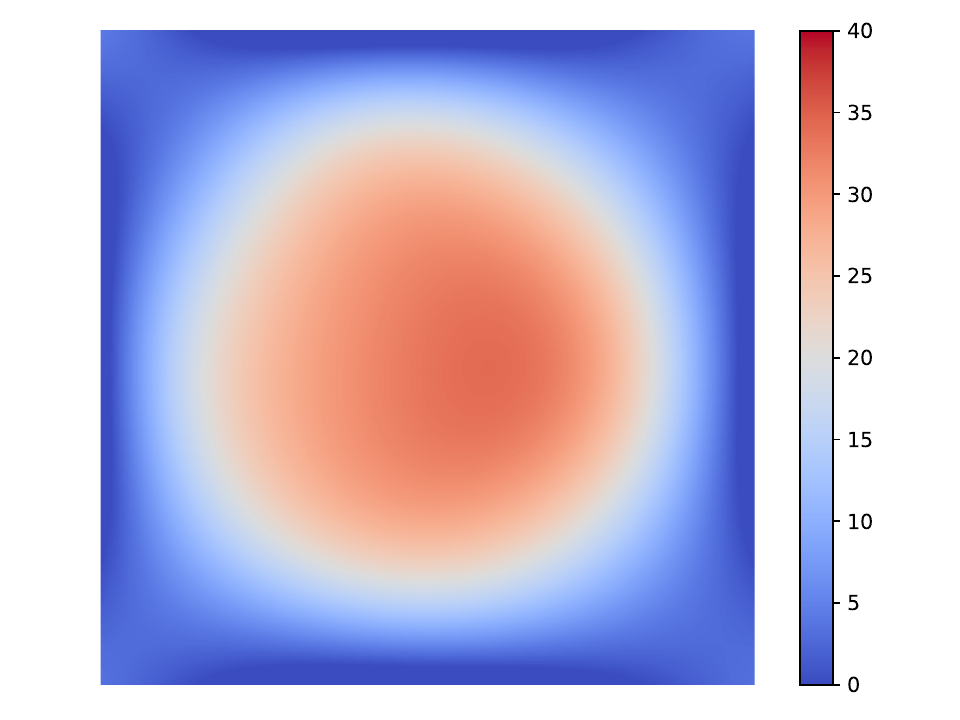} \\
            &\includegraphics[width=0.24\textwidth]{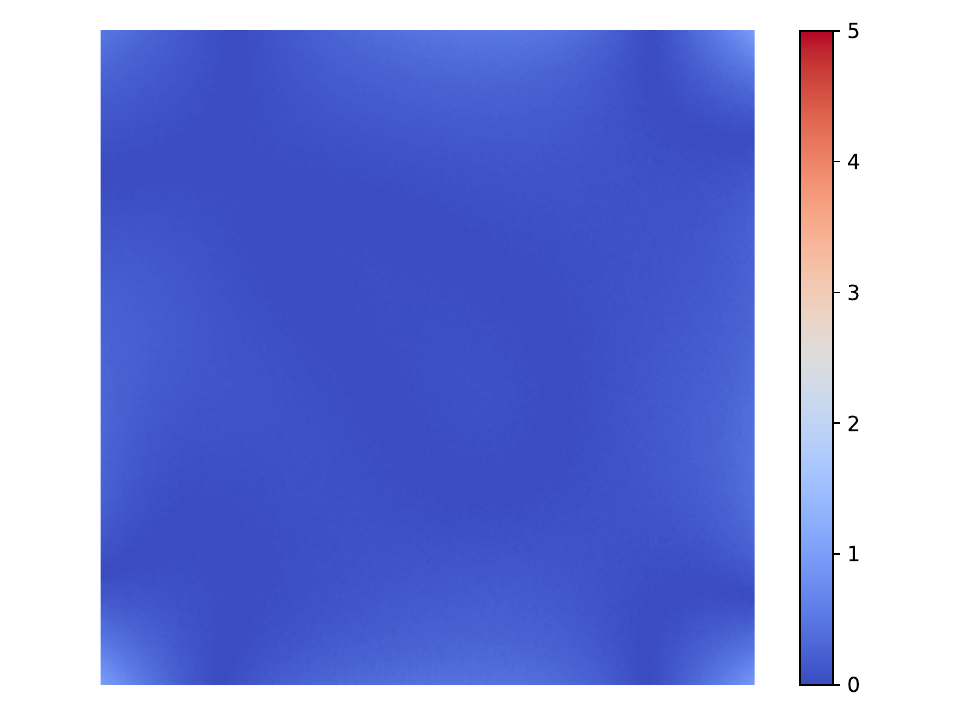}
            &\includegraphics[width=0.24\textwidth]{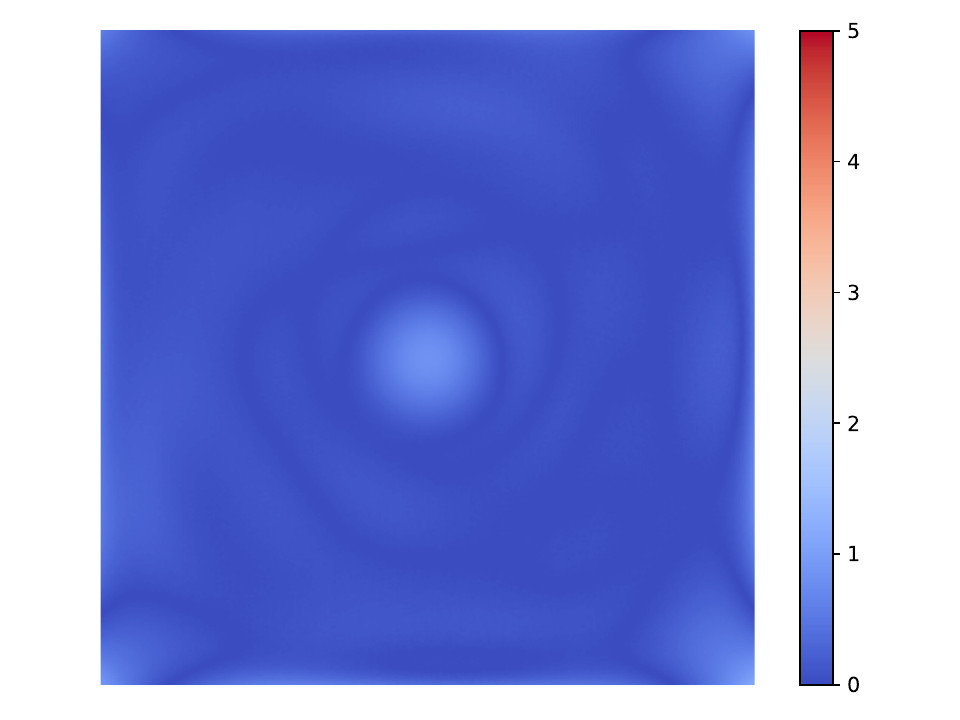}
            &\includegraphics[width=0.24\textwidth]{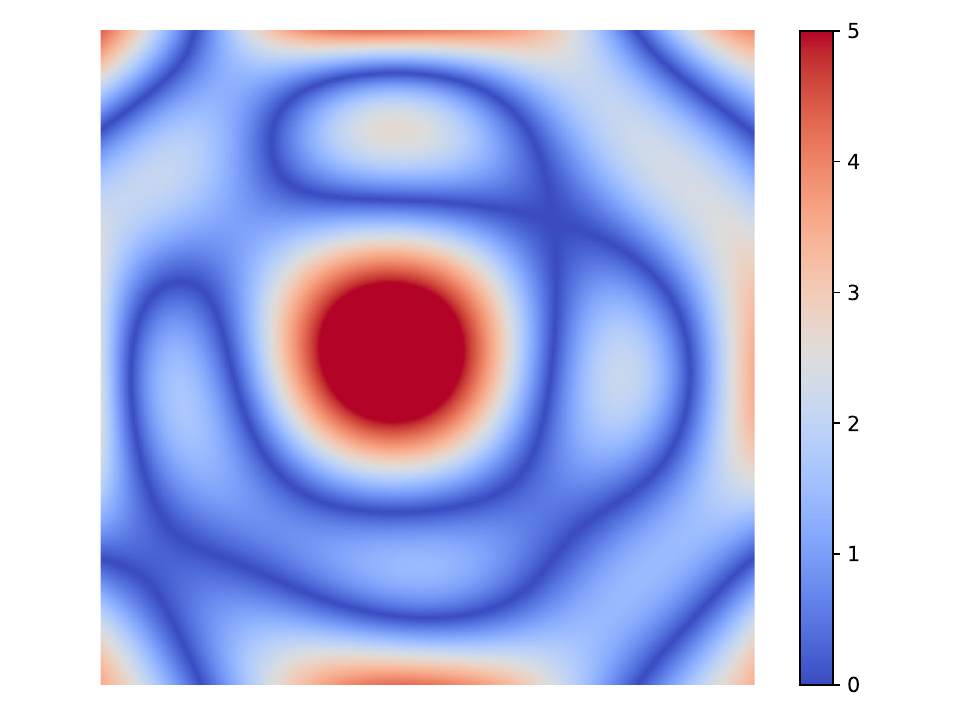}\\
             \includegraphics[width=0.24\textwidth]{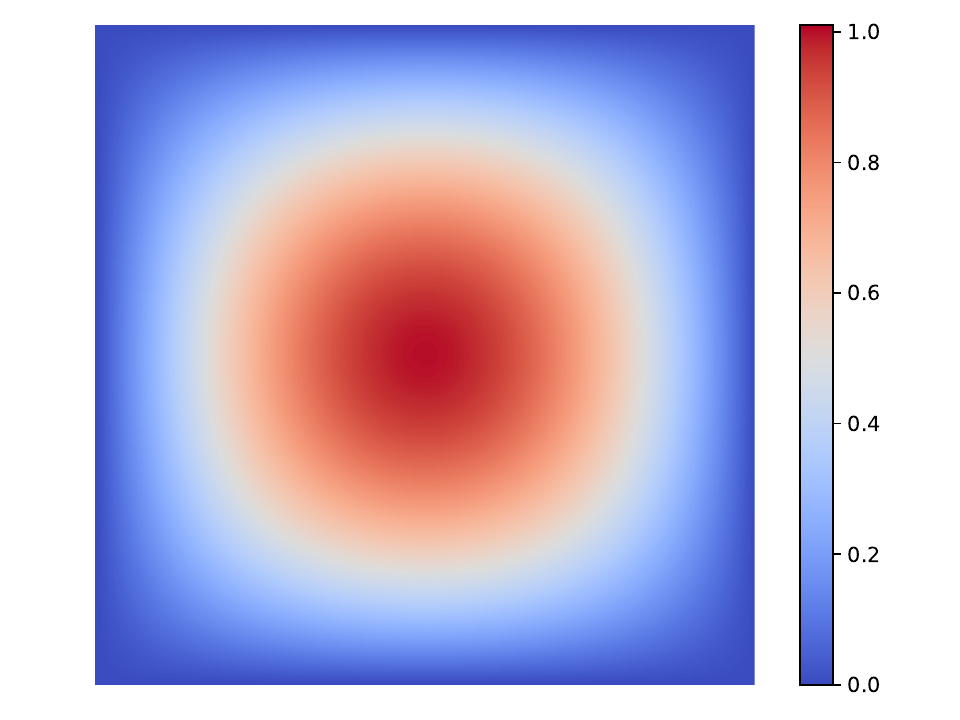}
            &\includegraphics[width=0.24\textwidth]{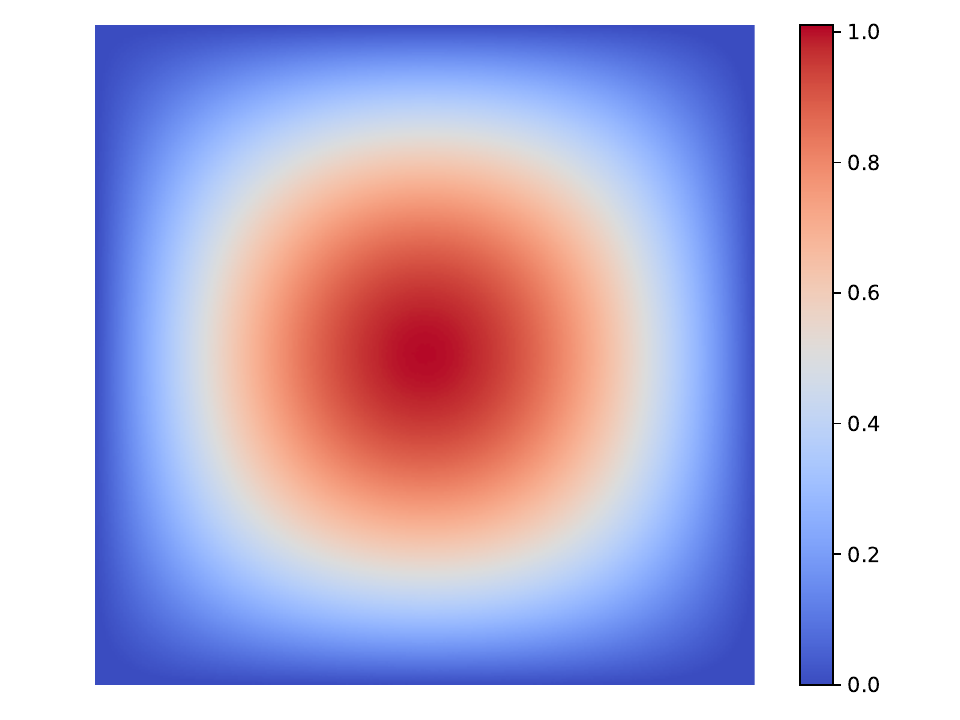}
            &\includegraphics[width=0.24\textwidth]{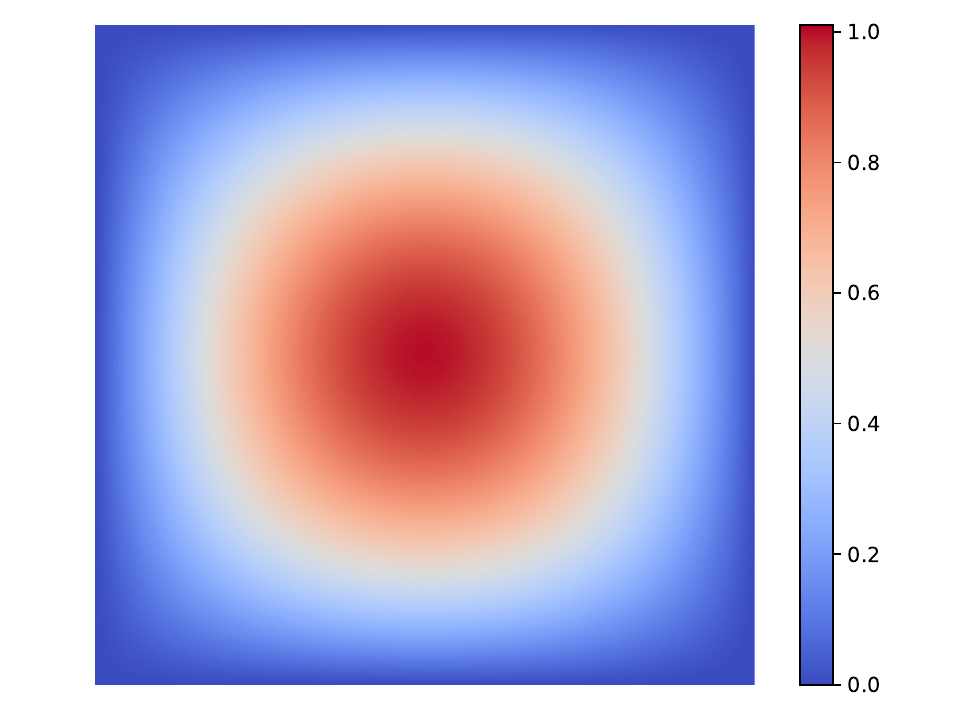}
            &\includegraphics[width=0.24\textwidth]{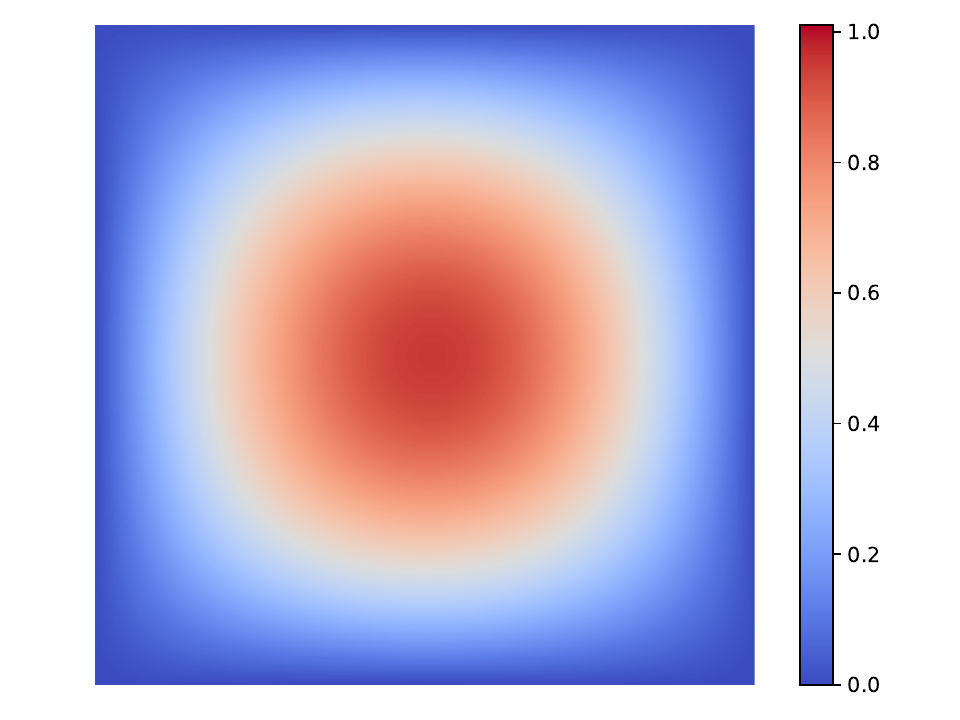}\\
            (a) exact & (b) OSNN & (c) AONN  & (d) ALM
    \end{tabular}
    \caption{The approximate optimal control $u^*$ (top), pointwise error $|u^*-\bar u|$ (mid), and approximate state $y^*$ (bottom) obtained by three NN-based methods for Example \ref{exam:hidim-unconstraint}, cross section at $x_3=x_4=0.5$.}
    \label{fig:hidim-unconstraint}
\end{figure}

\begin{figure}[hbt!]
    \centering
    \setlength{\tabcolsep}{0pt}
    \begin{tabular}{ccc}
       \includegraphics[width=0.33\textwidth]{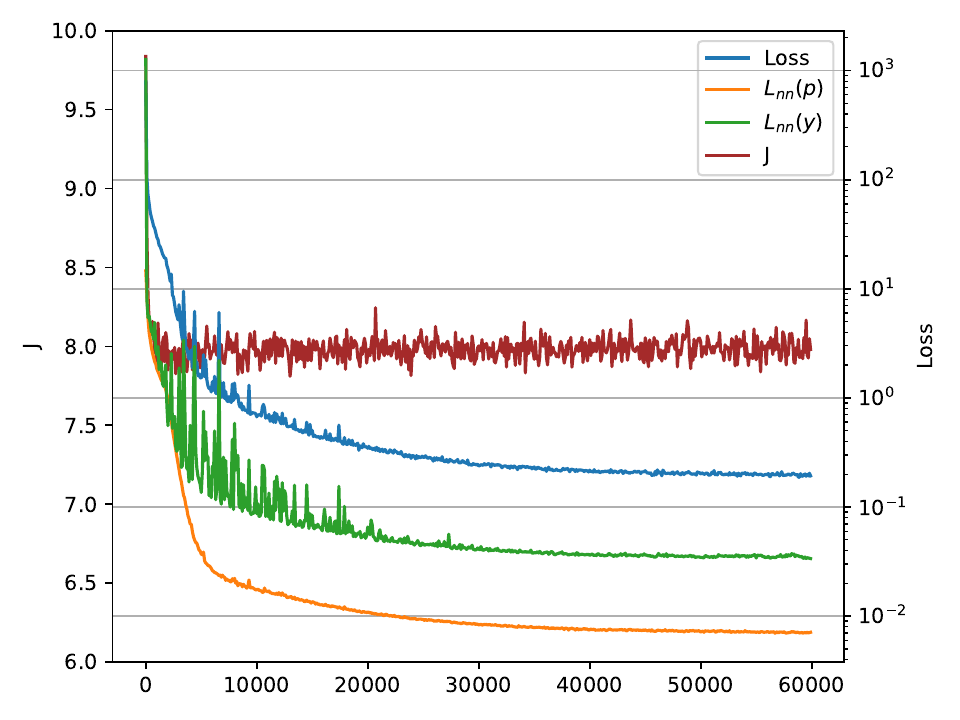}
        & \includegraphics[width=0.33\textwidth]{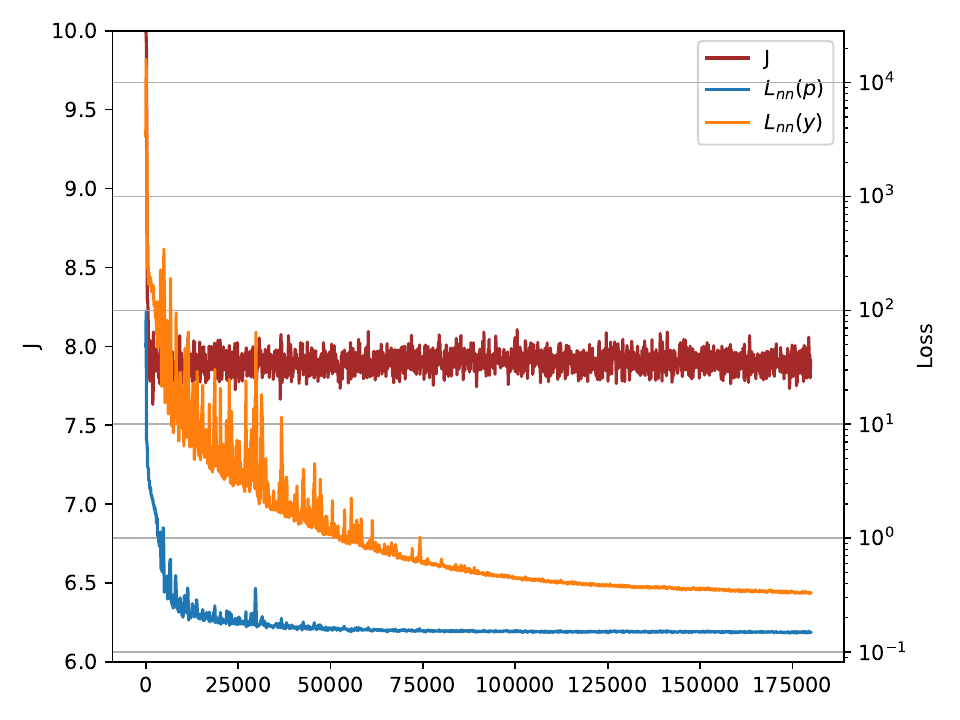}
        & \includegraphics[width=0.33\textwidth]{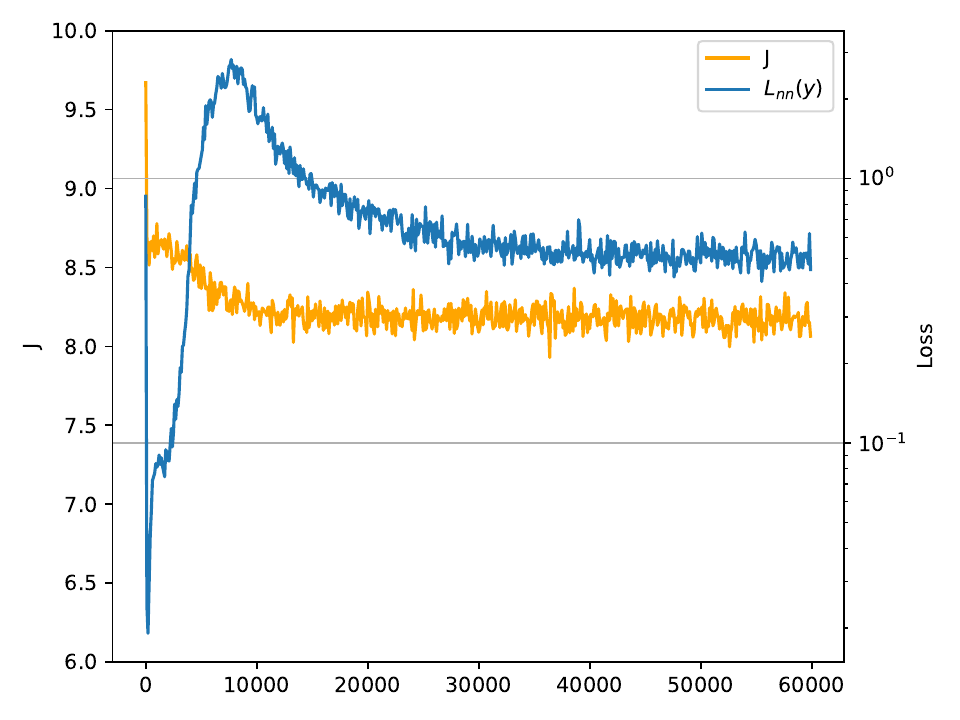}\\
        (a) OSNN & (b) AONN
        & (c) ALM
    \end{tabular}
    \caption{The training dynamics of the three methods, for Example \ref{exam:hidim-unconstraint}. Ticks on the left y-axis refer to cost objective $J$, and that on the right to losses $\widehat{\mathcal{L}}_{\rm nn}(y)$ and $\widehat{\mathcal{L}}_{\rm nn}(p)$.}
    \label{fig:hidim-unconst-dyn-exam1}
\end{figure}

In Fig. \ref{fig:hidim-unconst-dyn-exam1}, we show the training dynamics of the three methods. For OSNN, the optimizer converges rapidly and stably in terms of the loss value, so are the losses $\widehat{\mathcal{L}}_{\rm nn}(y)$ and $\widehat{\mathcal{L}}_{\rm nn}(p)$ for the state and adjoint. For AONN, it does minimize the loss steadily, but the sub-losses $\widehat{\mathcal{L}}_{\rm nn}(y)$ and $\widehat{\mathcal{L}}_{\rm nn}(p)$ exhibit alternating convergence behavior due to the alternating nature of the algorithm.
The first step leads to a drastic increase of the objective value, indicating an invalid descent direction. This issue is resolved after a few steps when the PDE solutions become more accurate. ALM enjoys fast initial decay in the loss but the convergence slows down greatly as the iteration proceeds, and also the loss $\widehat{\mathcal{L}}_{\rm nn}(y)$ decreases rapidly. 

Next, we consider a constrained problem.
\begin{example}\label{exam:4d-constraint}
In the polar coordinate $x=(r\cos\theta_1,r\sin(\theta_1)\cos \theta_2,r\sin \theta_1 \sin \theta_2 \cos \theta_3 ,r\sin \theta_1 \sin \theta_2 \sin \theta_3)$, $\Omega=\{(r,\theta_1,\theta_2,\theta_3):r\in[1,3],\theta_1,\theta_2,\theta_3\in[0,2\pi]\}$ is a hyper-annulus in $\mathbb{R}^4$ with outer radius $3$ and inner radius $1$. The control is subject to $u\in U: = \{u\in L^2(\Omega): u(r,\theta_1,\theta_2,\theta_3)\in [-0.5,0.7] \mbox{ a.e. in } \Omega\}$. $\lambda=0.01$, and data $y_d(r,\theta_1,\theta_2,\theta_3)=r^2+\lambda (5-\frac{9}{r^2})\cos\theta_1, f(r,\theta_1,\theta_2,\theta_3)=-8-{P}_{[-0.5,0.7]}(-(r-1)(r-3)\cos\theta_1)$ and $g(r)=r^2$. The exact state $\bar y$ and control $\bar u$ are given by
$\bar y(r,\theta_1,\theta_2,\theta_3)=r^2$ and $ \bar u(r,\theta_1,\theta_2,\theta_3)={P}_{[-0.5,0.7]}(-(r-1)(r-3)\cos\theta_1)$. 
\end{example}

\begin{figure}[h]
\centering
\setlength{\tabcolsep}{0pt}
\begin{tabular}{cccc}
  \includegraphics[width=0.24\textwidth]{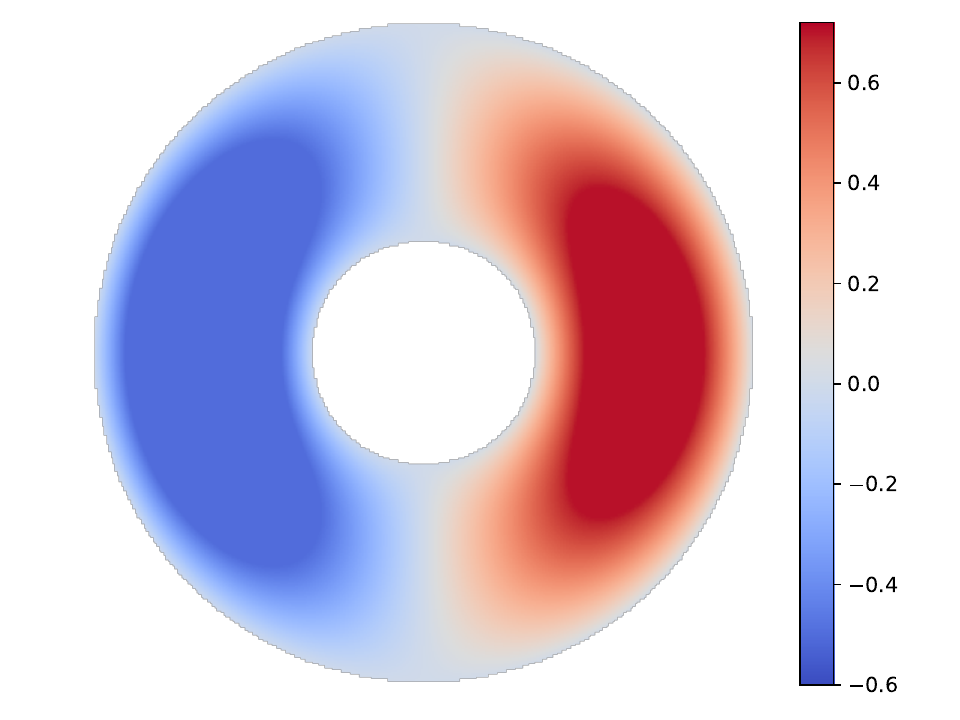}
& \includegraphics[width=0.24\textwidth]{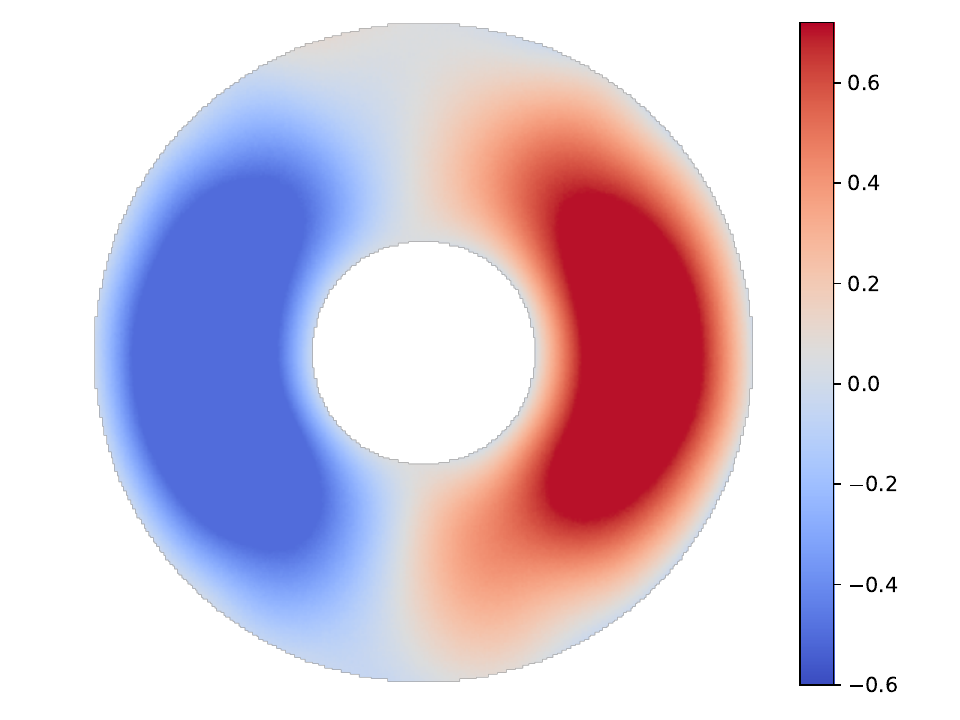}
& \includegraphics[width=0.24\textwidth]{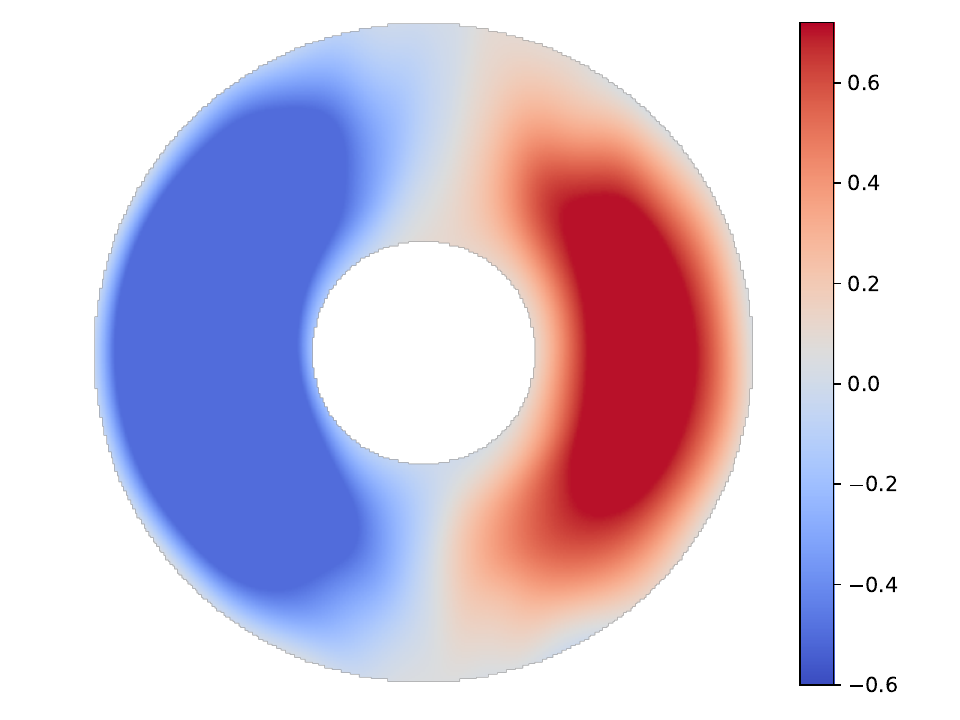}
& \includegraphics[width=0.24\textwidth]{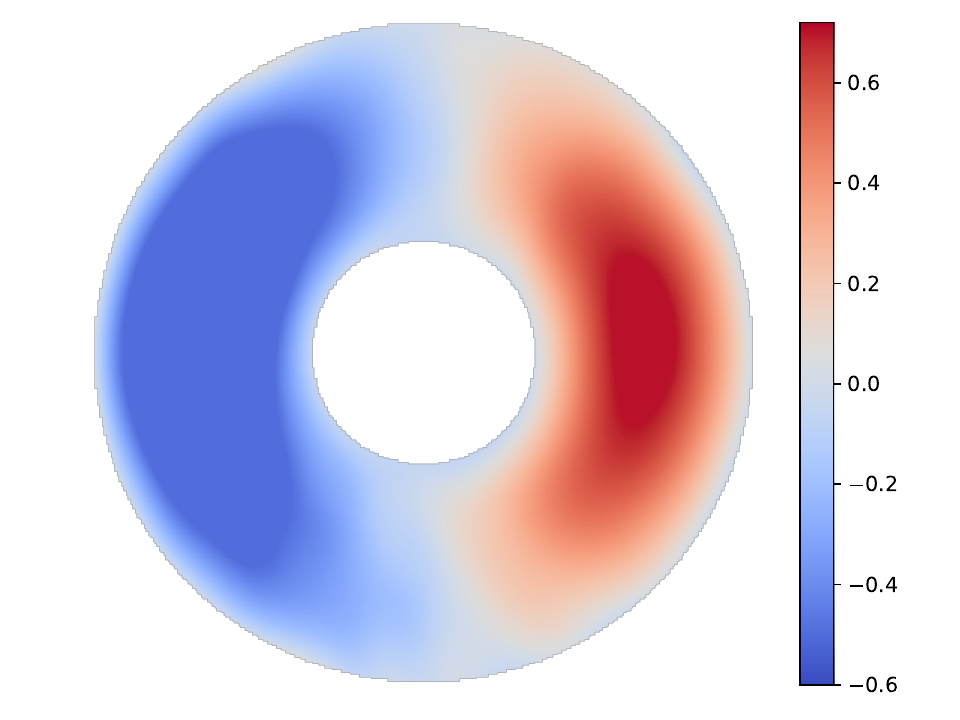}\\
& \includegraphics[width=0.24\textwidth]{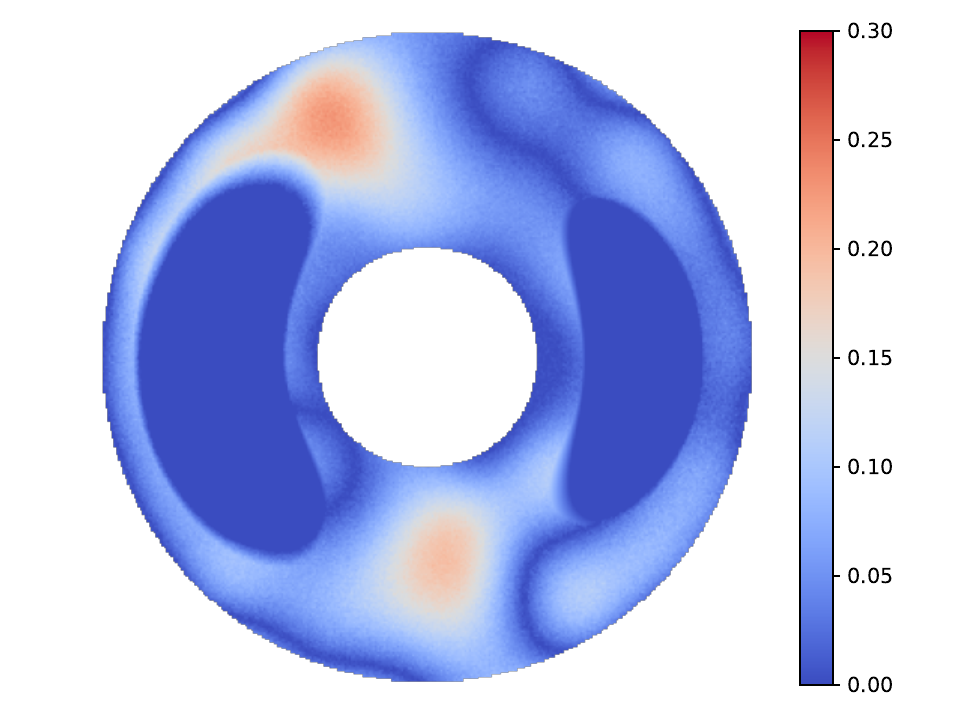}
& \includegraphics[width=0.24\textwidth]{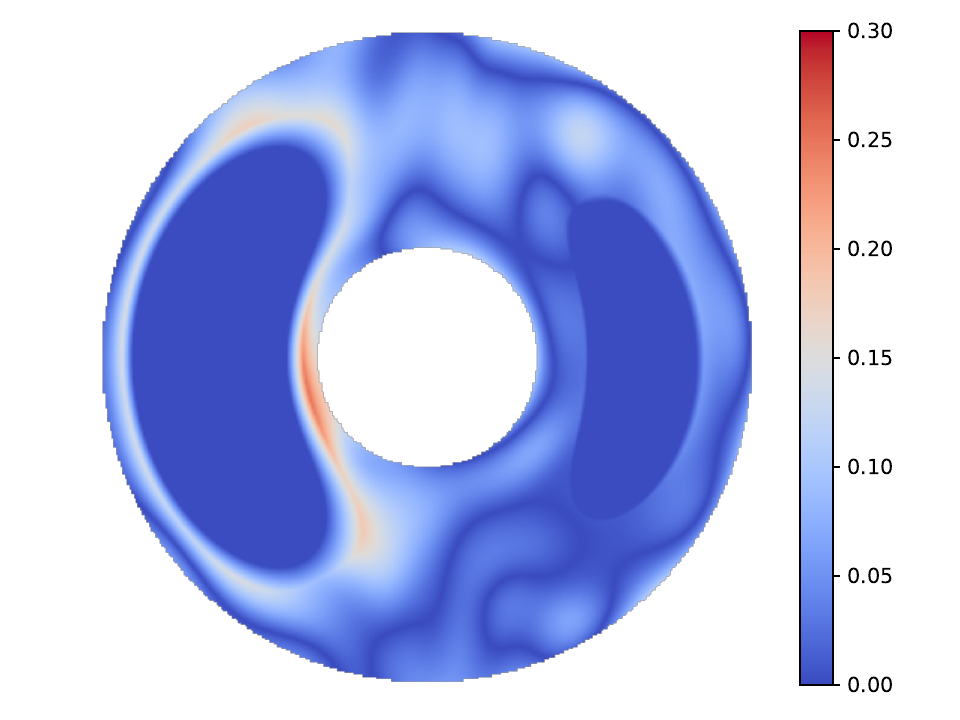}
& \includegraphics[width=0.24\textwidth]{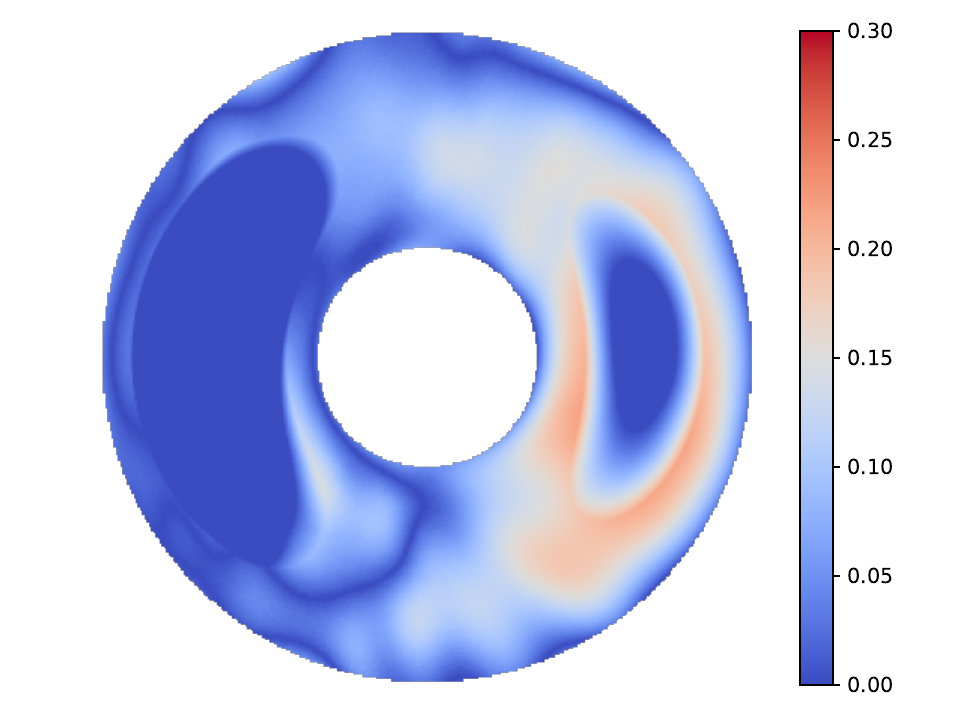}\\
  \includegraphics[width=0.24\textwidth]{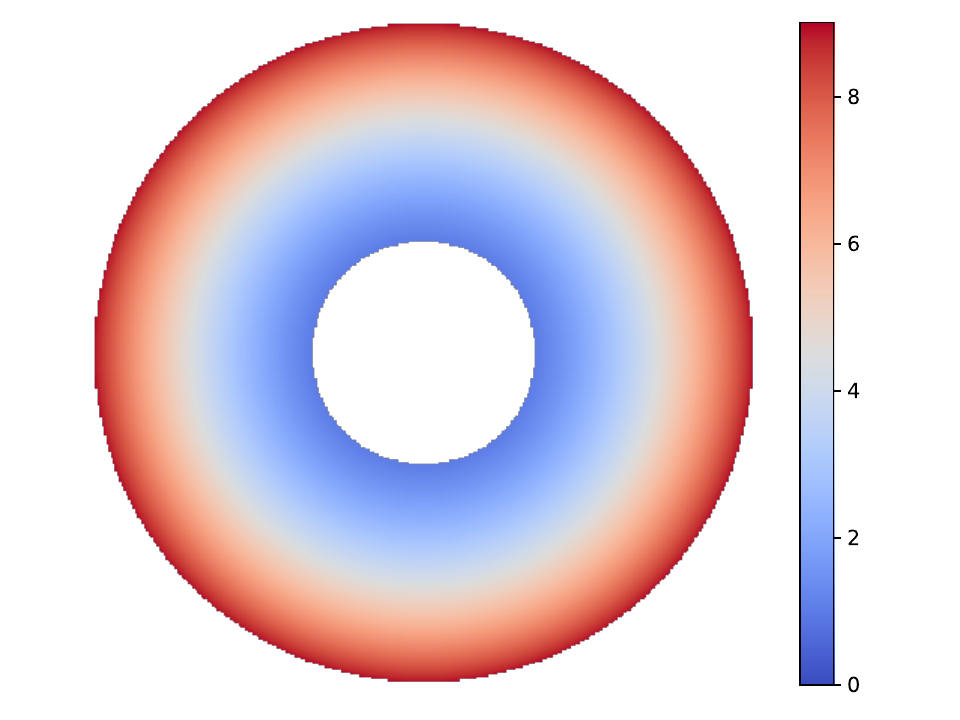}
& \includegraphics[width=0.24\textwidth]{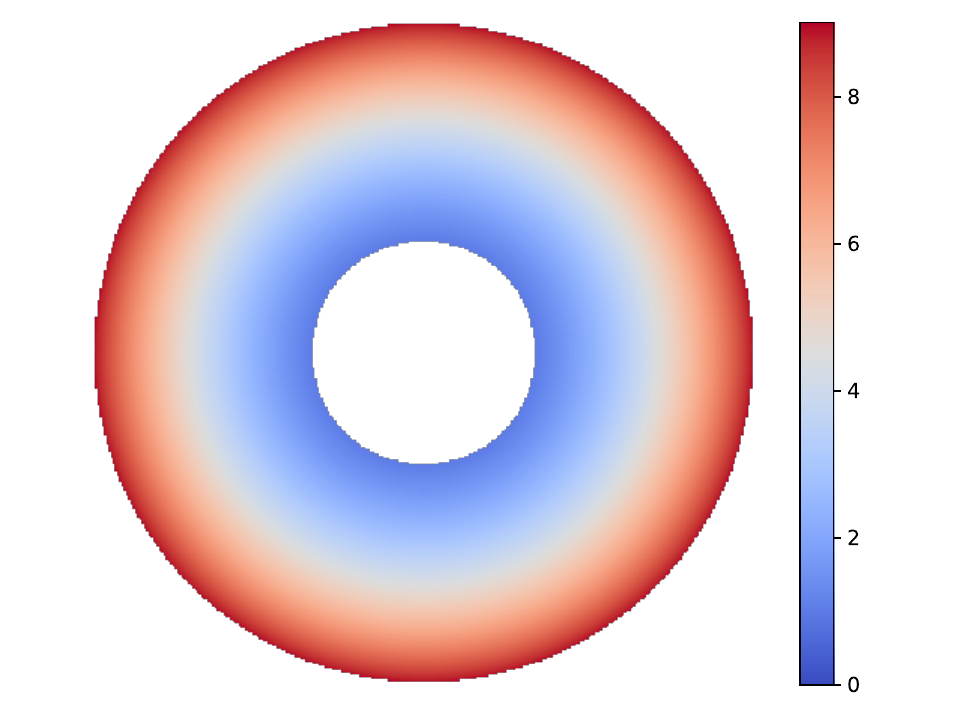}
& \includegraphics[width=0.24\textwidth]{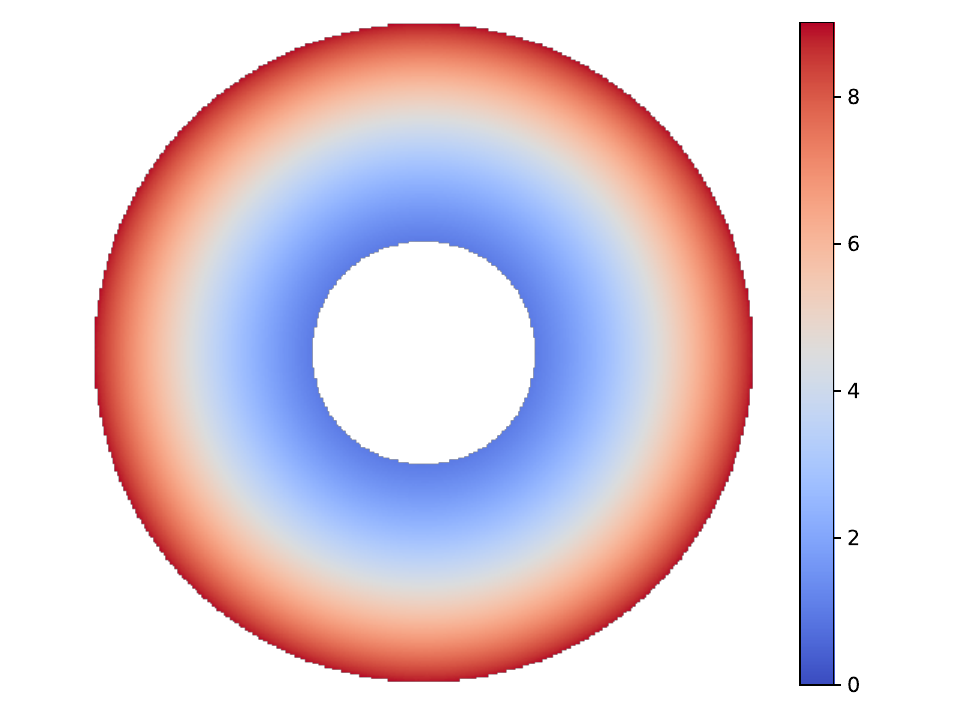}
& \includegraphics[width=0.24\textwidth]{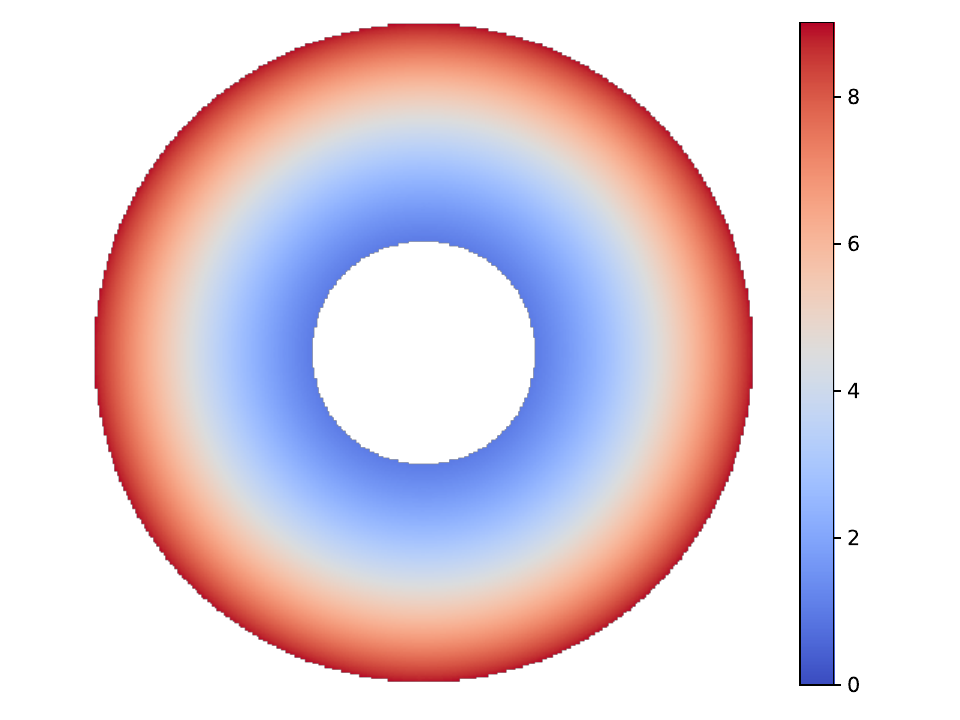}\\
(a) exact & (b) OSNN & (c) AONN  & (d) ALM
\end{tabular}
\caption{The approximate control $u^*$ (top), its pointwise error $|u^*-\bar u|$ (mid) and approximate state $y^*$ (bottom) obtained by three NN-based methods for Example \ref{exam:4d-constraint}, cross section at $x_3=0,x_4=0$.}
    \label{fig:4d-const}
\end{figure}

The experimental setting is identical with Example \ref{exam:hidim-unconstraint}. In ALM, we enforce the box constraint by including a penalty term with a weight $\mu' = 0.01$, cf. \eqref{eqn:pm-constraint}, and add a multiplier term and update it by Uzawa algorithm. Since the forward solver converges slower, we increase the number of iterations for all methods. The learning rates are halved after every 6000 iterations. The results are reported in Table \ref{tab:exam}(b) and Fig. \ref{fig:4d-const}. Since the optimal control $\bar u$ is non-smooth near the boundary of the constraint active area, ALM and AONN suffer from big errors therein, and AONN has to use a smaller step size. Table \ref{tab:exam}(b) shows that AONN has a good accuracy of the state approximation with 100 descent steps. For ALM, the projection penalty term requires one extra tunable weight. In Table \ref{tab:exam}(b), ALM approximates the state accurately, but the control error $e(u^*)$ is still large. The computing time (in second) for OSNN, AONN and ALM is 1.56e3, 4.42e3, and 1.70e3. While the errors are larger than others, ALM yields the best cost, implying that the PDE constraint is not well-satisfied. The approximate state $y^*$ by AONN is less accurate with a slightly better control $u^*$. Both OSNN and ALM  optimize two variables and thus requires more iterations for the PDE, but the control of OSNN converges faster.

Next we consider a semilinear example.
\begin{example}\label{exam:semilinear}
Consider the following semilinear PDE:
$-\Delta y+y +y^3=f+u$ in $ \Omega$ with $y = g$ on $\partial\Omega$,
with the domain $\Omega=(0,1)^4$.
Let $\lambda=0.01$. The data $y_d$, $g$ and $ f$ are taken such that the exact state $\bar y$ and control $\bar u$ are given by $\bar y(x)=\sum_{i=1}^3 e^{x_i(1-x_i)}\sin(\pi x_{i+1})+e^{x_4(1-x_4)}\sin(\pi x_1)$, and $\bar u(x)=-\lambda^{-1}\prod_{i=1}^4 x_i(1+\cos(\pi x_i))$.
\end{example}

\begin{figure}[h]
    \centering
    \setlength{\tabcolsep}{0pt}
    \begin{tabular}{cccc}
             \includegraphics[width=0.24\textwidth]{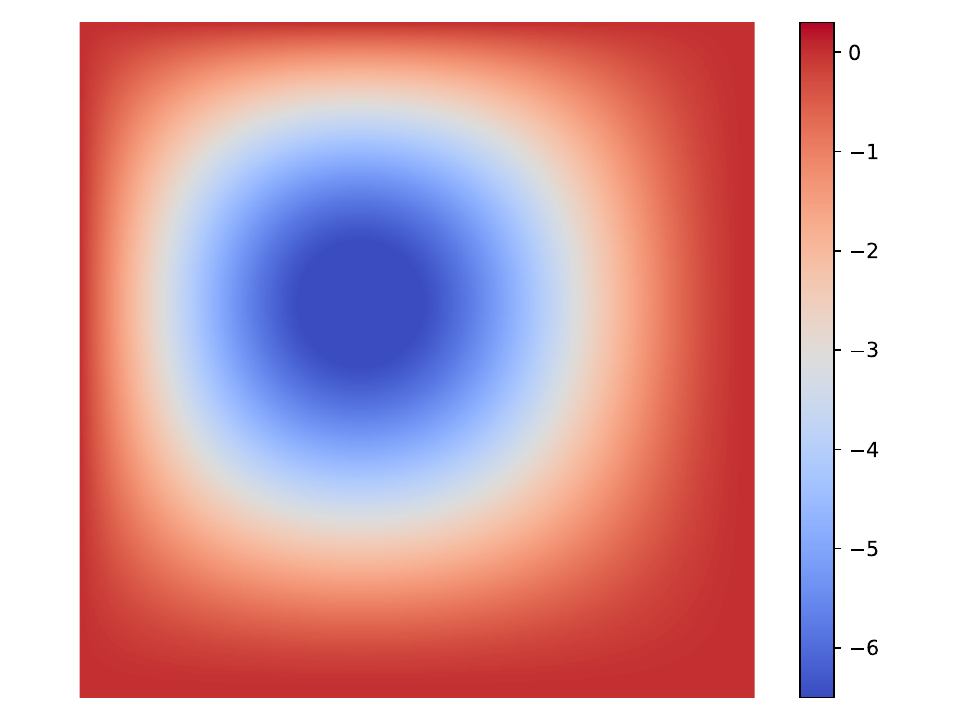}
            &\includegraphics[width=0.24\textwidth]{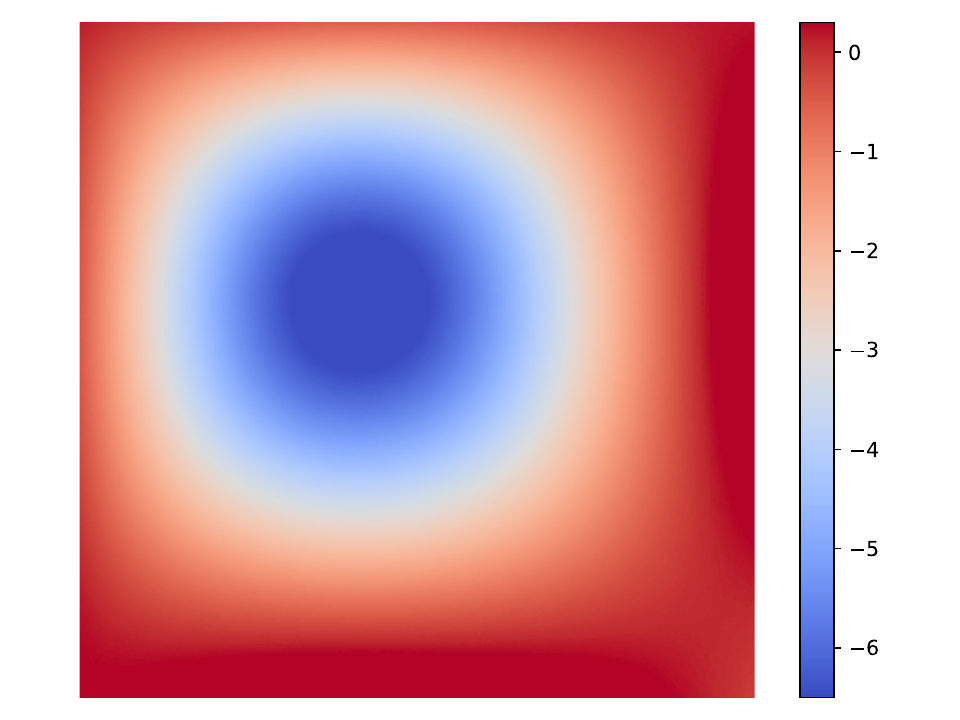}
            &\includegraphics[width=0.24\textwidth]{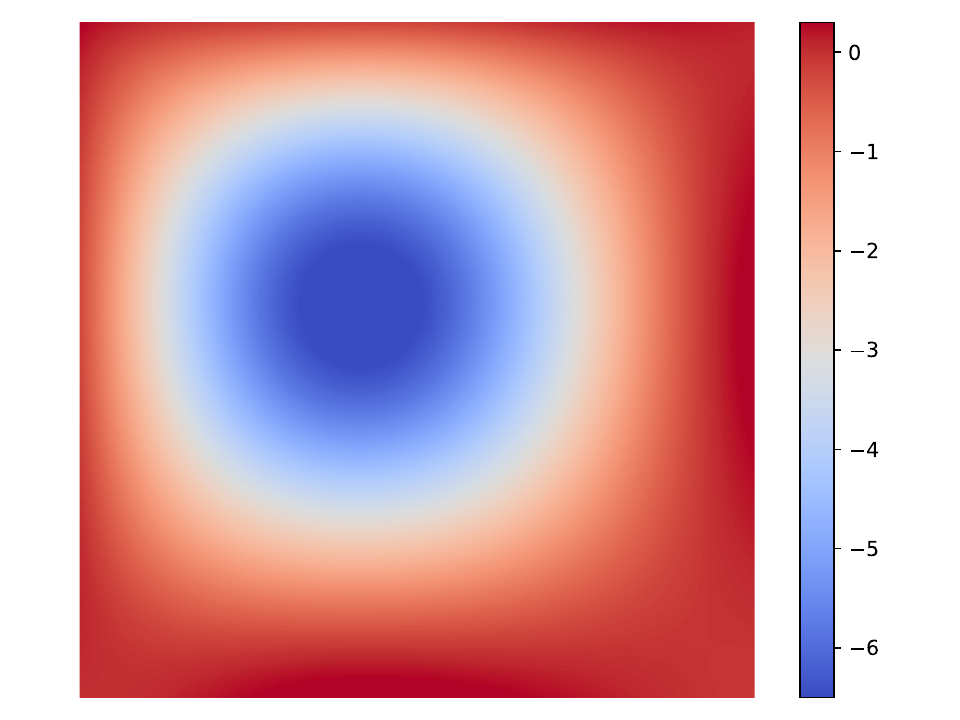}
            &\includegraphics[width=0.24\textwidth]{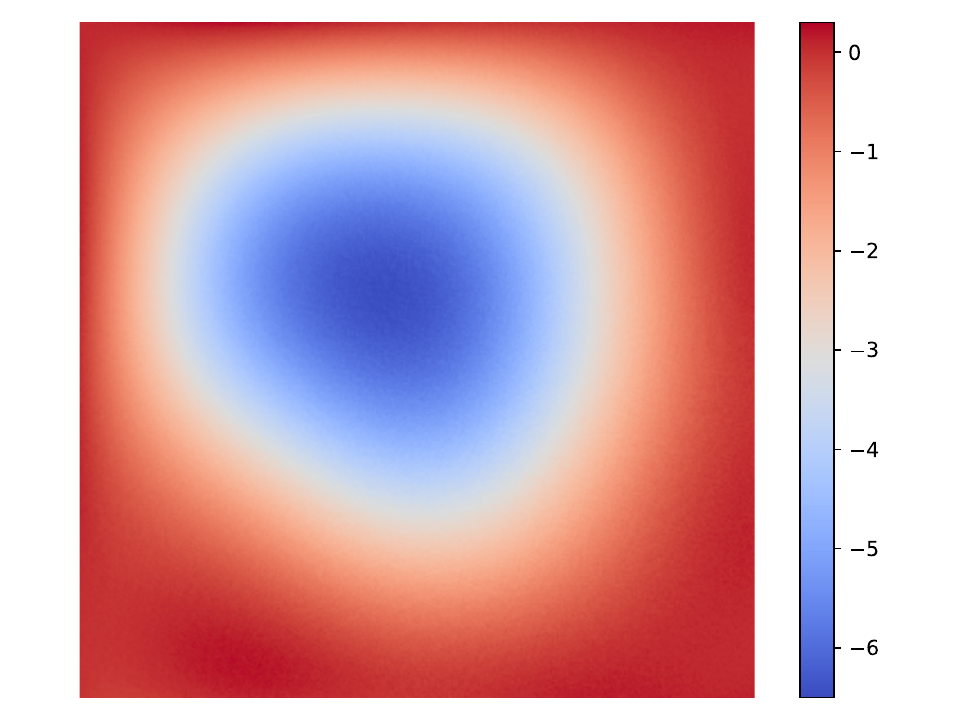}\\
            &\includegraphics[width=0.24\textwidth]{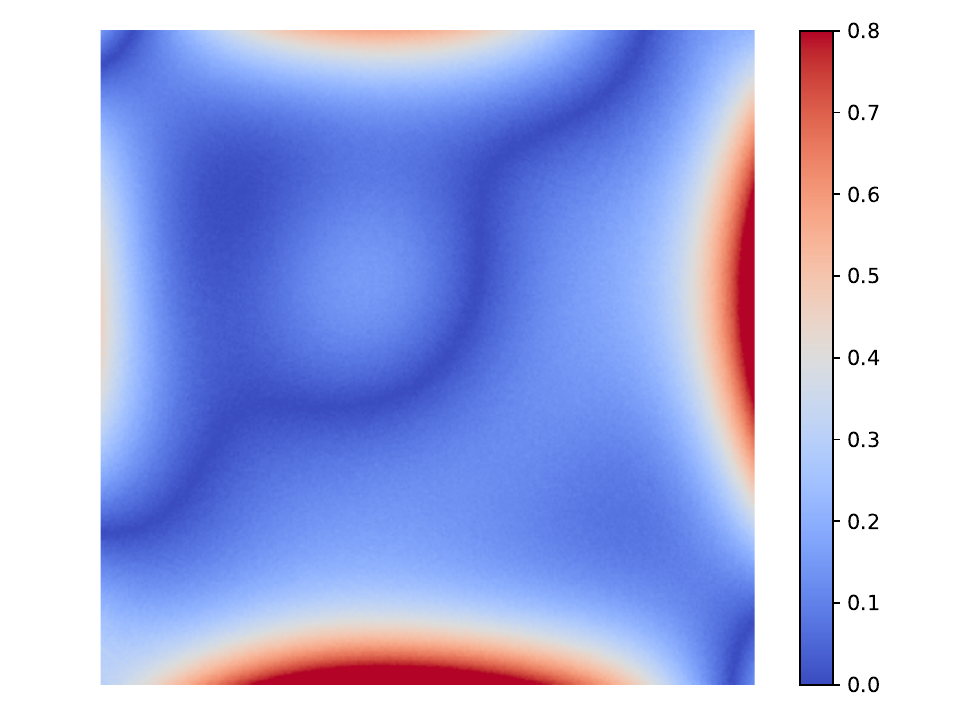}
            &\includegraphics[width=0.24\textwidth]{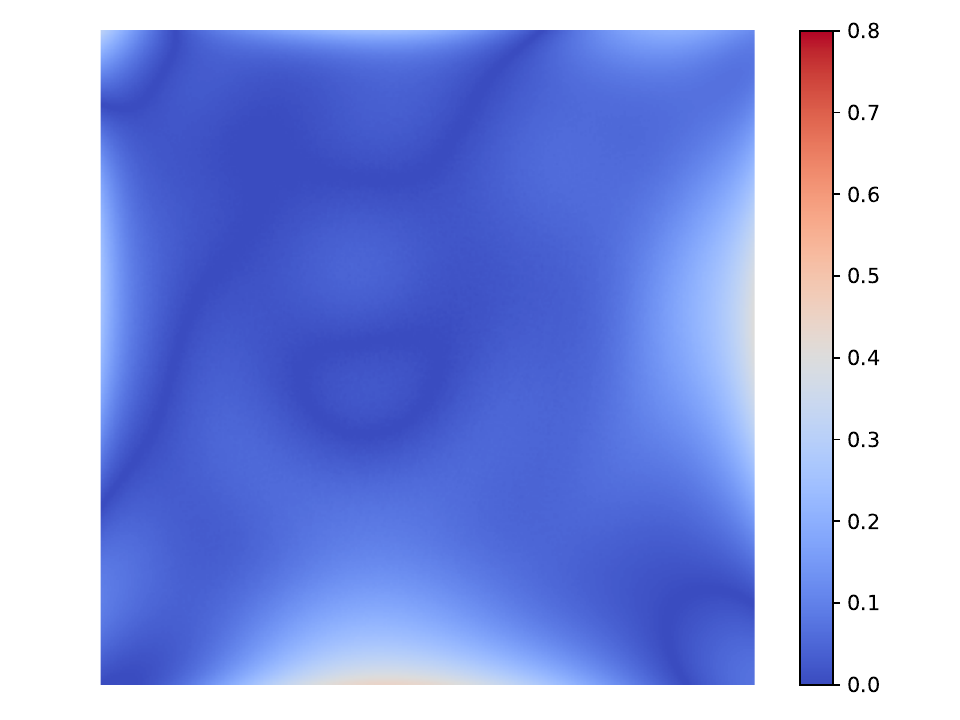}
            &\includegraphics[width=0.24\textwidth]{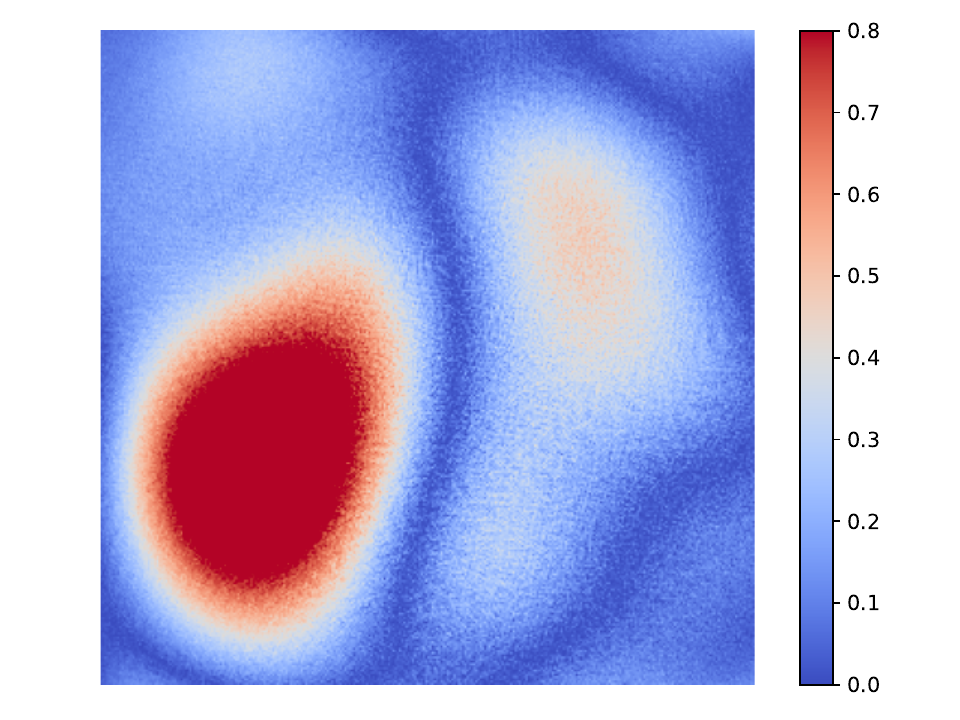}\\
             \includegraphics[width=0.24\textwidth]{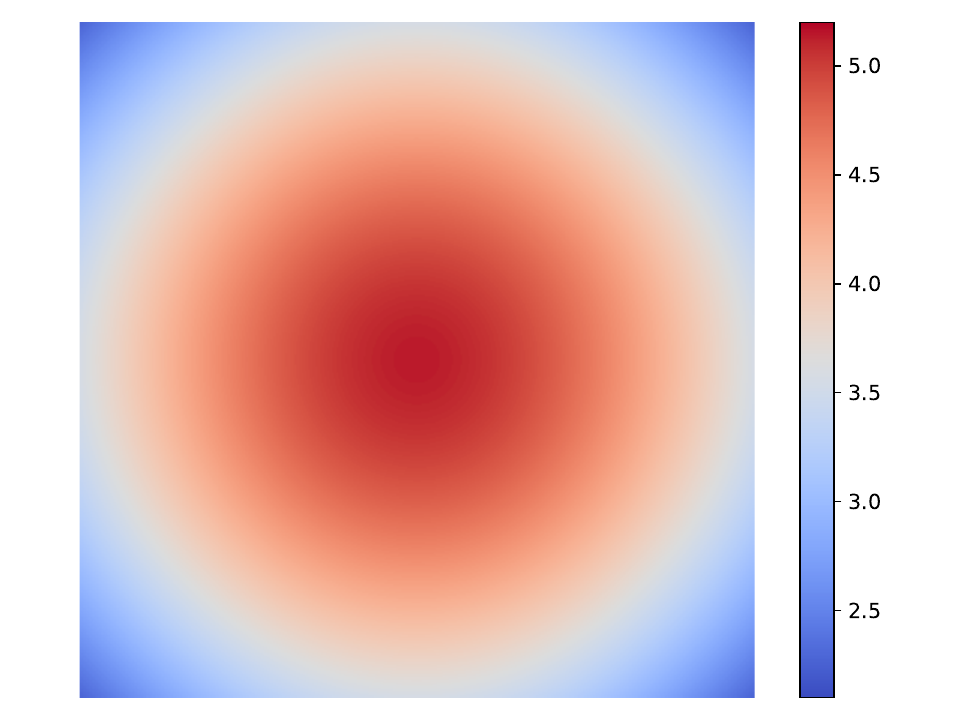}
            &\includegraphics[width=0.24\textwidth]{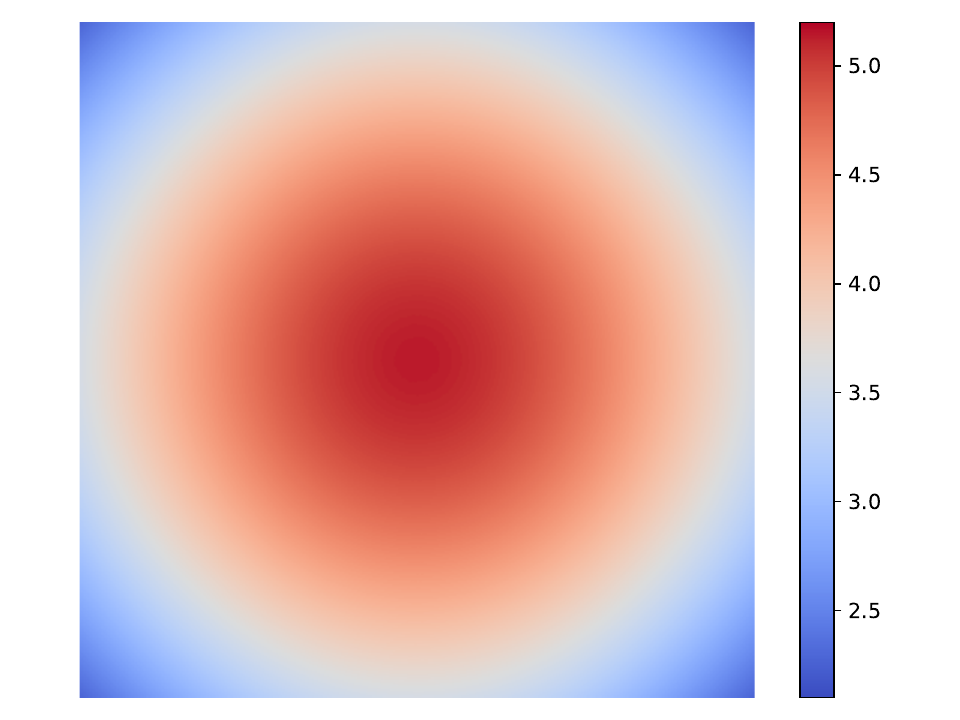}
            &\includegraphics[width=0.24\textwidth]{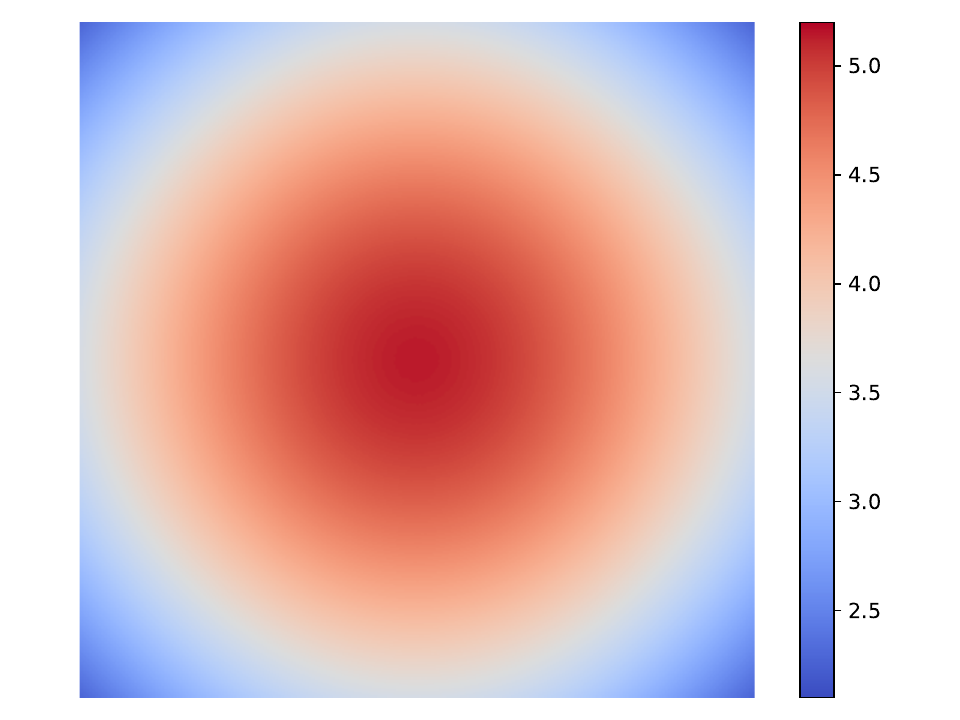}
            &\includegraphics[width=0.24\textwidth]{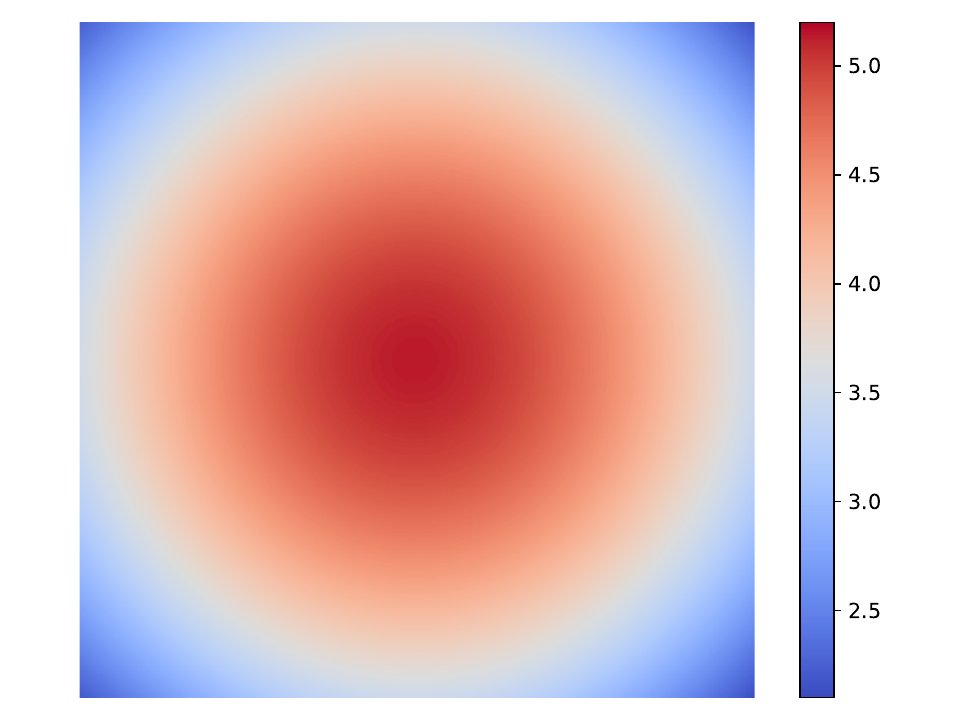}\\
            (a) exact & (b) OSNN & (c) AONN& (d) ALM
    \end{tabular}
    \caption{The approximate optimal control $u^*$ (top), its pointwise error $|u^*-\bar u|$ (middle), and the state $y^*$ (bottom) obtained by three different methods for Example \ref{exam:semilinear}, cross section at $(x_3,x_4)=(0.5,0.5)$.}
    \label{fig:semilinear}
\end{figure}

The results are reported in Table \ref{tab:exam}(c) and Fig. \ref{fig:semilinear}. Both OSNN and AONN approximate the optimal state $\bar y$ and control $\bar u$ accurately, but AONN has much larger control error $e(u^*)$. The efficiency of AONN depends heavily on the step size $s$.
The step size $s=20$ was determined in a trial-and-error manner in order to achieve good efficiency. ALM yields an accurate state approximation, which is however still inferior to OSNN and AONN. The computing time (in seconds) is only slightly higher than the unconstrained case: 4.39e3, 5.07e3 and 4.54e3 for OSNN, AONN and ALM, respectively.

Last, we investigate higher-dimensional problem to show the efficiency of OSNN.
\begin{example}\label{exam:highdimcpinn}
 $\Omega=(0,1)^6$, and  $\lambda=0.01$. The test cases are as follows: {\rm(i)} linear unconstrained (LU) $f=5\pi^2\prod_{i=1}^6\sin(\pi x_i)$ and $y_d=(-6\alpha\pi^4-1)\prod_{i=1}^6\sin(\pi x_i)$; {\rm(ii)} linear constrained (LC) {$f=6\pi^2\prod_{i=1}^6\sin(\pi x_i)-P_{[0,3]}(\pi^2 \prod_{i=1}^6 \sin(\pi x_i))$, $y_d = (-6\alpha\pi^4-1)\prod_{i=1}^6\sin(\pi x_i)$}; and {\rm(iii)} semilinar (semi) $f=6\pi^2\prod_{i=1}^6\sin(\pi x_i)+\prod_{i=1}^6\sin^3(\pi x_i)$, $y_d=(1+6\pi^2\lambda + \lambda)\prod_{i=1}^6\sin(\pi x_i)+\lambda \prod_{i=1}^6\sin^3(\pi x_i)$, and $q(x,y)=y^3$. All PDEs have a zero Dirichlet boundary condition. The exact states $\bar y$ are given by $\prod_{i=1}^6 \sin(\pi x_i)$, and exact controls $\bar u$ are given by $\pi^2\prod_{i=1}^6 \sin(\pi x_i)$ for {\rm(i)}, ${P}_{[0,3]}(\pi^2\prod_{i=1}^6 \sin(\pi x_i))$ for {\rm(ii)}, $\prod_{i=1}^6 \sin(\pi x_i)$ for {\rm(iii)}.
\end{example}

The NN architecture is identical with that for Example \ref{exam:hidim-unconstraint}. The boundary weight $\alpha_b$ is 100. We take 60000 points in $\Omega$ and 5000 points on $\partial\Omega$. The loss is optimized with Adam, with maximum 6e4 iterations. The learning rate is initialized to 1e-3, and halved after every 6e3 iterations, with a minimum 1e-3/1024. The results are shown in Table \ref{tab:exam}(d) and the training dynamics in Fig. \ref{fig:highdimcpinn}. The cost objective and PINN losses converge stably in all the three cases. OSNN yields accurate control and state approximations for linear problems, and the control error is slightly larger for the semilinear problem. These results show the potential of OSNN for solving optimal control problems.

\newcommand\newwd{0.16}
\begin{figure}[h]
\centering
\setlength{\tabcolsep}{0pt}
\begin{tabular}{cccccc}
  \includegraphics[width=\newwd\textwidth]{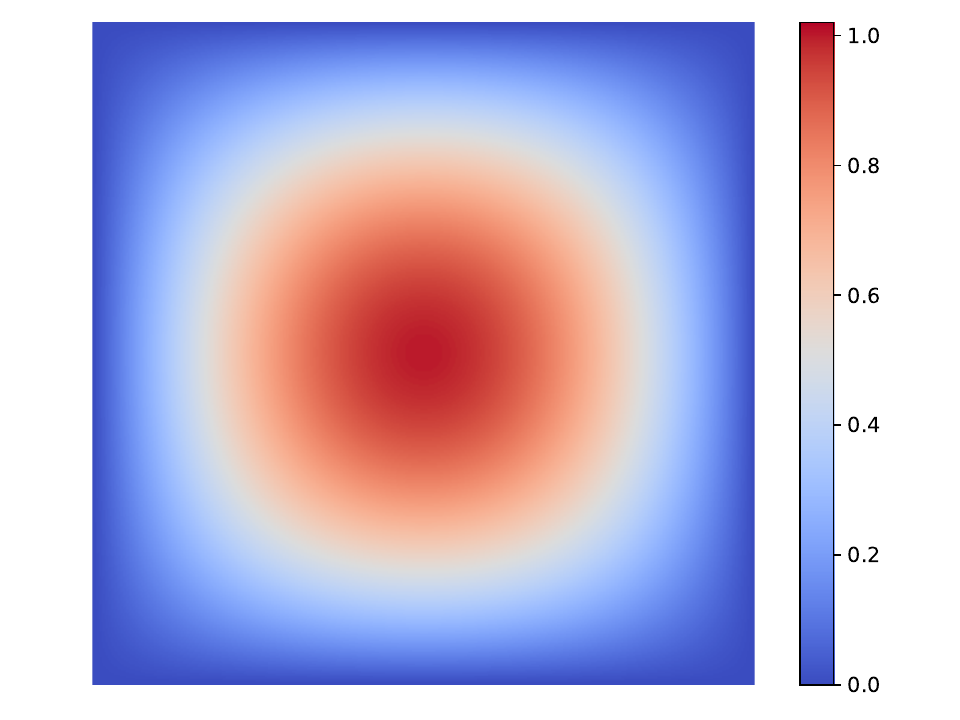}
& \includegraphics[width=\newwd\textwidth]{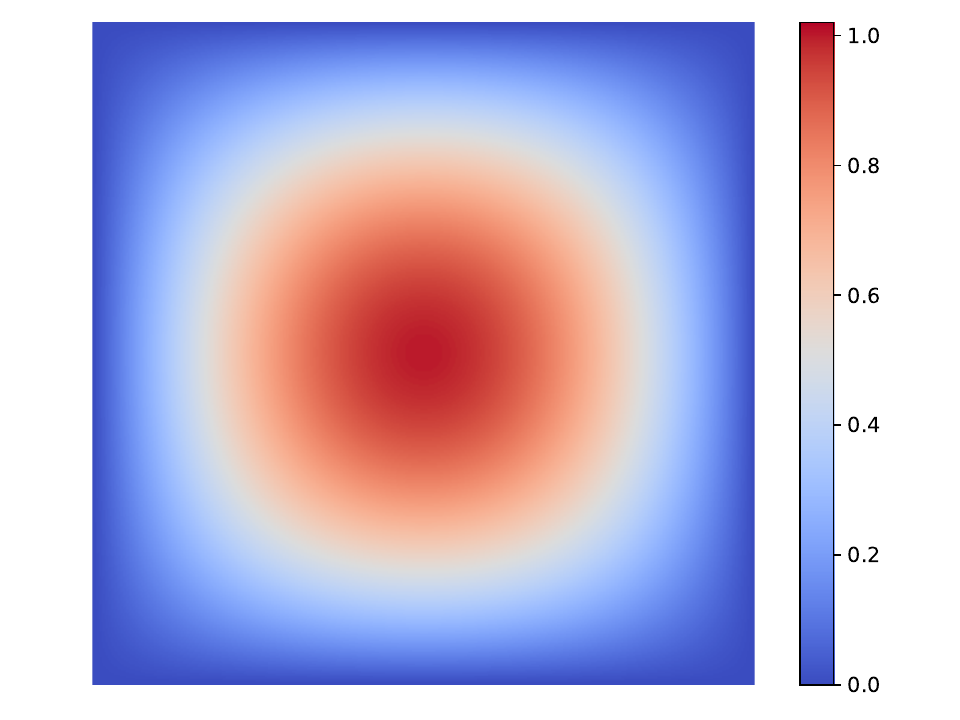}
& \includegraphics[width=\newwd\textwidth]{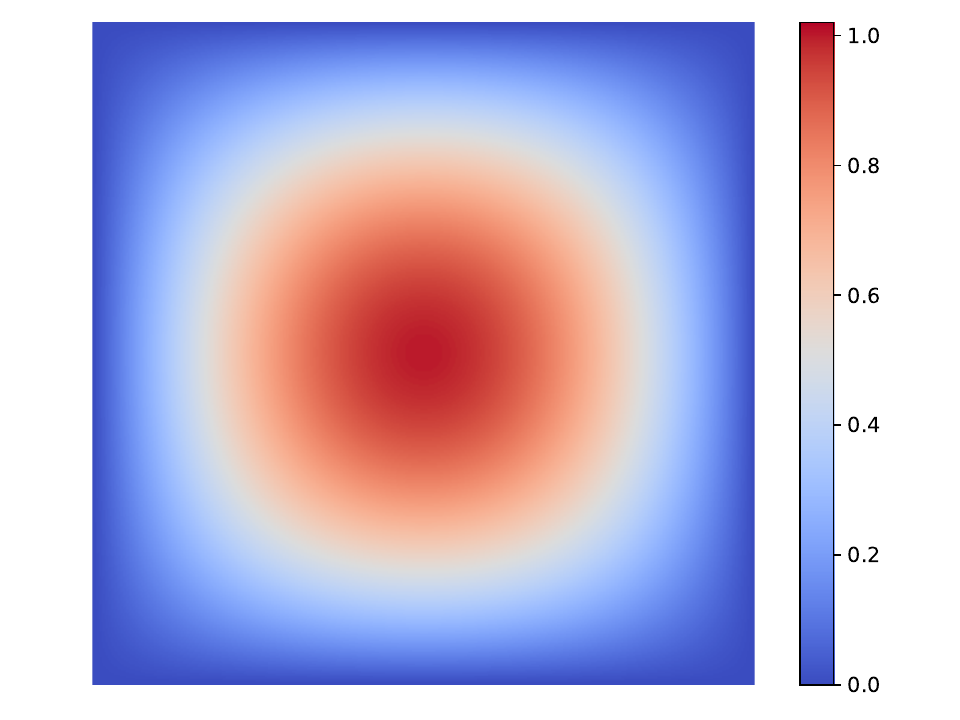}
& \includegraphics[width=\newwd\textwidth]{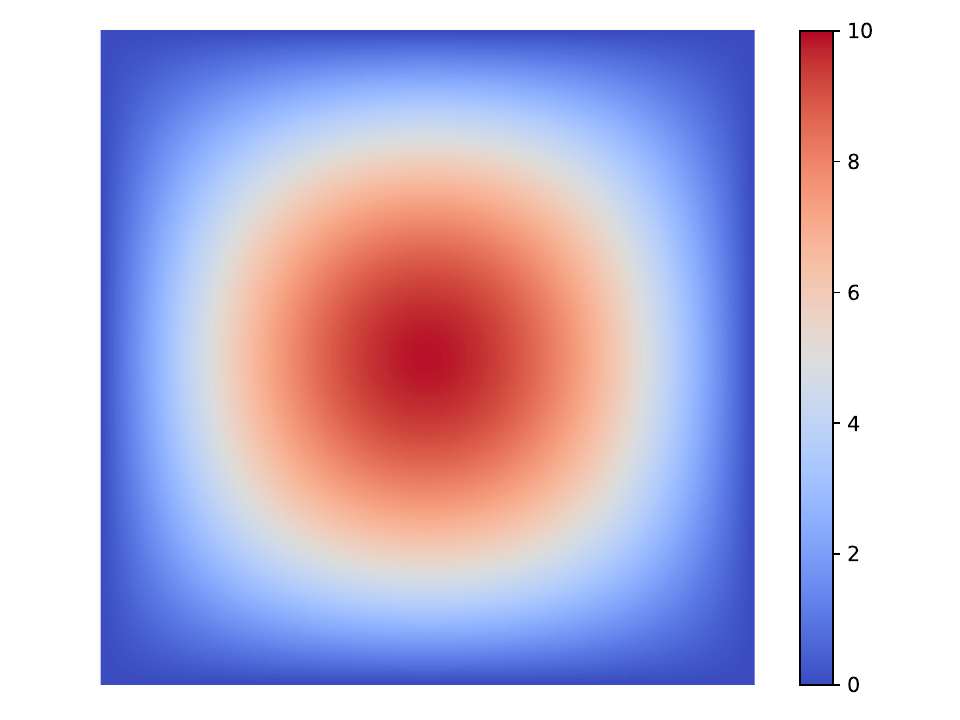}
& \includegraphics[width=\newwd\textwidth]{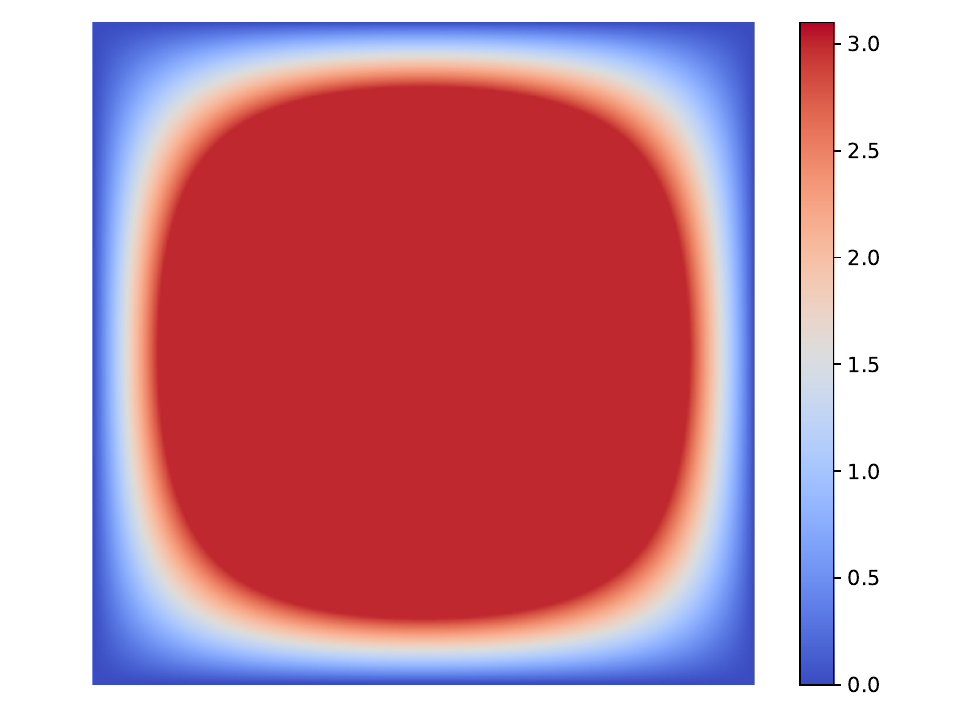}
& \includegraphics[width=\newwd\textwidth]{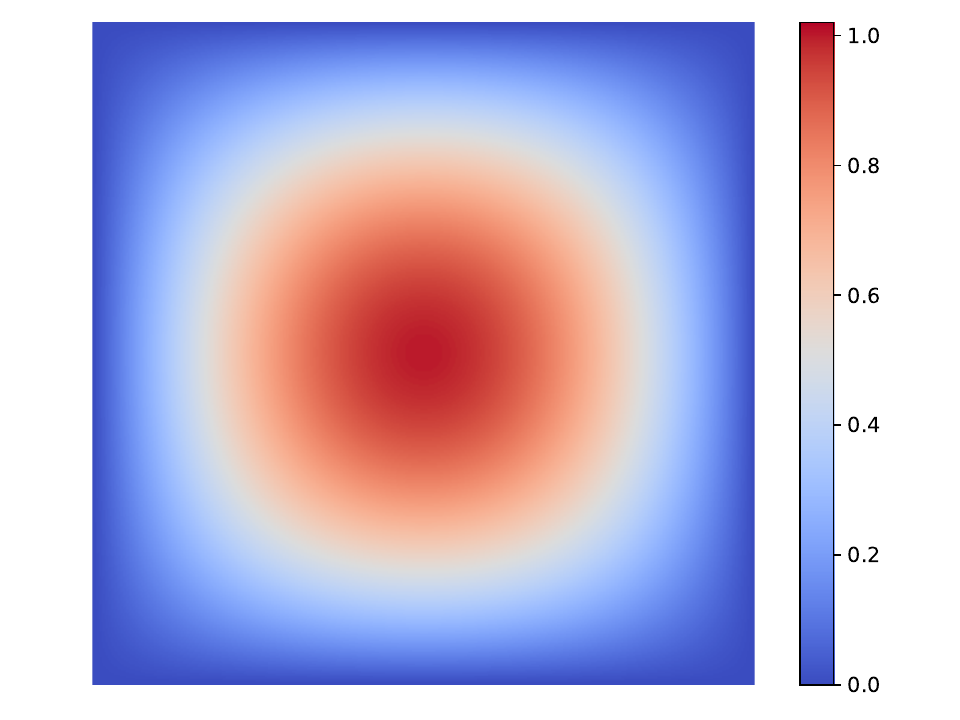}\\
  \includegraphics[width=\newwd\textwidth]{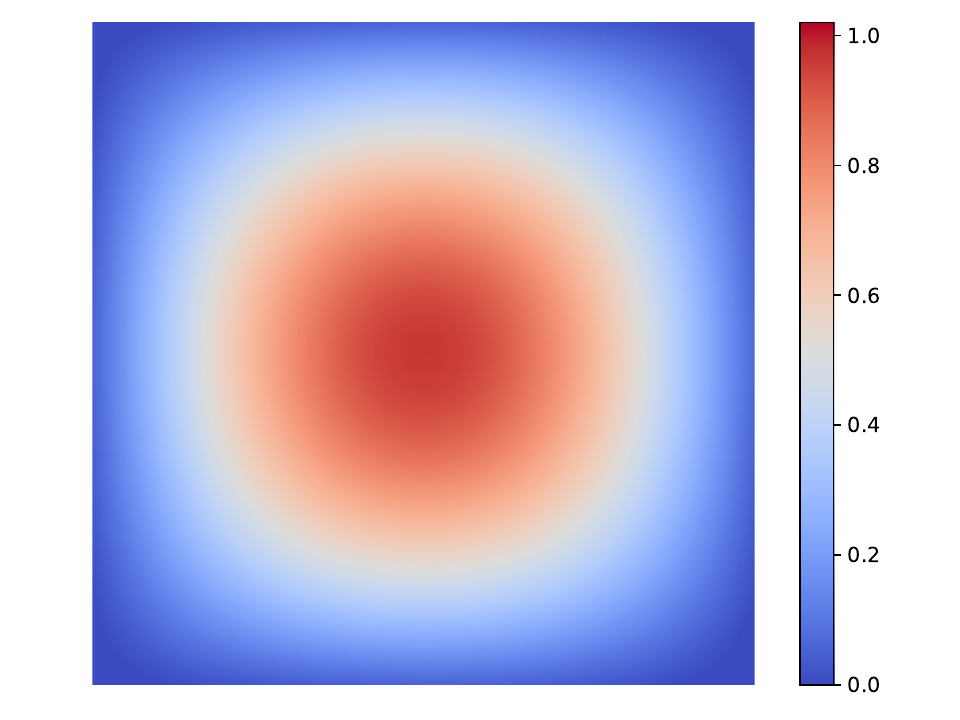}
& \includegraphics[width=\newwd\textwidth]{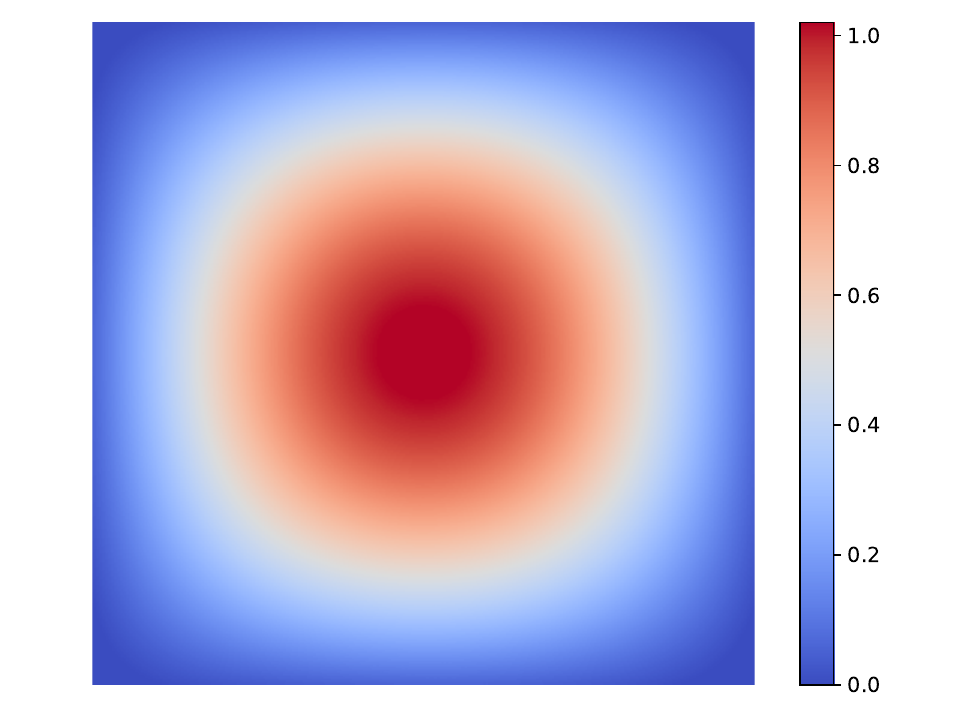}
& \includegraphics[width=\newwd\textwidth]{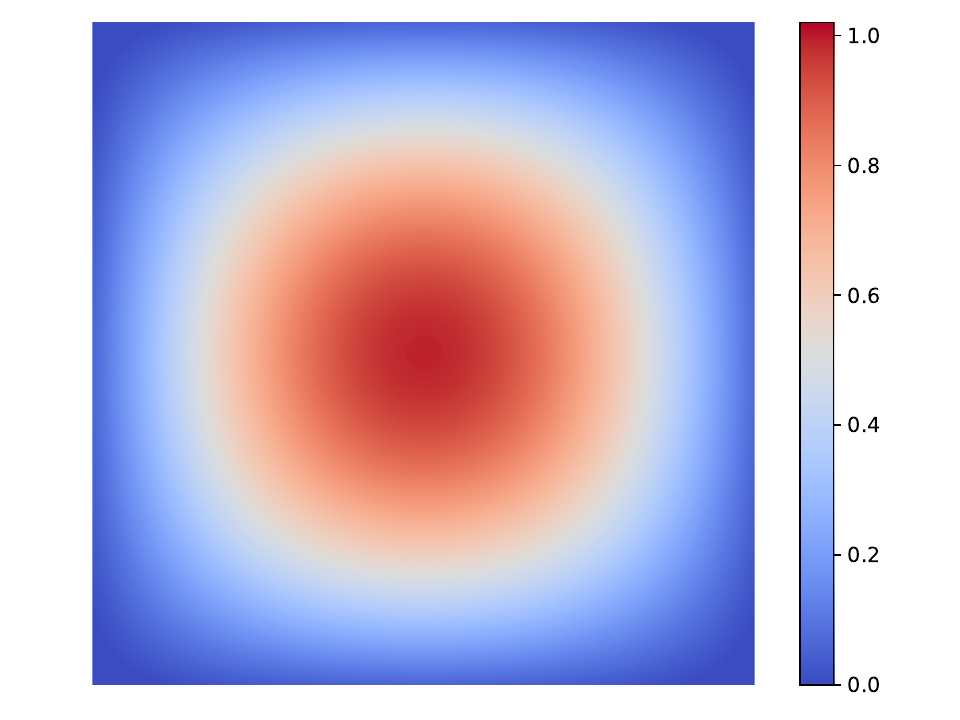}
& \includegraphics[width=\newwd\textwidth]{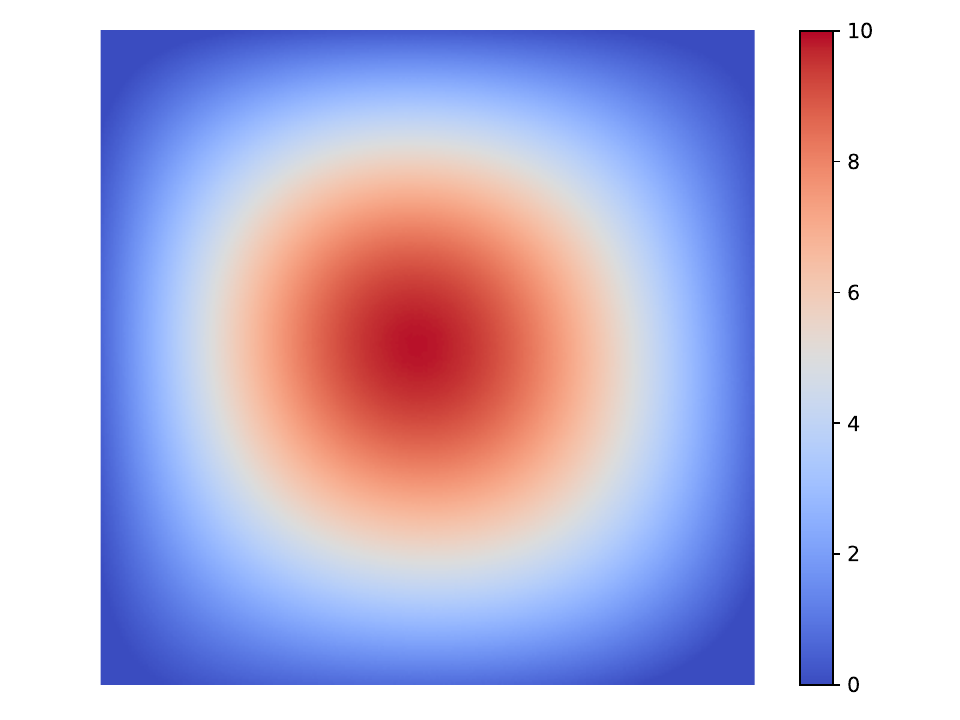}
& \includegraphics[width=\newwd\textwidth]{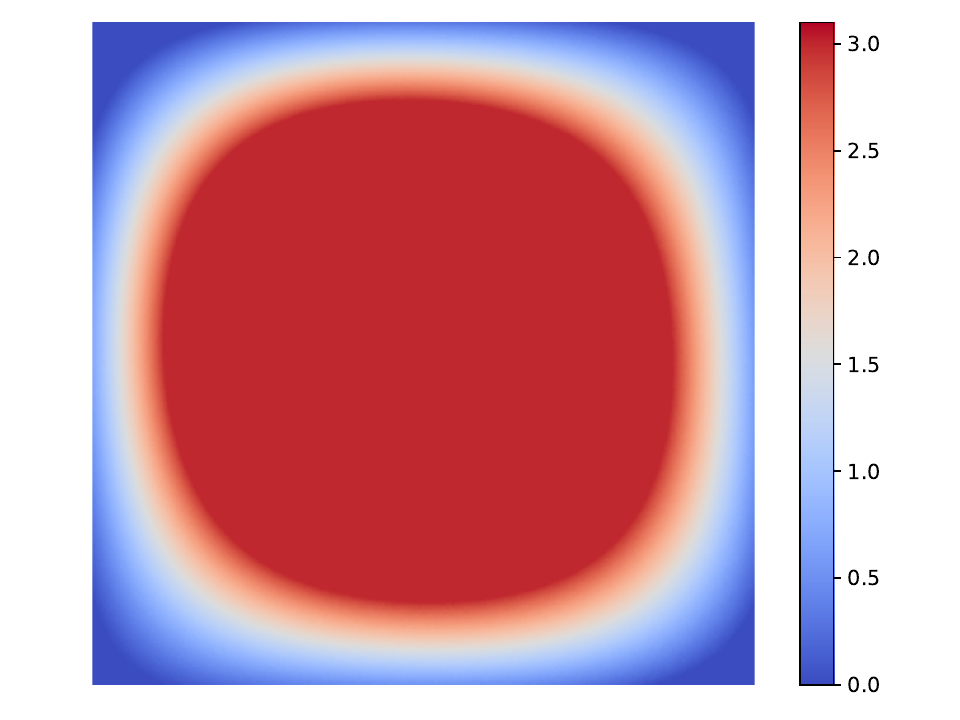}
& \includegraphics[width=\newwd\textwidth]{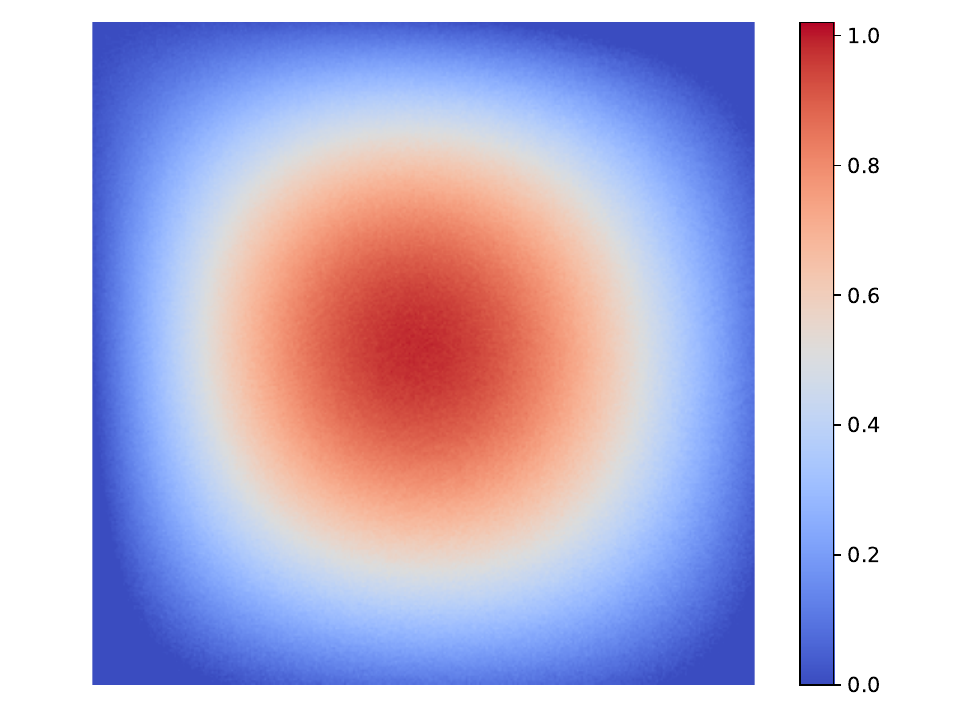}\\
  \includegraphics[width=\newwd\textwidth]{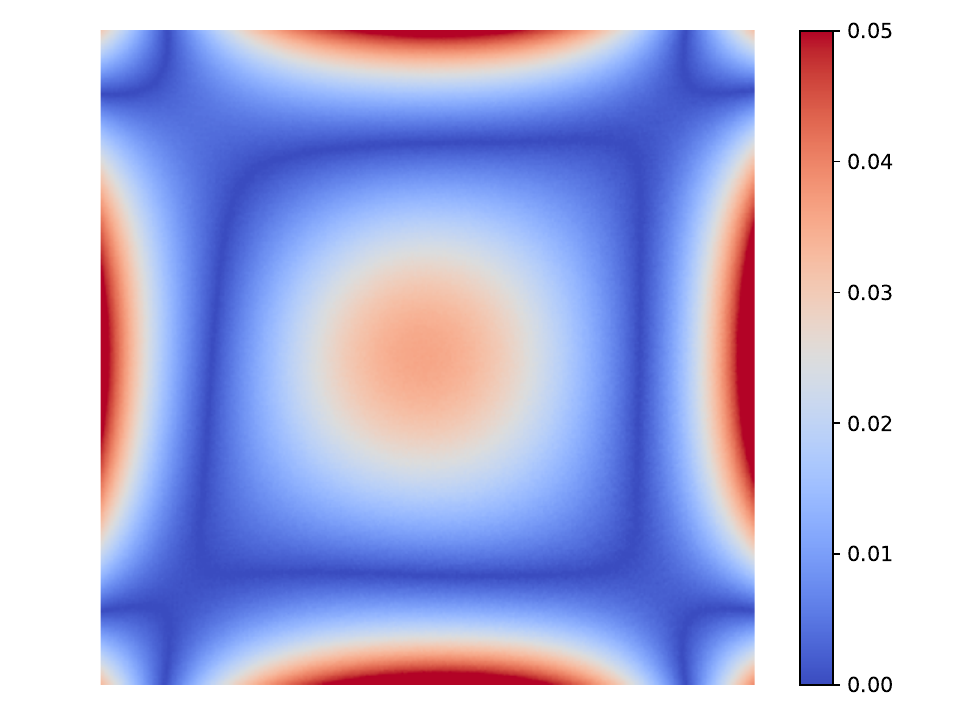}
& \includegraphics[width=\newwd\textwidth]{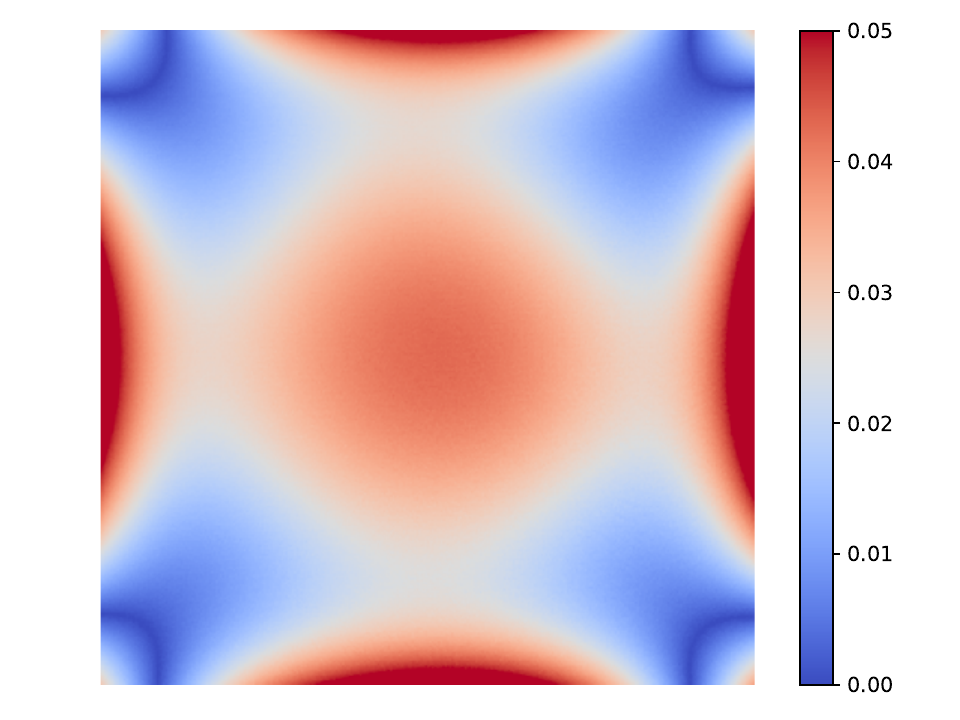}
& \includegraphics[width=\newwd\textwidth]{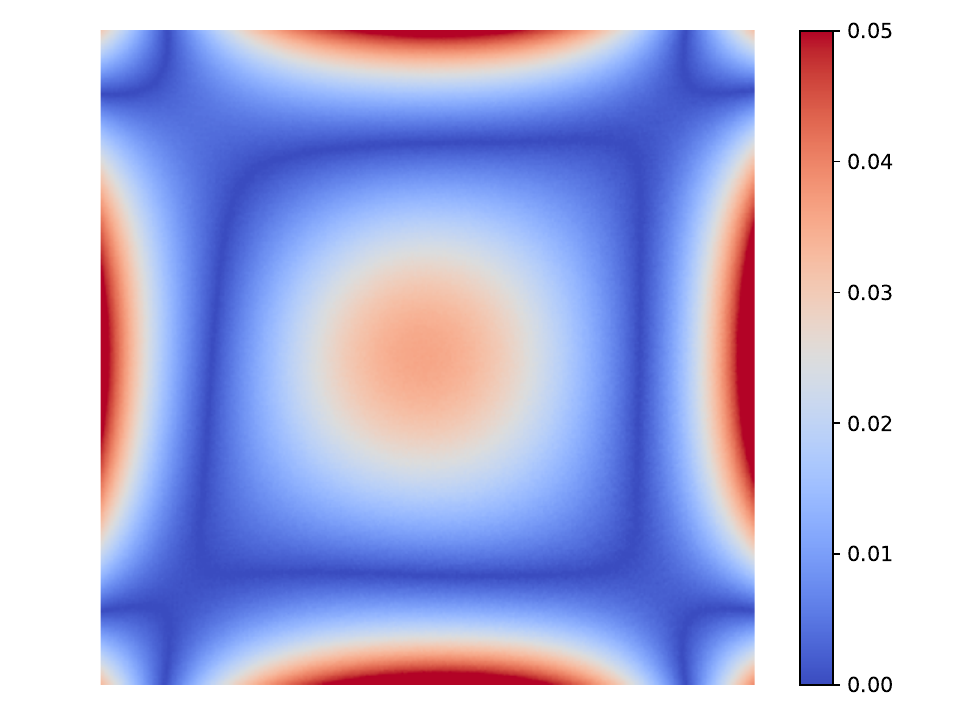}
& \includegraphics[width=\newwd\textwidth]{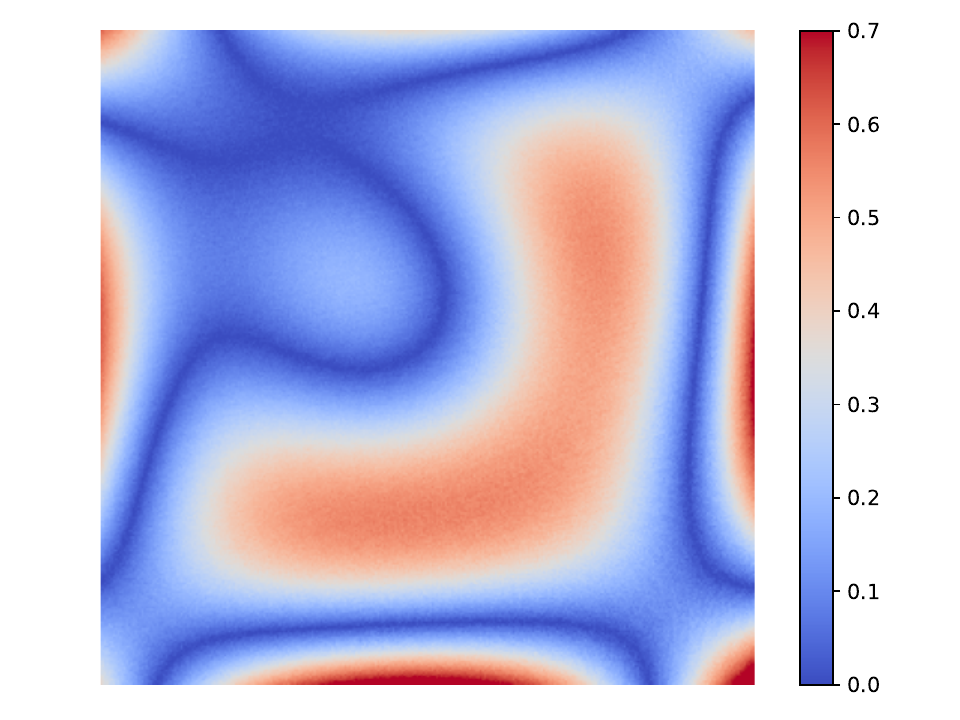}
& \includegraphics[width=\newwd\textwidth]{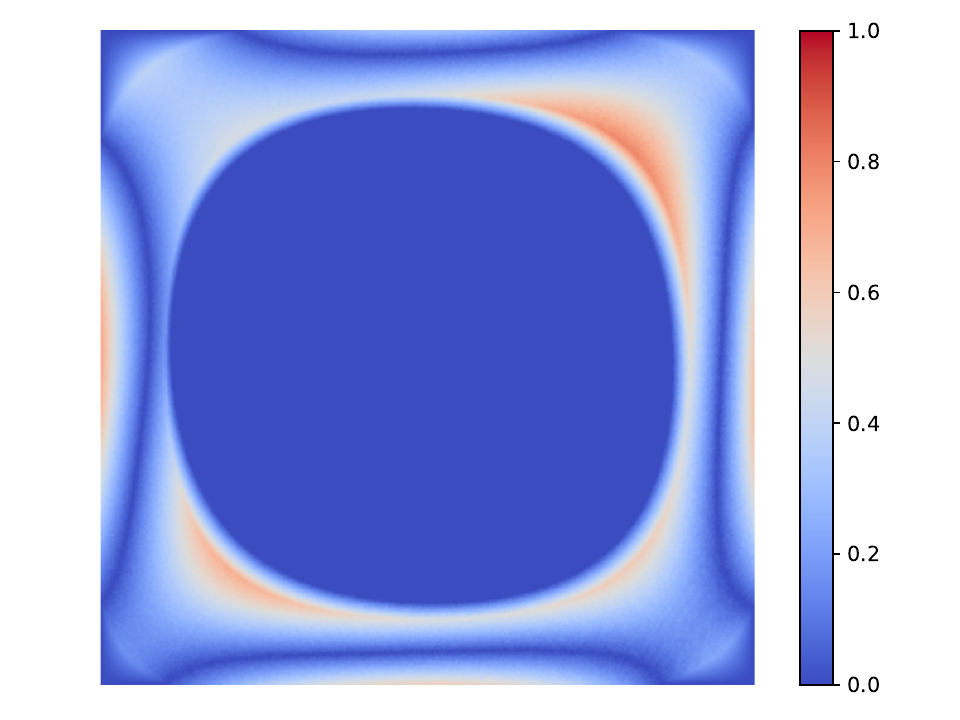}
& \includegraphics[width=\newwd\textwidth]{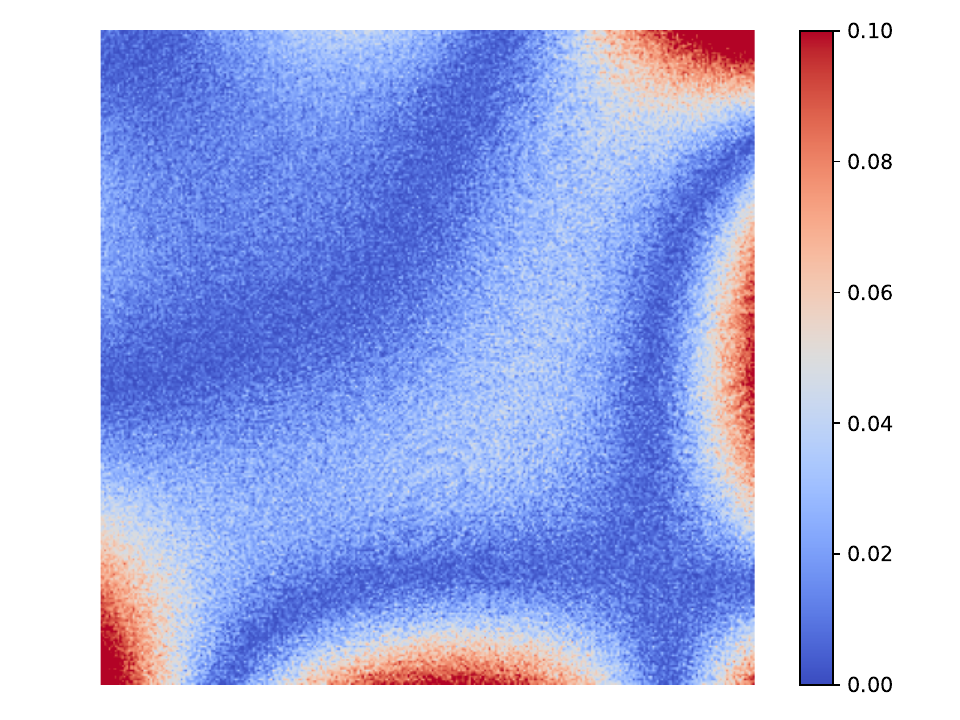}\\
(a) LU ($y$) & (b) LC ($y$) & (c) Semi ($y$) & (d) LU ($u$) & (e) LC ($u$) & (f) Semi ($u$)
\end{tabular}
\caption{The exact solutions (top), approximate
 solutions (middle) and pointwise errors (bottom) by OSNN for Example \ref{exam:highdimcpinn}, cross section at $x_3=x_4=x_5=x_6=0.5$.}
    \label{fig:highdimcpinn}
\end{figure}

\begin{figure}[hbt]
    \centering
    \setlength{\tabcolsep}{0pt}
    \begin{tabular}{ccc}
          \includegraphics[width=0.33\textwidth]{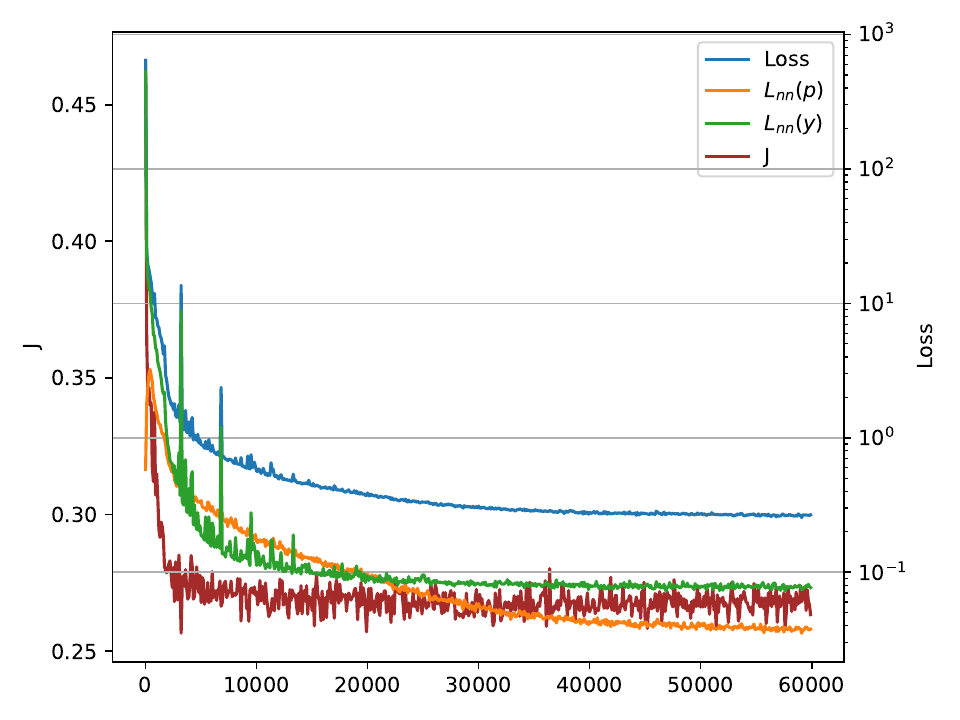}
        & \includegraphics[width=0.33\textwidth]{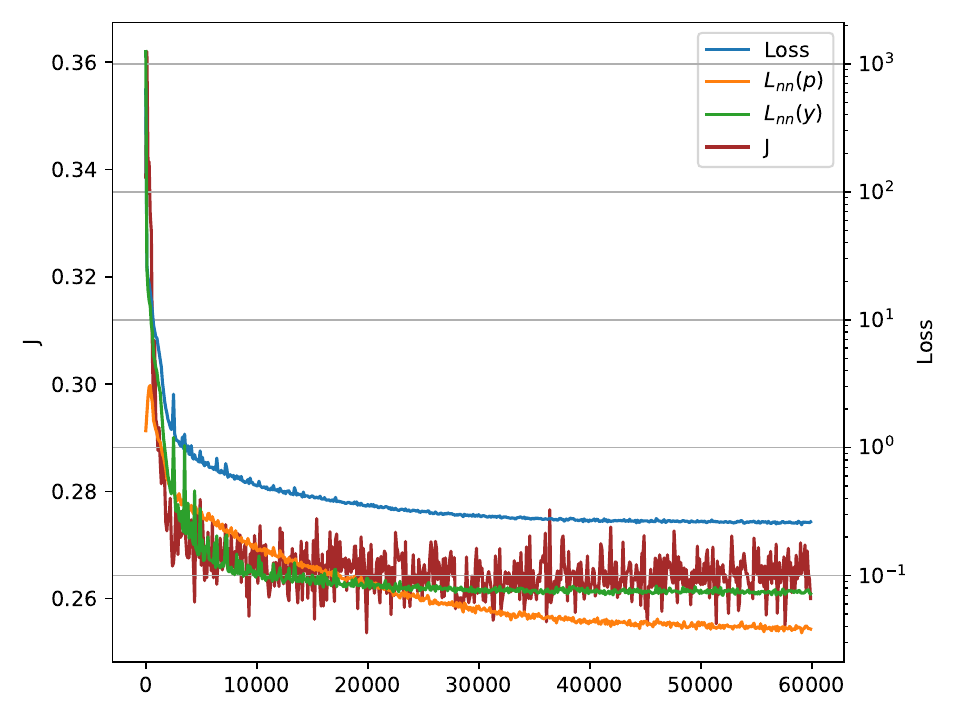}
        & \includegraphics[width=0.33\textwidth]{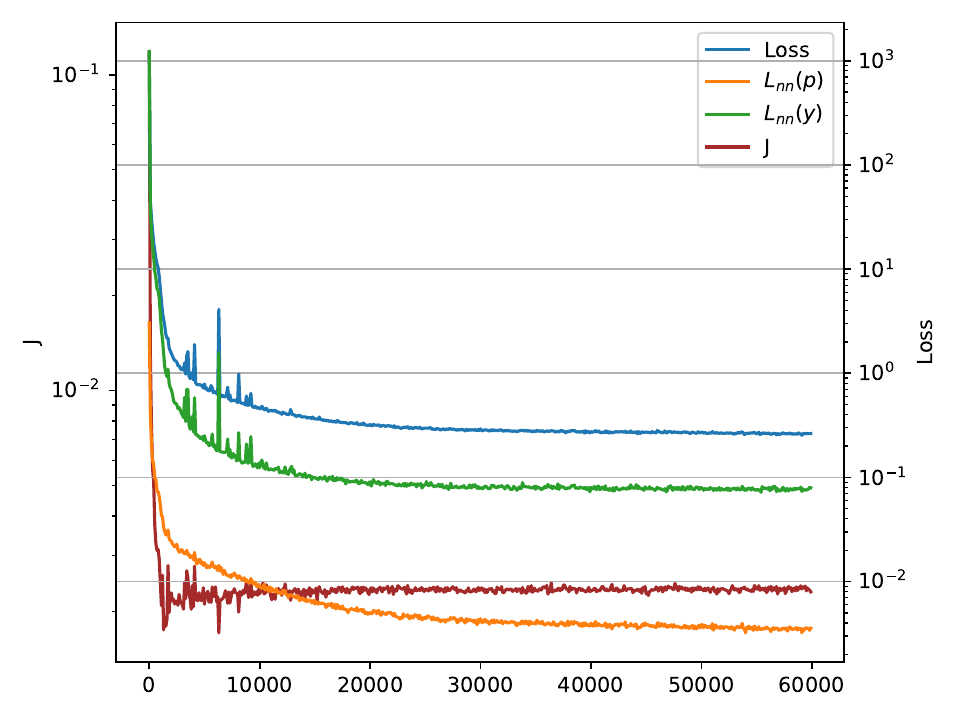}\\
        (a) LU & (b) LC
        & (c) Semi
    \end{tabular}
    \caption{The training dynamics for Example \ref{exam:highdimcpinn}. Ticks on the left y-axis refer to the objective $J$, and the right to  losses $\widehat{\mathcal{L}}_{\rm nn}(y)$ and $\widehat{\mathcal{L}}_{\rm nn}(p)$.}
    \label{fig:hidim-unconst-dyn}
\end{figure}

\section{Conclusion and discussions}
We have investigated a neural solver for distributed elliptic optimal control problems, with / without box constraint on the control. It is based on applying DNNs to a reduced first-order necessary optimality system, and straightforward to implement. We provide an error analysis of the approach, using the approximation theory of neural networks and offset Rademacher complexity, and derived an $L^2(\Omega)$ error bound on the control, state, and adjoint that is explicit in terms of the NN parameters (depth, width), and the number of Monte Carlo sampling points in the domain and on the boundary. We also presented several numerical examples to show its feasibility for linear and semilinear elliptic optimal control problems, and to illustrate its competitiveness with two existing NN-based algorithms. While all three NN based solvers have comparable accuracy, the proposed OSNN tend to be easier to tune and train and involving fewer tunable parameters.

The study leaves many important questions untouched, which we will investigate in future works. First, it would be interesting to provide error analysis for more complex cases, e.g., semilinear elliptic / parabolic optimal control problems. These cases would require new technical tools, e.g., to overcome nonlinearity. Also it is useful to analyze other quantities, e.g., the optimality gap $J(y(u_\sigma),u_\sigma)-J(\bar y,\bar u)$. Second, it is of much interest to evaluate the method on more realistic optimal control settings, including a comparative study traditional discretization methods. Third, it is of enormous importance to devise more accurate and efficient training methods. The Adam optimizer can only yield acceptable but not highly accurate approximations. This has largely limited the numerical experiments to qualitative evaluation. More accurate training algorithms may enable validating the error bound, which however is still unavailable  \cite{SiegelXu:2023}.
\appendix
\section{Technical estimates on $f_\theta$}\label{sec:techestm}

We derive several technical estimates, especially the bound and Lipschitz continuity
of $\partial_{x_p}^2f_\theta$ for $f_{\theta}\in \mathcal{N}_\rho(L,\boldsymbol{n}_L,R)$ in terms of
the NN parameters $\theta$. Let $\boldsymbol n_i$, $i=1,\ldots,L$, be the number of nonzero parameters on
the first $i$ layers, and $\boldsymbol{n}_L$ be the total number of nonzero weights.
First, we recall several bounds on $f_\theta$ and $\partial_{x_p}f_\theta$, which will
be used extensively below. Throughout, we denote $\pi_i= \prod_{j=1}^i n_j$, and let $\mathcal{P}=\mathcal{N}_{\rho}(L,\bs{n}_L,R)$.

\begin{lemma}\label{lem:rho-bdd}
The following estimates hold.
\begin{itemize}
  \item[{\rm(i)}] If $\rho(t)=\tanh(t)$, then
  $\|\rho^{(i)}\|_{L^\infty(\mathbb{R})}\leq1$, $i=0,1,2$, and $\|\rho'''\|_{L^\infty(\mathbb{R})}\leq 2$.
\item[{\rm(ii)}] If $\rho(t)=
\frac{1}{1+e^{-t}}$, then $\|\rho^{(i)}\|_{L^\infty(\Omega)}\leq 1$, $i=0,1,2,3$.
\end{itemize}
\end{lemma}
\begin{proof}
The lemma follows from direct computation. Indeed, for $\rho(t)=\tanh(t)$, we have
$\|\rho\|_{L^\infty(\mathbb{R})}\leq 1$ and
$\rho'(t) 
    =1-\rho(t)^2$,
$\rho''(t)  = -2\rho(t)(1-\rho(t)^2)$,
$\rho'''(t)=(6\rho(t)^2-2)(1-\rho(t)^2)$.
Similarly, for $\rho(t)=\frac{1}{1+e^{-t}}$, clearly $\rho(t)\in (0,1)$, and
$\rho'(t)  = \frac{e^{-t}}{(1+e^{-t})^2} = \rho(t)(1-\rho(t))\in (0,\tfrac12)$,
$  \rho''(t)  =
   \rho(t)(1-\rho(t))(1-2\rho(t))$, 
$\rho'''(t) = \rho(t)(1-\rho(t))[1-6\rho(t)+6\rho(t)^2]$.
Then the desired assertions follow immediately.
\end{proof}

\begin{lemma}\label{lem:est-f}
For any $f\equiv f_\theta,\tilde f\equiv f_{\tilde\theta}\in \mathcal{P}$, the following estimates hold: for any $p\in[d],q\in [n_\ell]$,
\begin{align}
  \|f\|_{L^\infty(\Omega)} & \leq n_{L-1}R,\\
  \|f_q^{(\ell)}-\tilde f_q^{(\ell)}\|_{L^\infty(\Omega)}& \leq \left\{\begin{aligned}
    \pi_{\ell-1}R^{\ell-1}\sum_{j=1}^{\boldsymbol{n}_{\ell}}|\theta_j-\tilde \theta_j|, &\quad \ell =1,\ldots, L-1,\\
    \sqrt{\boldsymbol{n}_L}\pi_{L-1} R^{L-1}\|\theta-\tilde \theta\|_{\ell^2}, &\quad \ell = L,
  \end{aligned}\right.\label{eqn:est-f-Lip}\\
  \|\partial_{x_p}f^{(\ell)}\|_{L^\infty(\Omega)}& \leq
   \pi_{\ell-1}R^\ell, \quad \ell = 1,\ldots,L,\label{eqn:est-f'-bdd}\\
   \|\partial_{ x_p} f_q^{(\ell)}-\partial_{x_p}\tilde f_q^{(\ell)}\|_{L^\infty(\Omega)} & \leq (\ell+1)\pi_{\ell-1}^2R^{2\ell-1}\sum_{j=1}^{\boldsymbol{n}_{\ell}}|\theta_j-\tilde \theta_j|,\quad \ell =1,\ldots,L-1.\label{eqn:est-f'-Lip}
\end{align}
\end{lemma}
\begin{proof}
These estimates are contained in \cite[Lemmas 5.9--5.11]{JiaoLai:2021error}. Note that the estimate
\eqref{eqn:est-f'-Lip} improves that in \cite[Lemma 5.11]{JiaoLai:2021error} by a factor of $R$, by
slightly improving the bound on page 17, line 7 of \cite{JiaoLai:2021error}.
\end{proof}

We also need the following uniform bound on $\partial_{x_p}^2f_\theta$.
\begin{lemma}\label{lem:est-f''-bdd}
Let $\mathcal{P}=\mathcal{N}_{\rho}(L,\bs{n}_{L},R)$. Then for any $p\in [d]$,
$\left| \partial_{x_p}^2 f^{(\ell)}_q\right| \leq \ell \pi^2_{\ell-1}R^{2\ell}$, $\ell = 1,2,\cdots,L$.
\end{lemma}
\begin{proof}
It follows from direct computation that
\begin{align}
  \partial_{x_p}^2f^{(\ell)}_q =& \brb{\sum_{j=1}^{n_{\ell-1}}a_{qj}^{(\ell)}\partial_{x_p}
  f_j^{(\ell-1)}}^2\rho''\brb{\sum_{j=1}^{n_{\ell-1}}a_{qj}^{(\ell)} f_j^{(\ell-1)}+b_q^{(\ell)}}\nonumber\\
    &+\rho'\brb{\sum_{j=1}^{n_{\ell-1}}
  a_{qj}^{(\ell)} f_j^{(\ell-1)}+b_q^{(\ell)}}\sum_{j=1}^{n_{\ell-1}}a_{qj}^{(\ell)} \partial^2_{x_p}f_j^{(\ell-1)}.\label{eqn:f''-exp}
\end{align}
When $\ell=1$, this identity and the uniform bound on $\rho''$ from Lemma \ref{lem:rho-bdd} imply
\begin{align}\label{eqn:est-f''-1}
  \|\partial_{x_p}^2f^{(1)}_q\|_{L^\infty(\Omega)} = \brb{\sum_{j=1}^{n_{0}}a_{qj}^{(1)}}^2\|\rho''\brb{\sum_{j=1}^{n_{0}}a_{qj}^{(1)} x+b_q^{(1)}}\|_{L^\infty(\mathbb{R})}\leq n_0^2R^2.
\end{align}
Next, we treat the case $\ell>1$. The bounds on $\rho'$ and $\rho''$ in Lemma \ref{lem:rho-bdd} and the a priori bound on $a_{qj}^{(\ell)}$ and the estimate \eqref{eqn:est-f'-bdd} imply
\begin{align*}
 \abs{\partial_{x_p}^2f^{(\ell)}_q}  &\leq \brb{\sum_{j=1}^{n_{\ell-1}}\abs{a_{qj}^{(\ell)}}}^2\cdot (\pi_{\ell-2}R^{\ell-1})^2+R\sum_{j=1}^{n_{\ell-1}}\abs{\partial_{x_p}^2 f^{(\ell-1)}_j(x)}\leq \pi_{\ell-1}^2R^{2\ell} +R\sum_{j=1}^{n_{\ell-1}} \abs{\partial_{x_p}^2f_j^{(\ell-1)}}.
\end{align*}
Then applying the recursion repeatedly, the estimate \eqref{eqn:est-f''-1} and mathematical induction yields
$|\partial_{x_p}^2f^{(\ell)}_q| \leq \ell \pi_{\ell-1}^2 R^{2\ell}$, $\ell = 1,\ldots, L-1$.
The case $\ell=L$ follows also from the definition of $f$ and the preceding estimate.
\end{proof}

The next result represents one of the main tools in establishing Rademacher complexity bound.
\begin{lemma}\label{lem:est-f''-Lip}
Let $f_\theta,f_{\tilde\theta}\in \mathcal{P}$, and define $\eta=1$ for sigmoid and $\eta=2$ for hyperbolic tangent.
 Then,
\begin{align*}
       \left|\partial^2_{x_p} f_\theta(x)-\partial^2_{x_p} f_{\tilde\theta}(x)\right| &\leq 2(L-1)L\eta\sqrt{\bs{n}_{L}}\pi_{L-1}^3R^{3L-3}\|{\theta -\tilde{\theta}}\|_{\ell^2}, \quad \forall p\in [d].
\end{align*}
\end{lemma}
\begin{proof}
The proof is based on mathematical induction, and it is divided into three steps.\\
(i) \textbf{Prove the bound for the case $\ell=1$}. By the identity \eqref{eqn:f''-exp}
and triangle inequality, for $\ell=1$,
\begin{align}
&\big|\partial_{x_p}^2 f^{(1)}_q -\partial^2_{x_p} \tilde{f}^{(1)}_q\big| = \Big| (a_{qp}^{(1)})^2\rho''\brb{\sum_{j=1}^{n_0}a_{qj}^{(1)}x_j+b^{(1)}_q}-(\tilde{a}_{qp}^1)^2\rho'' \brb{\sum_{j=1}^{n_0}\tilde{a}_{qj}^{(1)}x_j+\tilde{b}^{(1)}_q}\Big|\nonumber\\
\leq &\big|{(a_{qp}^{(1)})^2-(\tilde{a}_{qp}^{(1)})^2}\big|\Big|{\rho''\brb{\sum_{j=1}^{n_0} a_{qj}^{(1)}x_j+b^{(1)}_q}}\Big|+\big|{(\tilde{a}_{qp}^{(1)})^2}\big|
\Big|{\rho''\brb{\sum_{j=1}^{n_0}a_{qj}^{(1)}x_j+b^{(1)}_q}-\rho''\brb{\sum_{j=1}^{n_0}\tilde{a}_{qj}^{(1)}x_j+\tilde{b}^{(1)}_q}}\Big|\nonumber\\
\leq& 2R\big|{a_{qp}^{(1)}-\tilde{a}_{qp}^{(1)}}\big| + \eta R^2\sum_{j=1}^{n_0}\big|{a_{qj}^{(1)}-\tilde{a}_{qj}^{(1)}}\big|+\eta R^2\big|{b_q^{(1)}-\tilde{b}^{(1)}_q}\big|
\leq 3\eta R^2 \sum_{k=1}^{\bs{n}_1}\big|{\theta_k-\tilde{\theta}_k}\big|,\label{eqn:est-f''-Lip-1}
\end{align}
since by definition, $\eta\geq 1$.\\
(ii) \textbf{Derive the recursion relation.}
For $\ell=2,\ldots,L-1$, in view of the identity \eqref{eqn:f''-exp}, we have
\begin{align*}
\big|\partial^2_{x_p} f^{(\ell)}_q -\partial_{x_p}^2 \tilde{f}^{(\ell)}_q\big|
\leq& \Big|{\brb{\sum_{j=1}^{n_{\ell-1}}a_{qj}^{(\ell)}\partial_{x_p}f^{(\ell-1)}_j}^2- \brb{\sum_{j=1}^{n_{\ell-1}}\tilde{a}_{qj}^{(\ell)}\partial_{x_p}\tilde{f}^{(\ell-1)}_j}^2}\Big| \Big|\rho''\brb{\sum_{j=1}^{n_{\ell-1}}a_{qj}^{(\ell)}f^{(\ell-1)}_j+b^{(\ell)}_q}\Big|\\
&+\brb{\sum_{j=1}^{n_{\ell-1}}\tilde{a}_{qj}^{(\ell)}\partial_{x_p}\tilde{f}^{(\ell-1)}_j}^2
\Big|{\rho''\brb{\sum_{j=1}^{n_{\ell-1}}a_{qj}^{(\ell)}f^{(\ell-1)}_j+b^{(\ell)}_q} -\rho''\brb{\sum_{j=1}^{n_{\ell-1}}\tilde{a}_{qj}^{(\ell)}\tilde{f}^{(\ell-1)}_j+\tilde{b}^{(\ell)}_q}}\Big|\\
&+\Big|{\rho'\brb{\sum_{j=1}^{n_{\ell-1}}a_{qj}^{(\ell)}f^{(\ell-1)}_j+b^{(\ell)}_q}- \rho'\brb{\sum_{j=1}^{n_{\ell-1}}\tilde{a}_{qj}^{(\ell)}\tilde{f}^{(\ell-1)}_j+\tilde{b}^{(\ell)}_q}}\Big| \Big|\brb{\sum_{j=1}^{n_{\ell-1}}a_{qj}^{(\ell)}\partial_{x_p}^2f^{(\ell-1)}_j}\Big|\\
&+ \Big|\rho'\brb{\sum_{j=1}^{n_{\ell-1}}\tilde{a}_{qj}^{(\ell)}\tilde{f}^{(\ell-1)}_j+\tilde{b}^{(\ell)}_q}\Big| \Big|{\sum_{j=1}^{n_{\ell-1}}a_{qj}^{(\ell)}\partial_{x_p}^2f^{(\ell-1)}_j-\sum_{j=1}^{n_{\ell-1}}\tilde{a}_{qj}^{(\ell)}
\partial_{x_p}^2\tilde{f}^{(\ell-1)}_j}\Big| = \sum_{m=1}^4{\rm I}_m.
\end{align*}
Below we bound the four terms separately. To bound the term ${\rm I}_1$, by the estimate \eqref{eqn:est-f'-bdd},
\begin{align*}
  &\sum_{j=1}^{n_{\ell-1}}\brb{\big|{a_{qj}^{(\ell)}}\big|\big|{\partial_{x_p}f^{(\ell-1)}_j}\big|+\big|{\tilde{a}_{qj}^{(\ell)}}\big| \big|{\partial_{x_p}\tilde{f}^{(\ell-1)}_j}\big|}  \leq 2R \sum_{j=1}^{n_{\ell-1}} \pi_{\ell-2} R^{\ell-1} = 2\pi_{\ell-1} R^{\ell}.
\end{align*}
Likewise, by the estimates \eqref{eqn:est-f'-bdd} and \eqref{eqn:est-f'-Lip},
\begin{align*}
  &\sum_{j=1}^{n_{\ell-1}}\big|{a_{qj}^{(\ell)}\partial_{x_p}f^{(\ell-1)}_j-\tilde{a}_{qj}^{(\ell)}\partial_{x_p}\tilde{f}^{(\ell-1)}_j}\big|\\
 \leq & \sum_{j=1}^{n_{\ell-1}}\big|{a_{qj}^{(\ell)}-\tilde{a}_{qj}^{(\ell)}}\big|\big|{\partial_{x_p}f^{(\ell-1)}_j}\big| +\sum_{j=1}^{n_{\ell-1}}\big|{\tilde{a}_{qj}^{(\ell)}}\big|\big|{\partial_{x_p}f^{(\ell-1)}_j-\partial_{ x_p}\tilde{f}^{(\ell-1)}_j}\big|\\
  \leq &\pi_{\ell-2}R^{\ell-1}\sum_{j=1}^{n_{\ell-1}}\big|{a_{qj}^{(\ell)}-\tilde{a}_{qj}^{(\ell)}}\big|+ \sum_{j=1}^{n_{\ell-1}} R\cdot \Big(\ell \pi_{\ell-2}^2R^{2\ell-2} \sum_{k=1}^{\bs{n}_{\ell-1}}\big|{\theta_k-\tilde{\theta}_k}\big|\Big)\\
  \leq &\pi_{\ell-2}R^{\ell-1}\sum_{j=1}^{n_{\ell-1}}\big|{a_{qj}^{(\ell)}-\tilde{a}_{qj}^{(\ell)}}\big|+ \ell\pi_{\ell-2}\pi_{\ell-1} R^{2\ell-1} \sum_{j=1}^{\bs{n}_{\ell-1}}\big|{\theta_j-\tilde{\theta}_j}\big|.
\end{align*}
Thus we can bound the term ${\rm I}_1$ by
\begin{align*}
{\rm I}_1&\leq 2\pi_{\ell-2}\pi_{\ell-1} R^{2\ell-1} \sum_{j=1}^{n_{\ell-1}}\big|{a_{qj}^{(\ell)}-\tilde{a}_{qj}^{(\ell)}}\big|+2\ell \pi_{\ell-2}\pi_{\ell-1}^2R^{3\ell-1} \sum_{j=1}^{\bs{n}_{\ell-1}}\big|{\theta_j-\tilde{\theta}_j}\big|.
\end{align*}
Next, it follows from the estimate \eqref{eqn:est-f-Lip} that
\begin{align*}
  \Big|{\sum_{j=1}^{n_{\ell-1}}a_{qj}^{(\ell)}f^{(\ell-1)}_j-\tilde{a}_{qj}^{(\ell)}\tilde{f}^{(\ell-1)}_j
+b^{(\ell)}_q-\tilde{b}^{(\ell)}_q}\Big|
\leq \sum_{j=1}^{n_{\ell-1}}\big|{a_{qj}^{(\ell)}-\tilde{a}_{qj}^{(\ell)}}\big|
+\big|{b^{(\ell)}_q-\tilde{b}^{(\ell)}_q}\big|+\pi_{\ell-1}R^{\ell-1} \sum_{j=1}^{\bs{n}_{\ell-1}}\big|{\theta_j-\tilde{\theta}_j}\big|.
\end{align*}
Then for the term ${\rm I}_2$, from the estimate \eqref{eqn:est-f'-bdd}, we deduce
\begin{align*}
{\rm I}_2&\leq \eta\Big(\sum_{j=1}^{n_{\ell-1}}\big|{\tilde{a}_{qj}^{(\ell)}}\big|^2\Big)\Big( \sum_{j=1}^{n_{\ell-1}}\big|{\partial_{ x_p}\tilde{f}^{(\ell-1)}_j}\big|^2\Big)\Big(\Big|{\sum_{j=1}^{n_{\ell-1}}a_{qj}^{(\ell)}f^{(\ell-1)}_j-\tilde{a}_{qj}^{(\ell)}\tilde{f}^{(\ell-1)}_j
+b^{(\ell)}_q-\tilde{b}^{(\ell)}_q}\Big|\Big)\\
&\leq \eta \pi_{\ell-1}^2R^{2\ell}
\brb{\sum_{j=1}^{n_{\ell-1}}\big|{a_{qj}^{(\ell)}-\tilde{a}_{qj}^{(\ell)}}\big|
+\big|{b^{(\ell)}_q-\tilde{b}^{(\ell)}_q}\big|+\pi_{\ell-1}R^{\ell-1} \sum_{j=1}^{\bs{n}_{\ell-1}}\big|{\theta_j-\tilde{\theta}_j}\big|}.
\end{align*}
Similarly, Lemma \ref{lem:est-f''-bdd} implies
\begin{align*}
  \sum_{j=1}^{n_{\ell-1}}\big|a_{qj}^{(\ell)}{\partial^2_{x_p}f^{(\ell-1)}_j} \Big| \le R\cdot n_{\ell-1}\cdot (\ell-1)\pi_{\ell-2}^2R^{2(\ell-1)} = (\ell-1)\pi_{\ell-2}\pi_{\ell-1}R^{2\ell-1},
\end{align*}
and thus we can bound the term ${\rm I}_3$ by
\begin{align*}
     {\rm I}_3 
      &\leq (\ell-1)\pi_{\ell-2}\pi_{\ell-1}R^{2\ell-1}
      \brb{\sum_{j=1}^{n_{\ell-1}}\big|{a_{qj}^{(\ell)}-\tilde{a}_{qj}^{(\ell)}}\big|+\big|{b^{(\ell)}_q-\tilde{b}^{(\ell)}_q}\big|
      +\pi_{\ell-1}R^{\ell}\sum_{j=1}^{\bs{n}_{\ell-1}}\big|{\theta_j-\tilde{\theta}_j}\big|}.
\end{align*}
For the last term ${\rm I}_4$, using Lemma \ref{lem:est-f''-bdd} again, we have
    \begin{align*}
    {\rm I}_4&\leq \sum_{j=1}^{n_{\ell-1}}\big|{a_{qj}^{(\ell)}-\tilde{a}_{qj}^{(\ell)}}\big|\big|{\partial^2_{x_p}f^{(\ell-1)}_j}\big| +R\sum_{j=1}^{n_{\ell-1}}\big|{\partial^2_{\partial x_p}f^{(\ell-1)}_j-\partial^2_{x_p}\tilde{f}^{(\ell-1)}_j}\big|\\
    &\leq  (\ell-1)\pi_{\ell-2}\pi_{\ell-1}R^{2(\ell-1)} \sum_{j=1}^{n_{\ell-1}}\big|{a_{qj}^{(\ell)}-\tilde{a}_{qj}^{(\ell)}}\big| +R\sum_{j=1}^{n_{\ell-1}}\big|{\partial^2_{x_p}f^{(\ell-1)}_j-\partial^2_{x_p}\tilde{f}^{(\ell-1)}_j}\big|.
    \end{align*}
Combining the last four estimates gives the crucial recursion
\begin{align}\label{eqn:est-f''-Lip-0}
 \big|\partial_{x_p}^2 f^{(\ell)}_q-\partial_{x_p}^2 \tilde{f}^{(\ell)}_q\big| &\leq R\sum_{j=1}^{n_{\ell-1}}\big|{\partial^2_{x_p}f_j^{(\ell-1)}-\partial^2_{x_p}\tilde{f}_j^{(\ell-1)}}\big|+4\ell\eta\pi_{\ell-1}^3R^{3\ell-1} \sum_{j=1}^{\bs{n}_{\ell}}\big|{\theta_j-\tilde{\theta}_j}\big|.
 \end{align}
(iii) \textbf{Prove the intermediate case by mathematical induction.} Using the recursion \eqref{eqn:est-f''-Lip-0}, we claim that for $\ell=1,2,\ldots,L-1,$ there holds
\begin{align}\label{eqn:est-f''-Lip-claim}
        \big| \partial_{x_p}^2 f^{(\ell)}_q(x) -\partial_{x_p}^2 \tilde{f}^{(\ell)}_q(x) \big| &\leq 2\ell(\ell+1) \eta\pi_{\ell-1}^3R^{3\ell-1} \sum_{j=1}^{\bs{n}_{\ell}}\big|{\theta_j-\tilde{\theta}_j}\big|.
\end{align}
This is trivially true for $\ell=1$, by the estimate \eqref{eqn:est-f''-Lip-1}. For $\ell=2$, the recursion
\eqref{eqn:est-f''-Lip-0} and the estimate for $\ell=1$ in \eqref{eqn:est-f''-Lip-1} imply
    \begin{align*}
        \big| \partial_{x_p}^2 f^{(2)}_q(x) -\partial_{x_p}^2 \tilde{f}^{(2)}_q(x) \big|
        &\leq  3\eta n_1R^3\sum_{j=1}^{\bs{n}_1}\big|{\theta_k-\tilde{\theta}_k}\big|+8\eta R^5 \pi_1^3 \sum_{j=1}^{\bs{n}_{2}}\big|{\theta_j-\tilde{\theta}_j}\big|\\
       & \leq 11\eta\pi_1^3 R^5\sum_{j=1}^{\bs{n}_{2}}\big|{\theta_j-\tilde{\theta}_j}\big| \leq 2\cdot 2 \cdot 3 \eta\pi_1^3 R^5\sum_{j=1}^{\bs{n}_{2}}\big|{\theta_j-\tilde{\theta}_j}\big|,
    \end{align*}
and hence the claim holds for $\ell=2$. Now suppose it holds for some $2\leq \ell<L-1$. Then for the case $\ell+1$,
    \begin{align*}
        \left| \partial_{x_p}^2f^{(\ell+1)}_q(x) -\partial^2_{x_p}\tilde{f}^{(\ell+1)}_q(x) \right| &\leq R\sum_{j=1}^{n_{\ell}}\abs{\partial^2_{x_p}f_j^{(\ell)}-\partial^2_{x_p}\tilde{f}_j^{(\ell)}}+ 4(\ell+1) \eta \pi_\ell^3R^{3\ell+2}  \sum_{j=1}^{\bs{n}_{\ell+1}}\abs{\theta_j-\tilde{\theta}_j}\\
        &\leq n_{\ell}\cdot 2\ell(\ell+1)\eta \pi_{\ell-1}^3R^{3\ell-1} \sum_{j=1}^{\bs{n}_{\ell}}\abs{\theta_j-\tilde{\theta}_j}+4(\ell+1)\eta \pi_{\ell}^3R^{3\ell+2} \sum_{j=1}^{\bs{n}_{\ell+1}}\abs{\theta_j-\tilde{\theta}_j}\\
        &\leq 2(\ell+1)(\ell+2) \pi_\ell^3R^{3\ell+2} \sum_{j=1}^{\bs{n}_{\ell+1}}|\theta_j-\tilde{\theta}_j|,
    \end{align*}
which by mathematical induction implies the desired claim \eqref{eqn:est-f''-Lip-claim}.\\
(iv) \textbf{Obtain the final estimate}. Last, for $\ell =L$, direct computation shows
\begin{align*}
   \partial^2_{x_p} f -\partial^2_{x_p}\tilde{f} = \sum_{j=1}^{n_{L-1}}a^{(L)}_{1j}\partial_{x_p}^2 f^{(L-1)}_j - \sum_{j=1}^{n_{L-1}}\tilde a^{(L)}_{1j}\partial_{x_p}^2 \tilde f^{(L-1)}_j
\end{align*}
Thus, by the estimate \eqref{eqn:est-f''-Lip-claim} and Lemma \ref{lem:est-f''-bdd},
\begin{align*}
  \big|\partial^2_{x_p} f -\partial^2_{x_p}\tilde{f}\big| &\leq \sum_{j=1}^{n_{L-1}}|a^{(L)}_{1j}-\tilde a_{1j}^{(L)}||\partial_{x_p}^2 f^{(L-1)}_j| + \sum_{j=1}^{n_{L-1}}|\tilde a^{(L)}_{1j}|\partial_{x_p}^2 f^{(L-1)}_j-\partial_{x_p}^2 \tilde f^{(L-1)}_j|\\
   &\leq (L-1) \pi_{L-2}R^{2L-2}  \sum_{j=1}^{n_{L-1}}|a^{(L)}_{1j}-\tilde a_{1j}^{(L)}| + \sum_{j=1}^{n_{L-1}}R\cdot 2(L-1)L\eta\pi_{L-2}^3R^{3L-4} \sum_{j=1}^{\bs{n}_{L-1}}\big|{\theta_j-\tilde{\theta}_j}\big|\\
   & \leq 2(L-1)L\eta\pi_{L-1}^3R^{3L-3} \sum_{j=1}^{\bs{n}_{L}}\big|{\theta_j-\tilde{\theta}_j}\big|.
\end{align*}
This completes the proof of the lemma.
\end{proof}

\bibliographystyle{abbrv}
\bibliography{nn}

\end{document}